\documentclass[12pt,oneside,english]{amsart}
\usepackage[T1]{fontenc}
\usepackage[latin1]{inputenc}
\usepackage{graphicx}
\usepackage{amssymb}

\makeatletter

\newcommand{\noun}[1]{\textsc{#1}}

 \theoremstyle{plain}    
 \newtheorem{thm}{Theorem}[section]
 \numberwithin{equation}{section} 
 \numberwithin{figure}{section} 
 \theoremstyle{plain}
 \theoremstyle{definition}
 \newtheorem{defn}[thm]{Definition}
 \theoremstyle{remark}
 \newtheorem{rem}[thm]{Remark}
 \theoremstyle{definition}
  \newtheorem{example}[thm]{Example}
 \theoremstyle{plain}    
 \newtheorem{prop}[thm]{Proposition} 
 \theoremstyle{plain}    
 \newtheorem{conjecture}[thm]{Conjecture} 
 \theoremstyle{plain}    
 \newtheorem{lem}[thm]{Lemma} 
 \theoremstyle{plain}    
 \newtheorem{cor}[thm]{Corollary} 

\usepackage{amsmath}
\usepackage{amsfonts}
\usepackage{amssymb}
\usepackage{babel}
\input{diagrams}
\input xy
\xyoption{all}
\def\ms{\mathcal{S}}
\def\bb{\bullet}
\def\FF{\mathbb{F}}
\def\CC{\mathbb{C}}
\def\RR{\mathbb{R}}
\def\QQ{\mathbb{Q}}
\def\II{\mathbb{I}}

\def\ZZ{\mathbb{Z}}
\def\PP{\mathbb{P}}

\def\HH{\mathbb{H}}
\def\EE{\mathbb{E}}
\def\DD{\mathbb{D}}
\def\AA{\mathbb{A}}
\def\GG{\mathbb{G}}
\def\NN{\mathbb{N}}
\def\VV{\mathbb{V}}
\def\WW{\mathbb{W}}
\def\LL{\mathbb{L}}
\def\C{\mathcal{C}}
\def\A{\mathcal{A}}
\def\B{\mathcal{B}}
\def\D{\mathcal{D}}

\def\E{\mathcal{E}}
\def\F{\mathcal{F}}
\def\G{\mathcal{G}}
\def\H{\mathcal{H}}
\def\I{\mathcal{I}}
\def\J{\mathcal{J}}
\def\K{\mathcal{K}}
\def\L{\mathcal{L}}
\def\M{\mathcal{M}}

\def\N{\mathcal{N}}
\def\O{\mathcal{O}}

\def\R{\mathcal{R}}
\def\T{\mathcal{T}}

\def\V{\mathcal{V}}
\def\W{\mathcal{W}}
\def\TW{\widetilde{\W}}
\def\X{\mathcal{X}}
\def\Y{\mathcal{Y}}
\def\Z{\mathcal{Z}}

\def\a{\alpha}

\def\d{\partial}

\def\e{\epsilon}

\def\o{\Omega}

\def\sV{\mathfrak{V}}
\def\sW{\mathfrak{W}}

\def\sZ{\mathfrak{Z}}
\def\gs{\eta_{\ms}}

\def\rateq{\buildrel\text{rat}\over\equiv}

\def\homeq{\buildrel\text{hom}\over\equiv}
\def\nrateq{\mspace{7mu}/\mspace{-17mu}\rateq}
\def\nequiv{\mspace{7mu}/\mspace{-17mu}\equiv}

\def\sumint{\mspace{7mu}\int\mspace{-25mu}\sum}
\def\sumintt{\mspace{7mu}\int\mspace{-16mu}\Sigma}

\def\bb{\bullet}

\def\hm{Hom_{_{\text{MHS}}}}
\def\ext{Ext^1_{_{\text{MHS}}}}
\def\MHS{\text{MHS}}
\def\PMHS{\text{PMHS}}
\def\MHM{\text{MHM}}
\def\VMHS{\text{VMHS}}

\def\dlog{\text{dlog}}

\def\db{\partial_{\mathcal{B}}}
\def\m{\setminus}
\def\bx{\square}
\def\bxn{\bx^n}
\def\bxnx{\bx^n_X}
\def\fs{f^{\sharp}}
\def\tot{\mathbf{s}}

\usepackage{babel}
\makeatother
\begin{document}

\title{The Abel-Jacobi map for \\ higher Chow groups, \\II}

\author{Matt Kerr and James D. Lewis}

\maketitle
\tableofcontents{}

\section{\textbf{Introduction}}

This paper has its origin in two distinct, well-established problems:
that of detecting $K_{n}^{^{\text{{alg}}}}(X)$ classes --- especially
indecomposable ones --- in the regulator kernel for $X$ smooth projective;
and that of computing the limit of the regulator for a family of such
classes as $X$ degenerates. It has become a sequel to \cite{KLM}
because these problems benefit from two {}``dual'' extensions of
the regulator formulas introduced there.

To explain this, we recall the relation between Bloch's higher Chow
groups, $K$-theory and motivic cohomology \begin{equation} CH^p(X,n)_{_{\QQ}} \cong K_n^{^{\text{alg}}}(X)^{^{(p)}}_{_{\QQ}} \cong H_{\M}^{2p-n}(X,\QQ(p)) \end{equation} 
for $X$ smooth projective. For our purposes, the regulator is the
absolute Hodge cycle-class map \begin{equation} c_{_{\H}}^{_X} : \, CH^p(X,n)_{_{\QQ}} \rTo H^{2p-n}_{\H}(X,\QQ(p)),  \end{equation} 
and \cite{KLM} introduced an explicit morphism of complexes, from
cycles to {}``triples'' of currents, computing this. Nontrivial
cycles $\Z\in\ker(c_{_{\H}}^{_{X}})$ can be studied by \emph{spreading
out} ($X,\Z\mapsto\X_{U},\sZ_{U}$ as in $\S4.1$) and computing $c_{_{\H}}^{_{\X_{U}}}(\sZ_{U})$.
Here $\X_{U}$ is quasi-projective but still smooth, and so (1.1),
(1.2) still make sense; however, \cite{KLM} does not give a map of
complexes computing $c_{_{\H}}$ for quasi-projective varieties.

This is remedied in $\S3$, which as a byproduct also yields a formula
for\[
c_{_{\H,Y}}:\, CH^{p}(Y,n)\rTo H_{\H,Y}^{2p-n-2}(X,\QQ(p+1))\]
when $Y\subset X$ is a complete normal crossing divisor (NCD). But
this is the {}``wrong'' regulator for singular varieties, as far
as many arithmetic applications or (finite) limits under degeneration
are concerned. Indeed, for \emph{general} $Y$ \[
CH^{p}(Y,n)\cong G_{n}(Y)_{_{\QQ}}^{^{(p)}}\cong H_{2d_{Y}-2p+n}^{\M,B.M.}(Y,\QQ(d_{Y}-p))\]
is motivic Borel-Moore \emph{homology}; whereas \begin{equation} c^{_Y}_{_{\H}} : \, H_{\M}^{2p-n}(Y,\QQ(p)) \rTo H^{2p-n}_{\H}(Y,\QQ(p)) \end{equation} 
is the interesting map. Again, a morphism of complexes computing this
--- for many types of (complete) singular and relative varieties ---
is the subject of $\S\S8.1-3$. The construction involves partially
dualizing that of $\S3$, as well as considerably greater technical
difficulties.

Applications of the second construction get relatively short shrift
in this paper, but $\S8.4$ does compute (1.3) for higher cycles on
a degenerate elliptic curve and $K3$ surface --- examples related
to irrationality of $\zeta(2)$ and $\zeta(3)$, respectively, see
\cite{Ke5}. (The cycles actually do deform in each case as the variety
{}``smooths''.) The formulas of $\S8$ are also used in \cite{DK}.

On the other hand, almost half the paper ($\S\S4-7$) concerns applications
of $c_{_{\H}}$ for smooth quasi-projective varieties, to define and
compute \emph{higher Abel-Jacobi invariants} for higher Chow cycles%
\footnote{to be clear: \cite{KLM} studied usual $AJ$ on higher cycles; \cite{GG,Le1,Ke5}
studied higher $AJ$ on usual (algebraic) cycles; here we study higher
$AJ$ on higher cycles, which includes both.%
} (on projective $X$) in the kernel of (1.2). These are, roughly speaking,
graded pieces of $c_{_{\H}}^{_{\X_{U}}}(\sZ_{U})$ under the Leray
filtration for a smooth morphism $\X_{U}\to U$ (in the limit as we
shrink $U$). These higher invariants appear in various guises: in
$\S4$, as a generalization of the classical notion of normal functions;
in $\S5$, as K\"unneth components of the \cite{KLM} map; and in
$\S7$, as generalized membrane integrals on $X$ in the spirit of
\cite{GG} and \cite[sec. 11]{Ke1}. 

The actual double-complex formalism used in the formula for $c_{_{\H}}$
in $\S3$, only shows up once in the higher $AJ$ discussion --- in
$\S5.2$ to prove a very precise result on nontrivial \emph{decomposable}
cycles in the regulator kernel (Theorem 5.1). These are detected by
composing the higher $AJ$'s with \emph{higher residue maps}, a feature
not present when spreading out usual algebraic cycles. But the most
beautiful results on (higher) cycles are to be found in $\S6$, where
we prove that certain regulator-trivial exterior products (Theorem
6.2) and Pontryagin products (Theorem 6.3) of higher Chow cycles are
strongly indecomposable, i.e not in the image of\[
CH^{p-1}(X_{K},n-1)\otimes K^{*}\rTo CH^{p}(X_{K},n).\]
The {}``indecomposable higher $AJ$'' invariant used to detect these
(Theorem 6.1) is an improvement of an earlier invariant of the second
author \cite{Le2}, and Theorem 6.3 is likewise a strengthening of
a result of Collino and Fakhruddin \cite{CF} which used \cite{Le2}.

Section $2$ is partly expository: we show how to compute $H_{\H}^{q}(\W,\QQ(p)$
($\W$ smooth quasi-projective) as the cohomology of a double complex,
then place a {}``weight'' filtration on it and compute the graded
pieces (Prop. 2.7). This is all just to clear the runway for $\S3$.

We want this paper (in conjunction with \cite{KLM}) to be useful
as a handbook of sorts on various cycle-class maps, so we have attempted
a {}``modular'' format as far as possible. With the following exceptions,
a well-informed reader can probably start at the beginning of any
section: $\S3$ depends on $\S2$, $\S5.2$ on $\S3.1$, and $\S6$
and $\S7$ each on $\S5.1$. Our definition of {}``well-informed''
includes familiarity with (or at least reference to) $\S\S5.1-5.7$
of \cite{KLM}, which is used almost everywhere.

\subsection{Notation}

Here we forewarn the reader of some potentially confusing idiosyncracies.
Most importantly: \textbf{cycle groups throughout this paper are automatically
$\mathbf{\otimes\mathbb{{Q}}}$}, in particular the higher precycles
$Z^{p}(X,m)$, $Z_{\RR}^{p}(X,m)$ (see $\S8.2$ for a defnition of
the latter). These are not to be confused with the {}``topological''
$C^{\infty}$ chains resp. cycles with coefficients in $(2\pi i)^{q}\QQ$,
written $C_{top}^{r}(X_{\CC}^{an},\QQ(q))$ resp. $Z_{top}^{r}(X_{\CC}^{an},\QQ(q))$.
Here $r$ refers to $real$ codimension (in contrast to $p$). The
sheaf of $C^{\infty}$ forms of degree $q$ is denoted $\Omega_{X^{\infty}}^{q}$;
currents $\D_{X}^{2\dim(X)-q}$ are the distributions that act on
them.

A higher Chow precycle $\Z\in Z^{p}(X,m)_{_{\db\text{{-closed}}}}$
has rational equivalence class $\left\langle \Z\right\rangle \in CH^{p}(X,m)$
and homological equivalence class $\left[\Z\right]\in\hm(\QQ(0),H^{2p-n}(X,\QQ(p))$.
When it is not crucial, we are sometimes lax about the $\Z\leftrightarrow\left\langle \Z\right\rangle $
distinction.

Given a Hodge structure $\H$, $F_{h}^{q}\H$ is defined to be the
largest sub-Hodge structure of $\H$ contained in $\H\cap F^{q}\H_{\CC}$.
If $\W$ is smooth quasi-projective with good compactification $\overline{\W}$,
then for our various cycle/cohomology groups, an underline denotes
the image of the restriction from $\overline{\W}$ to $\W$: e.g.,
$\underline{CH^{p}}(\W,m):=\text{{im}}\{ CH^{p}(\overline{\W},m)\to CH^{p}(\W,m)\}.$

The following is more standard: for a complex $(\K^{\bb},d)$, the
truncation $\tau_{\bb\leq q}\K^{\bb}$ is $\K^{\bb}$ for $\bb<q$,
0 for $\bb>q$, and $\{\ker(d)\subset\K^{q}\}$ for $\bb=q$. $D^{b}(\text{{---}})$
is the bounded derived category, and $\tot^{\bb}(\text{{---}})$ is
the simple complex associated to a double complex.

\subsection{Acknowledgements}

This paper was conceived at the May 2005 Antalya Algebra Days in Antalya,
Turkey; both authors would like to thank the organizers for a stimulating
atmosphere. We would also like to thank Donu Arapura for sharing \cite{Ar2},
and Shuji Saito for helpful conversations.

\section{\textbf{Absolute Hodge cohomology}}

Let $X$ be a smooth projective analytic variety $/\CC$, $Y=Y_{1}\cup\cdots\cup Y_{N}\subset X$
a normal-crossing divisor (NCD) with smooth components (and smooth
intersections of all orders), and $\AA\subseteq\RR$ a subring such
that $\AA\otimes_{\ZZ}\QQ$ is a field. The principal aim of this
section is to produce explicit complexes whose cohomologies compute
the absolute $\AA$-Hodge cohomology groups\[
H_{\H}^{q}(X,\AA(p)),\,\, H_{\H}^{q}(X\m Y,\AA(p)),\,\, H_{\H,Y}^{q}(X,\AA(p))\]
as defined in \cite[sec. 2]{Ja1}. At least half of the section is
expository, albeit brief and sketchy.

\begin{defn}
\cite{Ja1},\cite{BZ} A mixed $\AA$-Hodge complex consists of the
following data:
\end{defn}
\begin{itemize}
\item A bounded below complex $K_{\AA}^{\bb}$ of $\AA$-modules (in the
derived category), such that $H^{p}(K_{\AA}^{\bb})$ is an $\AA$-module
of finite type for all $p$.
\item A bounded below filtered complex $(K_{\AA\otimes\QQ}^{\bb},W)$ of
$\AA\otimes_{\ZZ}\QQ$-vector spaces, and an isomorphism $K_{\AA\otimes\QQ}^{\bb}\rTo^{\sim}K_{\AA}^{\bb}\otimes\QQ$
in the derived category.
\item A bifiltered complex $(K_{\CC}^{\bb},W,F)$ of $\CC$-vector spaces,
and a filtered isomorphism $\alpha:\,(K_{\CC}^{\bb},W)\rTo^{\sim}(K_{\AA\otimes\QQ}^{\bb},W)\otimes\CC$
in the filtered derived category.
\end{itemize}
$\mspace{100mu}$ Further,

\begin{itemize}
\item For every $m\in\ZZ$,\[
Gr_{m}^{W}K_{\AA\otimes\QQ}^{\bb}\to(Gr_{m}^{W}K_{\CC}^{\bb},F)\]
is a (polarizable) $\AA\otimes\QQ$-Hodge complex of weight $m$,
i.e. the differentials of $Gr_{m}^{W}K_{\CC}^{\bb}$ are strictly
compatible with the induced filtration $F$, and $F$ induces a pure
(polarizable) $\AA\otimes\QQ$-Hodge structure of weight $m+r$ on
$H^{r}(Gr_{m}^{W}K_{\AA\otimes\QQ}^{\bb})$ for $r\in\ZZ$.
\end{itemize}
${}$\\
A mixed $\AA$-Hodge complex gives rise to a diagram: \begin{equation} \xymatrix{ & 'K_{\AA\otimes \QQ}^{\bb}  & & ('K_{\CC}^{\bb},W) \\ K_{\AA}^{\bb} \ar [ru]^{\alpha_1} & & (K_{\AA\otimes \QQ}^{\bb},W) \ar [lu]_{\alpha_2} \ar [ru]^{\beta_1} & & (K_{\CC}^{\bb},W,F), \ar [lu]_{\beta_2} \\ }  \end{equation}
where $\alpha_{j},\,\beta_{j}$, $j=1,2$, are morphisms of complexes,
$\alpha_{2}$ being a quasi-isomorphism, $\beta_{1}$ a filtered morphism,
and $\beta_{2}$ a filtered quasi-isomorphism. By the work of Deligne
and Beilinson, the construction of mixed $\AA$-Hodge complexes is
equivalent to the construction of mixed $\AA$-Hodge structures.

Let $\mu:\,(M^{\bb},d_{M})\to(N^{\bb},d_{N})$ be a morphism of complexes.
The cone of $\mu$ is the complex:\[
Cone\{ M^{\bb}\rTo^{\mu}N^{\bb}\}=C_{\mu}^{\bb}:=M[1]^{\bb}\oplus N^{\bb},\]
with differential $D(a,b)=(-d_{M}(a),\,\mu(a)+d_{N}(b))$. The absolute
Hodge cohomology $H_{\H}^{\ell}(K^{\bb})$ (of a mixed $\AA$-Hodge
complex $K^{\bb}$) is given by $H^{\ell}$ of \begin{equation} Cone \{ K_{\AA}^{\bb} \oplus \hat{W}_0 K_{\AA \otimes \QQ}^{\bb} \oplus (\hat{W}_0 \cap F^0)K_{\CC}^{\bb} \rTo^{(\alpha,\beta )} {}'K_{\AA \otimes \QQ}^{\bb} \oplus \hat{W}_0 {} 'K_{\CC}^{\bb} \} [-1], \end{equation} 
where\[
(\alpha,\beta)\left(\xi_{\AA},\xi_{\QQ},\xi_{\CC}\right)=\left(\alpha_{1}\xi_{\AA}-\alpha_{2}\xi_{\QQ},\,\beta_{1}\xi_{\QQ}-\beta_{2}\xi_{\CC}\right),\]
and $\hat{W}_{\bb}$ is the filtration decal\'ee \cite[sec. 2]{De2}.
More precisely, set $'W^{\bb}:=W_{-\bb}$, $\hat{W}_{\bb}:=(\text{{Dec}}\,'W)^{-\bb}$,
where \[
(\text{{Dec}}\,'W)^{\ell}K^{m}:=\ker\left\{ 'W^{\ell+m}K^{m}\rTo^{d_{K}}\frac{K^{m+1}}{'W^{\ell+m+1}K^{m+1}}\right\} .\]

We henceforth assume $K_{\AA}^{\bb}={}'K_{\AA}^{\bb}$, $K_{\CC}^{\bb}={}'K_{\CC}^{\bb}$
with $\alpha_{2},\beta_{2}$ identity maps, and $K_{\AA\otimes\QQ}^{\bb}=K_{\AA}^{\bb}\otimes_{\ZZ}\QQ$
with $\alpha_{1}$ the inclusion.

\begin{rem}
Suppose also $\AA=\QQ$, so $\alpha_{1}=\,$identity. In this case,
there is a quasi-isomorphism\\
\xymatrix{&  (*) \; \; \;\;\;\; Cone \{ \hat{W}_0 K_{\QQ}^{\bb} \oplus F^0\hat{W}_0 K_{\CC}^{\bb} \rTo^{\beta_1 -\beta_2} \hat{W}_0 K_{\CC}^{\bb} \} [-1] \;\;\;\;\;\;\;\; \ar [d]^{\simeq} \\ &   Cone \{ K_{\QQ}^{\bb} \oplus \hat{W}_0 K_{\QQ}^{\bb} \oplus F^0 \hat{W}_0 K_{\CC}^{\bb} \rTo^{(\alpha,\beta)} K_{\QQ}^{\bb} \oplus \hat{W}_0 K_{\CC}^{\bb} \} [-1] }\\
given by $(\xi_{1},\xi_{2},\xi_{3})\mapsto(\xi_{1},\xi_{1},\xi_{2},0,\xi_{3}).$
(The differential in $(*)$ is $D(\xi_{1},\xi_{2},\xi_{3})=(-d\xi_{1},-d\xi_{2},d\xi_{3}+\beta_{1}\xi_{1}-\beta_{2}\xi_{2})$.)
We refer to $(*)$ as the normalized AH complex. It has an obvious
(non-quasi-isomorphic) embedding in the related {}``Deligne'' (or
weak AH) complex \[
Cone\left\{ K_{\QQ}^{\bb}\oplus F^{0}K_{\CC}^{\bb}\rTo^{\beta_{1}-\beta_{2}}K_{\CC}^{\bb}\right\} [-1],\]
 which simply forgets the weights.
\end{rem}
We do the basic example for general $\AA$, and then restrict to the
case of Remark 2.2; though we do explain (Remark 2.6) how to extend
the result of the Proposition to arbitrary $\AA$.

\begin{example}
(Absolute $\AA$-Hodge cohomology of $X$.) Let \[
('K_{\AA}^{\bb}=)K_{\AA}^{\bb}:=C(X,\AA(p))[2p]^{\bb}=C^{2p+\bb}(X,\AA(p)),\]
\[
K_{\AA\otimes\QQ}^{\bb}:=C^{2p+\bb}(X,\AA(p))\otimes_{\ZZ}\QQ,\]
\[
('K_{\CC}^{\bb}=)K_{\CC}^{\bb}:=\D(X)(p)[2p]^{\bb},\]
with $K_{\AA\otimes\QQ}^{\bb}\rTo^{\beta_{1}}K_{\CC}^{\bb}$ given
by integration: $(2\pi i)^{p}\gamma\mapsto(2\pi i)^{p}\delta_{\gamma}.$
$F^{\bb}$ is the Hodge filtration (on $\D_{X}^{\bb}$, twisted by
$p$); and $W_{r}K_{\cdots}^{\bb}:=\left\{ \begin{array}{cc}
K_{\cdots}^{\bb} & r\geq0\\
0 & r<0\end{array}\right.$ (the stupid weight filtration) yields e.g.\small \[
\hat{W}_{0}K_{\AA\otimes\QQ}^{m}=\ker\left\{ W_{-m}K_{\AA\otimes\QQ}^{m}\to\frac{K_{\AA\otimes\QQ}^{m+1}}{W_{-m-1}K_{\AA\otimes\QQ}^{m+1}}\right\} =\left\{ \begin{array}{cc}
K_{\AA\otimes\QQ}^{m}, & m<0\\
\ker(\d)\subset K_{\AA\otimes\QQ}^{0}, & m=0\\
0, & m>0\end{array}\right..\]
\normalsize Writing out (2.2) explicitly gives a complex which is\[
C^{2p+\bb}(X,\AA(p))\,\oplus\, C^{2p+\bb-1}(X,\AA(p))\otimes_{\ZZ}\QQ\]
for $\bb>1$,\[
C^{2p+1}(X,\AA(p))\,\oplus\, C^{2p}(X,\AA(p))\otimes_{\ZZ}\QQ\,\oplus\,\{\ker(d)\subset\D^{2p}(X)\}\]
for $\bb=1$,\[
C^{2p}(X,\AA(p))\,\oplus\,\{\ker(\d)\subset C^{2p}(X,\AA(p))\otimes_{\ZZ}\QQ\}\,\oplus\,\{\ker(d)\subset F^{p}\D^{2p}(X)\}\]
\[
\oplus C^{2p-1}(X,\AA(p))\otimes_{\ZZ}\QQ\,\oplus\,\D^{2p-1}(X)\]
for $\bb=0$, and\[
C^{2p+\bb}(X,\AA(p))\,\oplus\, C^{2p+\bb}(X,\AA(p))\otimes_{\ZZ}\QQ\,\oplus\, F^{p}\D^{2p+\bb}(X)\]
\[
\oplus\, C^{2p+\bb-1}(X,\AA(p))\otimes_{\ZZ}\QQ\,\oplus\,\D^{2p+\bb-1}(X)\]
for $\bb<0$. For $q\in\ZZ$, its $q^{\text{{th}}}$ cohomology is
$H_{\H}^{2p+q}(X,\AA(p))$. This maps to $H_{\D}^{2p+q}(X,\AA(p))$
by forgetting all $\hat{W}_{0}$'s in the definition of the $Cone$
complex; the map is an isomorphism (only) in degrees $q\leq0$.
\end{example}
Now for each $\ell\in\NN$, let $Y^{\ell}$ denote the union of the
$\ell$-fold intersections of various components of $Y$, so that
$\widetilde{Y^{\ell}}$ is their disjoint union. Set $Y^{0}=X$ and
$Y^{\ell}=\emptyset$ for $\ell<0$. The proper hypercovering $\widetilde{Y^{\bb}}\to Y$
(together with the inclusion $Y\hookrightarrow X$) leads to double
complexes\[
\D^{2p+2i+j}(\widetilde{Y^{-i}})(p+i)=:DR(p)^{i,j}\,,\,\,\,\, C^{2p+2i+j}(\widetilde{Y^{-i}},\QQ(p+i))=:\B(p)^{i,j}\]
with differentials $Gy$ ($i\mapsto i+1$, see formula in $\S3$)
and $d$ resp. $\d$ ($j\mapsto j+1$). (Here $\B$, $DR$ stand for
{}``Betti'', {}``de Rham''.) We have associated simple complexes
with total differential $\DD$, {}``weight'' filtrations ($\K:=\B$
or $DR$)\[
\tot^{k}\K(p):=\oplus_{i(\leq0)}\K(p)^{i,\, k-i}\supset\oplus_{\ell\leq i(\leq0)}\K(p)^{i,\, k-i}=:('W)^{\ell}\tot^{k}\K(p),\]
 and Hodge filtration\[
F^{q}\tot^{k}DR(p):=\oplus_{i}F^{q}\D^{2p+2i+j}(\widetilde{Y^{-i}})(p+i)\cong\oplus_{i}F^{q+p+i}\D^{2p+2i+j}(\widetilde{Y^{-i}}).\]

The function of these complexes is to compute cohomology of $X\m Y$:\[
H^{*}(F^{q}\tot^{\bb}DR(p))\cong F^{q+p}H^{*+2p}(X\m Y,\CC),\,\,\, H^{*}(\tot^{\bb}\B(p))\cong H^{*+2p}(X\m Y,\QQ(p)).\]
The only nontrivial observation involved in proving this is that \begin{equation} \frac{\D^{\bb}(X)}{\sum_j \iota^{Y_j}_{_{*}} \D^{\bb}(Y_j)[-2]} \to \frac{\D^{\bb}(X)}{\D^{\bb}(\text{on} Y)(X)} \cong \D^{\bb}(\log Y)(X) \end{equation} 
is an $F^{\bb}$-filtered quasi-isomorphism \cite{Ja1}, so that $H^{*}$
of {[}$F^{p}$ of{]} the l.h.s. computes $[F^{p}]H^{*}(X\m Y,\CC)$.
Then one looks at the spectral sequence $\{ E_{0}^{i,j}=DR(p)^{i,j},\, d_{0}=Gy\}$
computing $H^{*}(\tot^{\bb}DR(p))$; its $E_{1}^{i,j}=0$ for $i\neq0$,
and $E_{1}^{0,\bb}$ is the l.h.s. of $(2.3)$ shifted by $[2p]$.

Now we explain briefly the function of the decal\'ee filtration \[
(\text{{Dec}}('W))^{\ell}\tot^{k}DR(p)=\ker\left\{ 'W^{\ell+k}\tot^{k}DR(p)\rTo^{\DD}\frac{\tot^{k+1}DR(p)}{'W^{\ell+k+1}\tot^{k+1}DR(p)}\right\} \]
\small \[
=\left\{ \ker(d)\subset\D^{2p+2k+\ell}(\widetilde{Y^{-\ell-k}})(p+k+\ell)\right\} \,\bigoplus\,\oplus_{\ell+k<i(\leq0)}\D^{2p+i+k}(\widetilde{Y^{-i}})(p+i)\]
\normalsize \[
=\tot^{k}\left\{ \tau_{j\leq-\ell}DR(p)^{i,j}\right\} .\]
The usual weight filtration on cohomology is defined by \begin{equation} \begin{matrix} W_{2p+k+\ell}H^{2p+k}(X\m Y,\CC) \; := \; 'W^{-\ell}H^k(\tot^{\bb}DR(p)) \\ \\ := \; \text{im} \{ H^k( 'W^{-\ell}\tot^{\bb} DR(p)) \to H^k(\tot^{\bb} DR(p)) \} . \end{matrix} \end{equation} 
To relate $'W$ and $\text{{Dec}}('W)$, consider their associated
spectral sequences (both $\implies$ $H^{*}(\tot^{\bb}DR(p))$), with
$0^{\text{{th}}}$ pages\[
\E_{0}^{i,j}:=Gr_{'W}^{i}\tot^{i+j}DR(p)\cong DR(p)^{i,j},\,\,\,\,\,\hat{\E}_{0}^{i',j'}:=Gr_{\text{{Dec}}('W)}^{i'}\tot^{i'+j'}DR(p),\]
and $d_{0}=Gr\,\DD$. The first is just the Gysin spectral sequence,
which by \cite{De1} degenerates at $\E_{2}$; hence the second spectral
sequence degenerates at $\hat{\E}_{1}$, since by \cite[sec. 2]{De2}
$\hat{\E}_{r}^{-j,\, i+2j}\rTo^{\cong}\E_{r+1}^{i,j}$ $\forall\, r\geq1$.
So $\DD$ on $\tot^{\bb}DR(p)$ is strictly compatible with $\text{{Dec}}('W)^{\bb}$,
and \begin{equation} \begin{matrix} H^k \left( \text{Dec}('W)^{-\ell -k}\tot^{\bb} DR(p) \right) \cong \\ \\ \text{im} \{ H^k \left( \text{Dec}('W)^{-\ell-k} \tot^{\bb} DR(p) \right) \to H^k(\tot^{\bb} DR(p)) \} . \end{matrix} \end{equation} 
The above compatibility of spectral sequences also implies agreement
of the last lines of (2.4) and (2.5). Writing $\hat{W}_{0}=\text{{Dec}}('W)^{0}$
and specializing to $k=-\ell$, we have shown that \begin{equation} H^k \left( \hat{W}_0 \tot^{\bb} DR(p) \right) \cong W_{2p}H^{2p+k}(X\m Y,\CC) . \end{equation}
This argument can be applied also to $\tot^{\bb}\B(p)$ and $F^{*}\tot^{\bb}DR(p)$.

To define absolute Hodge cohomology, we take our mixed $\QQ$-Hodge
complex $K^{\bb}$ to be given by\[
K_{\QQ}^{\bb}:=\tot^{\bb}\B(p)\,\,\,,\,\,\,\,\,\,\, K_{\CC}^{\bb}:=\tot^{\bb}DR(p),\]
with $W_{\bb}(='W^{-\bb})$ and $F^{\bb}$ as above, and $\beta_{1}$
given by integration (on each $\widetilde{Y^{\ell}}$).

\begin{defn}
$H_{\H}^{2p+k}(X\m Y,\QQ(p)):=H_{\H}^{k}(K^{\bb})$.
\end{defn}
This is computed by the $k^{\text{{th}}}$ cohomology of $(*)$ in
Remark 2.2, or \begin{equation} \begin{matrix} Cone \left\{ \begin{matrix} \tot^{\bb} \left( \tau_{j\leq 0} \B(p) \right) \; \; \; \; \; \oplus \\  [(F^0\tot^{\bb}DR(p)) \cap \tot^{\bb} (\tau_{j\leq 0}DR(p))] \end{matrix} \rTo^{\beta_1 -\beta_2} \tot^{\bb} \left( \tau_{j\leq 0 } DR(p) \right)  \right\} [-1] \\  \\ = \tot^{\bb} \H(p)^{\bb,\bb}, \end{matrix} \end{equation} 
where \begin{equation} \H(p)^{i,j} := \left\{ \begin{matrix}  0 & j>1 \\ \\ \{\ker(d) \subset \D^{2(p+i)}(\widetilde{Y^{-i}}) \} & j=1 \\ \\  \begin{pmatrix} \{ \ker(\d) \subset C^{2(p+i)}(\widetilde{Y^{-i}},\QQ(p+i)) \} \\ \oplus \; \; \{ \ker(d) \subset F^{p+i}\D^{2(p+i)}(\widetilde{Y^{-i}}) \} \\ \oplus \; \; \; \; \D^{2(p+i)-1} (\widetilde{Y^{-i}}) \end{pmatrix} & j=0 \\ \\ \begin{pmatrix} C^{2p+2i+j}(\widetilde{Y^{-i}},\QQ(p+i)) \oplus F^{p+i}\D^{2p+2i+j}(\widetilde{Y^{-i}}) \\ \oplus \; \; \D^{2p+2i+j-1}(\widetilde{Y^{-i}}) \end{pmatrix} & j<0 \end{matrix} \right. \end{equation} 
has differentials $D$ (vertical $Cone$ differential) and $Gy$ (horizontal).
Note that this is $0$ for $i>0$. Elements of $\H(p)^{i,j}$ are
always written $(a,b,c)$, with $a,b$ both zero if $j=1$.

\begin{prop}
The (Gysin) spectral sequence $\H(p)_{\bb}^{\bb,\bb}$ associated
to $\H(p)^{\bb,\bb}$ (with $\H(p)_{0}^{i,j}:=\H(p)^{i,j}$, $d_{0}:=D$)
converges to $H_{\H}^{2p+*}(X\m Y,\QQ(p))$.
\end{prop}
\begin{rem}
This does not in general degenerate at $\H_{2}$; see the proof of
Prop. 2.7 below. It is easy to modify $\H(p)^{\bb,\bb}$ to compute
$H_{\H}^{2p+*}(X\m Y,\AA(p))$ for more general $\AA$. Take the $0^{\text{{th}}}$
column to be the complex from Example 2.3, computing $H_{\H}^{2p+*}(X,\AA(p))$,
and similarly replace the other columns by complexes computing AH
cohomologies of the $\widetilde{Y^{-i}}$.
\end{rem}
The presentation (2.7) of $\tot^{\bb}\H(p)$ as a $Cone$, together
with (2.6) (and variants), leads to a short-exact sequence (s.e.s.)
\begin{equation} \begin{matrix} 0 \to \ext (\QQ(0), H^{2p-m-1}(X \m Y,\QQ(p))) \rTo^{\alpha} H^{2p-m}_{\H}(X \m Y,\QQ(p)) \\ \;\;\;\;\;\;\;\; \rTo^{\beta} \hm (\QQ(0), H^{2p-m}(X\m Y,\QQ(p)) \to 0 . \end{matrix} \end{equation} 
There is also a localization long-exact sequence (l.e.s.) \begin{equation} \to H^{2p-m}_{\H} (X,\QQ(p)) \to H^{2p-m}_{\H} (X\m Y,\QQ(p)) \to H^{2p-m+1}_{\H, Y}(X,\QQ(p)) \to \end{equation} 
produced as follows. If we define double complexes\[
\H_{Y}(p)^{i,j}:=\left\{ \begin{array}{cc}
\H(p)^{i,j} & i<0\\
0 & i\geq0\end{array}\right.,\,\,\,\,\,\,\,\,\H_{X}(p)^{i,j}:=\left\{ \begin{array}{cc}
\H(p)^{0,j} & i=0\\
0 & i\neq0\end{array}\right.,\]
then $\tot\H_{Y}(p)[-1]^{\bb}$, $\tot^{\bb}\H_{X}(p)=\H(p)^{0,\bb}$
have respective $k^{\text{{th}}}$ cohomologies $H_{\H,Y}^{2p+k}(X,\QQ(p)$,
$H_{\H}^{2p+k}(X,\QQ(p))$. The s.e.s. $\H_{X}(p)^{\bb,\bb}\rInto\H(p)^{\bb,\bb}\rOnto$$\linebreak$$\H_{Y}(p)^{\bb,\bb}$
of double-complexes translates to a s.e.s. of complexes\[
\tot^{\bb}\H_{X}(p)\rInto\tot^{\bb}\H(p)\rOnto(\tot\H_{Y}(p)[-1])[+1]^{\bb},\]
which produces (2.10).

The spectral sequence of Proposition 2.5 induces a {}``weight''
filtration with graded pieces $Gr_{j}^{W}H_{\H}^{2p+i+j}(X\m Y,\QQ(p))\cong\H(p)_{\infty}^{i,j}$,
which will concern us (in a broader context) in $\S3$. For now, we
show how to compute these graded pieces. To grease the skids for this
endeavor (and for describing the $AJ$ map in $\S3$), define double
complexes\[
\W_{\QQ}(p)^{i,j}:=\tau_{j\leq0}\B(p)^{i,j}=\left\{ \begin{array}{cc}
0 & j>0\\
\\\{\ker(\d)\subset C^{2(p+i)}(\widetilde{Y^{-i}},\QQ(p+i))\} & j=0\\
\\C^{2p+2i+j}(\widetilde{Y^{-i}},\QQ(p+i)) & j<0\end{array}\right.,\]
$\W(p)^{i,j}:=\tau_{j\leq0}DR(p)^{i,j}$, with upward shifts $\TW_{[\QQ]}(p)^{i,j}:=\W_{[\QQ]}(p)^{i,j-1}$
(and $\widetilde{\B}(p)^{i,j}:=\B(p)^{i,j-1}$). These are each the
$0^{\text{{th}}}$ page of a spectral sequence with $d_{0}=$ vertical
differential. By \cite{De1}, $DR(p)_{\bb}^{\bb,\bb}$ and $\B(p)_{\bb}^{\bb,\bb}$
degenerate at the $2^{\text{{nd}}}$ page; so therefore do the $\W$-spectral
sequences, with e.g.\[
\TW_{\QQ}(p)_{\infty}^{i,j}=\TW_{\QQ}(p)_{2}^{i,j}\cong Gr_{j-1}^{W}H^{2p+i+j-1}(X\m Y,\QQ(p)),\]
and\[
\TW_{\QQ}(p)_{1}^{i,j}\cong H^{2p+2i+j-1}(\widetilde{Y^{-i}},\QQ(p+i)).\]
Also note that besides being a sub- double complex of $DR(p)^{\bb,\bb-1}$,
$\TW(p)^{\bb,\bb}$ injects into $\H(p)^{\bb,\bb}$ by sending $\xi\in\W(p)^{i,j-1}$
to $(0,0,\xi)\in\H(p)^{i,j}$. This produces a map of spectral sequences.

Now write $H^{[k]}:=H^{2p+k}(X\m Y,\QQ(p))$ and consider the s.e.s.
of MHS with $j\leq0$:\[
0\to Gr_{j-1}^{W}H^{[k]}\to\frac{W_{0}}{W_{j-2}}H^{[k]}\to\frac{W_{0}}{W_{j-1}}H^{[k]}\to0.\]
 This leads to a l.e.s.: \SMALL \begin{equation} \begin{matrix} 0 \to \hm (\QQ(0),\frac{W_0}{W_{j-2}}H^{[k]}) \rInto \hm (\QQ(0),\frac{W_0}{W_{j-1}}H^{[k]}) \rTo^{\delta_j^{[k]}} \\ \\ \ext (\QQ(0),Gr^W_{j-1}H^{[k]}) \to \ext (\QQ(0), \frac{W_0}{W_{j-2}}H^{[k]}) \rOnto \ext (\QQ(0),\frac{W_0}{W_{j-1}}H^{[k]}) \to 0, \end{matrix} \end{equation} \normalsize 
and we set $\text{{im}}(\delta_{j}^{[k]})=:\Xi_{j}^{[k]}$. Note that
$k=i+j-1$ $\implies$ $Gr_{j-1}^{W}H^{[k]}\cong\TW(p)_{2}^{i,j}$.

\begin{prop}
For $j=0$ we have a s.e.s.\[
0\to\frac{\ext(\QQ(0),Gr_{-1}^{W}H^{[i-1]})}{\Xi_{0}^{[i-1]}}\to\H(p)_{\infty}^{i,0}\to\hm(\QQ(0),H^{[i]})\to0.\]
Otherwise,\[
\H(p)_{\infty}^{i,j}\cong\left\{ \begin{array}{cc}
0, & j>1\,\,(\text{{or\,}}\, i>0)\\
\\\ext(\QQ(0),Gr_{0}^{W}H^{[i]}), & j=1\\
\\\frac{\ext(\QQ(0),Gr_{j-1}^{W}H^{[i+j-1]})}{\Xi_{j}^{[i+j-1]}}, & j<0\end{array}\right..\]

\end{prop}
\begin{proof}
We first determine $\H(p)_{1}^{\bb,\bb}$ and $\H(p)_{2}^{\bb,\bb}$.
Recall that $H^{*}$ of $\H(p)^{0,\bb}$computes $H_{\H}^{2p+*}(X,\QQ(p))$;
it follows easily from this that the other columns compute AH cohomology
of the $\widetilde{Y^{-i}}$:\SMALL \[
0\to\ext(\QQ(0),\TW_{\QQ}(p)_{1}^{i,j})\to\begin{array}[t]{c}
H_{\H}^{2p+2i+j}(\widetilde{Y^{-i}},\QQ(p+i))\\
||\\
\H(p)_{1}^{i,j}\end{array}\to\hm(\QQ(0),\W_{\QQ}(p)_{1}^{i,j})\to0\]
\normalsize Now $\W_{\QQ}(p)_{1}^{i,j}$ is of pure weight $j$.
Suppose $j\neq0$, so that the r.h. term vanishes and\[
\H(p)_{1}^{i,j}\cong\frac{\TW(p)_{1}^{i,j}}{F^{0}\TW(p)_{1}^{i,j}+\TW_{\QQ}(p)_{1}^{i,j}}.\]
The $\TW(p)_{1}^{i,j}$ are polarized HS's; by semisimplicity, $d_{1}$-cohomology
passes under the $Ext$ and \begin{equation} \H(p)_2^{i,j} \cong \ext (\QQ(0), \TW_{\QQ}(p)^{i,j}_2 ) \lOnto \TW (p)^{i,j}_2 . \end{equation} 
On the other hand, for $j=0$, the s.e.s. of complexes\[
0\to\ext(\QQ(0),\TW_{\QQ}(p)_{1}^{\bb,0})\to\H(p)_{1}^{\bb,0}\to\hm(\QQ(0),\W_{\QQ}(p)_{1}^{\bb,0})\to0\]
yields a $d_{1}$-cohomology l.e.s.\[
\to\hm(\QQ(0),\W_{\QQ}(p)_{2}^{i-1,0})\rTo^{\delta_{0}^{[i-1]}}\ext(\QQ(0),\TW_{\QQ}(p)_{2}^{i,0})\to\H_{2}^{i,0}\]
\[
\to\hm(\QQ(0),\W_{\QQ}(p)_{2}^{i,0})\rTo^{\delta_{0}^{[i]}}\ext(\QQ(0),\TW_{\QQ}(p)_{2}^{i+1,0})\to\]
where we have again used semisimplicity (to pass $d_{1}$ under $Ext$
and $Hom$). Using (2.11), we get a s.e.s. \begin{equation} 0 \to \frac{\ext (\QQ(0),\TW_{\QQ}(p)^{i,0}_2)}{\Xi_{0}^{[i-1]}} \rTo^{\alpha} \H^{i,0}_2 \rTo^{\beta} \hm (\QQ(0),\frac{W_0}{W_{-2}}H^{[i]}) \to 0. \end{equation} 

Now the higher differentials $d_{\ell}$ ($\ell\geq2$) are all $0$
on $\TW(p)_{\bb}^{\bb,\bb}$, which means that the surjectivity (2.12)
extends to the entire map of spectral sequences $\TW(p)_{r}^{\bb,\bb}\to\H(p)_{r}^{\bb,\bb}$
($r\geq2$) off the $0^{\text{{th}}}$ row. Hence the only nontrivial
higher differentials on $\H(p)_{\bb}^{\bb,\bb}$ are the\[
d_{-j+1}:\,\H_{-j+1}^{i+j-1,\,0}\to\H_{-j+1}^{i,\, j}\]
 for $j<0$. (For $j=1$, the proof is now finished.) Moreover, since
$\TW(p)_{2}^{i,0}\to\H(p)_{2}^{i,0}$ factors through, and surjects
onto, the l.h. term of (2.13), we have $d_{2}|_{\text{{im}}(\alpha)}=0$
and so $d_{2}$ factors through $\beta$:\\
\xymatrix{& & \H_2^{i,0} \ar [r]^{d_2 \mspace{150mu}} \ar @{->>} [rd]^{\beta} & \H^{i+2,-1}_2 \cong \ext (\QQ(0),Gr^W_{-2}H^{[i]}) \\ & & & \hm (\QQ(0),\frac{W_0}{W_{-2}}H^{[i]}) \ar [u]_{\delta^{[i]}_{-1}} . \\ } \\
Using again (2.11) we get a s.e.s.\\
\xymatrix{0 \to \frac{\ext (\QQ(0),\TW_{\QQ}(p)^{i,0}_3)}{\Xi_{0}^{[i-1]}} \ar [r]^{\mspace{100mu} \alpha} & \H^{i,0}_3 \ar [r]^{\beta_3 \mspace{130mu}} & \hm (\QQ(0),\frac{W_0}{W_{-3}}H^{[i]}) \to 0 , \\ & \ker(d_2) \ar @{=} [u] \ar @{->>} [r] & \ker(\delta^{[i]}_{-1}) \ar @{=} [u] \\ }\\
and iterating the argument completes the proof for $j=0$. (At $\infty$,
we have r.h. term $\hm(\QQ(0),W_{0}H^{[i]})=\hm(\QQ(0),H^{[i]})$.)
This argument also shows that each $d_{-j+1}$ factors\\
\xymatrix{& & \H_{-j+1}^{i+j-1,0} \ar [r]^{d_{-j+1}} \ar @{->>} [rd]^{\beta_{-j+1}} & \H^{i,j}_{-j+1} \\ & & & \hm (\QQ(0),\frac{W_0}{W_{j-1}}H^{[i+j-1]}) \ar [u]_{\delta^{[i+j-1]}_{j}} , \\ } \\
so that $\text{{im}}(d_{-j+1})=\text{{im}}(\delta_{j}^{[i+j-1]})=\Xi_{j}^{[i+j-1]}$
and we are done.
\end{proof}

\section{\textbf{The weight-filtered Abel-Jacobi map}}

Having dispensed with the preliminaries of $\S2$, we are ready to
explicitly describe the absolute Hodge cycle-class map\[
c_{\H}^{_{p,m}}:\, CH^{p}(U,m)\to H_{\H}^{2p-m}(U,\QQ(p))\]
for smooth quasi-projective $U$. We do this first in terms of a map
of double-complexes ($\S3.1$), then in $\S3.2$ via generalized membrane
integrals. The double-complex construction is then applied in a discussion
of the Beilinson-Hodge conjecture (in $\S3.3$).

\subsection{Gysin spectral sequences and higher residue maps}

We start by laying out a number of complexes: given a smooth projective
algebraic $W/k\subseteq\CC$, set\[
C_{\M}^{\bb}(W,\QQ(q)):=Z_{\RR}^{q}(W,2q-\bb).\]
Also writing $W$ for $W_{\CC}^{\text{{an}}}$, define\\
\[
C_{\B}^{\bb}(W,\QQ(q)):=C_{top}^{\bb}(W,\QQ(q))\,,\,\,\,\,\,\,\,\,\,\,\,\, C_{[F]DR}^{\bb}(W,\CC):=[F^{q}]\D^{\bb}(W)\,,\]
\\
\[
C_{\D}^{\bb}(W,\QQ(q)):=Cone\left\{ C_{top}^{\bb}(W,\QQ(q))\oplus F^{q}\D^{\bb}(W)\to\D^{\bb}(W)\right\} [-1]^{\bb}\]
\[
=C_{top}^{\bb}(W,\QQ(q))\oplus F^{q}\D^{\bb}(W)\oplus\D^{\bb-1}(W)\,,\]
\\
\[
C_{\H}^{\bb}(W,\QQ(q)):=\left\{ \begin{array}{cc}
0\,, & \bb>2q+1\\
\D_{d\text{{-cl}}}^{2q}(W)\,, & \bb=2q+1\\
Z_{top}^{2q}(W,\QQ(q))\oplus F^{q}\D_{d\text{{-cl}}}^{2q}(W)\oplus\D^{2q-1}(W)\,, & \bb=2q\\
C_{\D}^{\bb}(W,\QQ(q))\,, & \bb<2q\end{array}\right.\]
 \\
and let $\K$ stand for $\B$, $DR$, $FDR$, $\D$, $\H$, or $\M$.
(The respective differentials $\partial_{top}$, $d$, $d$, $D$,
$D$, and $\partial_{\B}$ are then uniformly denoted $d_{\K}$.)
We have $H^{*}(C_{[F]DR}^{\bb}(W,\CC))=[F^{q}]H_{DR}^{*}(W,\CC)$,
and for all other $\K$\[
H^{*}(C_{\K}^{\bb}(W,\QQ(q)))\cong H_{\K}^{*}(W,\QQ(q)).\]

Now let $\Y=\bigcup\Y_{i}\subset\X$ be a NCD (in $\X$ smooth projective$/k$),
$\Y_{I}:=\bigcap_{i\in I}\Y_{i}$, and $\Y^{\ell}:=\bigcup_{|I|=\ell}\Y_{I}$,
so that $\widetilde{\Y^{\ell}}\cong\amalg_{_{|I|=\ell}}\Y_{I}$. (By
convention, $\Y^{0}:=\X$.) In what follows, we exclude $\K=[F]DR$
for uniformity of notation, though similar results hold. (In particular,
the spectral sequences $DR(p)$, $FDR(p)$ would degenerate at the
$2^{\text{{nd}}}$ page.) 

Define (mostly 3rd quadrant) double complexes \[
\K_{\X\m\Y}(p)_{0}^{a,b}:=C_{\K}^{2p+2a+b}(\widetilde{\Y^{-a}},\QQ(p+a))\]
 with vertical differential $d_{\K}$, horizontal differential%
\footnote{ignore the $2\pi i$ (in $Gy$) for $\K=\M$.%
}\[
Gy:=2\pi i\sum_{|I|=-a}\sum_{i\in I}(-1)^{\left\langle i\right\rangle _{I\m\{ i\}}}(\iota_{_{\Y_{I}\subset\Y_{I\m\{ i\}}}})_{_{*}}\]
 ($\left\langle i\right\rangle _{J}:=$position in which $i$ occurs
in $J$), and total differential $\DD:=d_{\K}+(-1)^{b}Gy$ (on the
associated simple complex $s^{\bb}\K(p))$.

\begin{example}
$\M(p)_{0}^{\bb,\bb}$ sits in the $3^{\text{{rd}}}$ quadrant, with
upper-right-hand $3\times3$\xymatrix{& \ar [r] & Z^{p-2}(\widetilde{\Y^2}) \ar [r]_{Gy} & Z^{p-1}(\widetilde{\Y^1}) \ar [r]_{Gy} & Z^p(\X) \\ & \ar [r] & Z^{p-2}(\widetilde{\Y^2},1) \ar [u]^d \ar [r]_{Gy} & Z^{p-1}(\widetilde{\Y^1},1) \ar [u]^d \ar [r]_{Gy} & Z^p(\X,1) \ar [u]^d \\ & \ar [r] & Z^{p-2}(\widetilde{\Y^2},2) \ar [r]_{Gy} \ar [u]^d & Z^{p-1}(\widetilde{\Y^1},2) \ar [r]_{Gy} \ar [u]^d & Z^p(\X,2) \, \, . \ar [u]^d \\ & & \ar [u] & \ar [u] & \ar [u] \\ }
\end{example}
Taking $d_{0}:=d_{\K}$ gives Gysin spectral sequences $\K(p)_{r}^{a,b}\begin{array}[t]{c}
\Longrightarrow\\
^{(*=a+b)}\end{array}$$\linebreak$$H^{*}(s^{\bb}\K(p))\cong H_{\K}^{2p+*}(\X\m\Y,\QQ(p)).$
(Omitting the $0$th column would compute $H_{\K,\Y}^{2p+*+1}(\X,\QQ(p))$.)
One finds in each case $\K(p)_{1}^{a,b}\cong$$\linebreak$$H_{\K}^{2p+2a+b}(\widetilde{\Y^{-a}},\QQ(p+a))$;
and for $\K=\B$ we have degeneration at $\B(p)_{2(=\infty)}$ with
$\B(p)_{\infty}^{a,b}=Gr_{b}^{W}H_{(\B)}^{2p+a+b}(\X\m\Y,\QQ(p))$,
where $W_{\bb}$ is the weight filtration of Deligne \cite{De1}.
To generalize the weight filtration to arbitrary $\K$, set $\tilde{W}_{j}s^{\bb}\K(p):=s^{\bb}\K(p)^{\bb\geq-j,\,\bb}$
and define\[
W_{*+j}H_{\K}^{2p+*}(\X\m\Y,\QQ(p)):=im\left\{ H^{*}(\tilde{W}_{j}s^{\bb}\K(p))\to H^{*}(s^{\bb}\K(p))\right\} ,\]
so that $Gr_{b}^{W}H_{\K}^{2p+a+b}(\X\m\Y,\QQ(p))\cong\K_{\X\m\Y}(p)_{\infty}^{a,b}$;
in particular, \[
\K_{\X\m\Y}(p)_{\infty}^{0,b}\cong\underline{H_{\K}^{2p+b}}(\X\m\Y,\QQ(p)):=im\left\{ H_{\K}^{2p+b}(\X,\QQ(p))\to H_{\K}^{2p+b}(\X\m\Y,\QQ(p))\right\} \]
 (lowest weight). While the Gysin spectral sequences for $\K=\M,\H,\D$
do not degenerate at $\K(p)_{2}$, for $\H,\D$ one can show all $d_{r\geq2}:\,\K(p)_{r}^{a,b<0}\to\K(p)_{r}^{a+r,\, b-r+1}$
are zero; precise formulas for the $\H(p)_{\infty}^{a,b}$ were worked
out in Prop. 2.7.

More generally and perhaps more usefully, one can map all $\K_{\infty}$
terms to easily describable targets. First, to compute $H_{\K}$ of
$\Y^{k}\m\Y^{k+1}\cong\widetilde{\Y^{k}}\m\left(\Y^{k+1}\times_{_{\Y^{k}}}\widetilde{\Y^{k}}\right)$
(with good compactification $\widetilde{\Y^{k}}$), we can use\\
\begin{equation} \\
\K_{\Y^k \m \Y^{k+1}}(q)^{\alpha,\beta}_0 \, = \, C_{\K}^{2q+2\alpha+\beta} \left( \widetilde{\Y^{k-\alpha} \times_{\Y^k} \widetilde{\Y^k}}, \QQ(q+\alpha) \right) . \\
\end{equation}\\
Now while the quotient double complex $\K^{\leq-k}(p)^{a,b}:=\left\{ \begin{array}{cc}
0 & a>-k\\
\K(p)^{a,b} & a\leq-k\end{array}\right.$ (of $\K(p)_{0}^{a,b}$) does not compute $H_{\K,\Y^{k}}^{2p+*}(\X,\QQ(p))$
for $k>1$, it does map to (3.1). More precisely, there are maps of
double complexes\[
\K_{\X\m\Y}(p)_{0}^{a,b}\twoheadrightarrow\K_{\X\m\Y}^{\leq-k}(p)_{0}^{a,b}\to\K_{\Y^{k}\m\Y^{k+1}}(p-k)_{0}^{a+k,\, b}\]
inducing maps of spectral sequences, hence of $\infty$ terms\[
\K_{\X\m\Y}(p)_{\infty}^{-k,b}\to\K_{\Y^{k}\m\Y^{k+1}}(p-k)_{\infty}^{0,b}\cong\underline{H_{\K}^{2p-2k+b}}(\Y^{k}\m\Y^{k+1},\QQ(p-k)).\]
This defines $higher$ $residue$ $maps$\[
Res^{-a}:\, Gr_{b}^{W}H_{\K}^{2p+a+b}(\X\m\Y,\QQ(p))\to\underline{H_{\K}^{2p+2a+b}}(\Y^{-a}\m\Y^{-a+1},\QQ(p+a)),\]
which are in general neither injective nor surjective. More concretely,
if we let $U^{-a}:=$ connected component of $\Y^{-a}\m\Y^{-a+1}$,
the $Res^{-a}$ into $\underline{H_{\K}^{2p+2a+b}}(U^{-a},\QQ(p+a))$
is computed by taking $(-a)$ successive residues along a {}``flag''
of connected components $U^{i}\subseteq\Y^{i}\m\Y^{i+1}$ with $\overline{U^{-a}}\subsetneq\overline{U^{-a-1}}\subsetneq\cdots\subsetneq\overline{U^{1}}\subsetneq\overline{U^{0}}=\X$.

Consider the {}``realization'' maps (of double complexes)\[
\M(p)_{0}^{a,b}\to\H(p)_{0}^{a,b}\to\D(p)_{0}^{a,b}\to\B(p)_{0}^{a,b},\]
with the first given by sending $\sZ\in Z_{\RR}^{p+a}(\widetilde{\Y^{-a}},-b)$
to$\linebreak$ $(2\pi i)^{p+a+b}\left((2\pi i)^{-b}T_{\sZ},\Omega_{\sZ},R_{\sZ}\right)$,
and the latter two obvious. (Here $T_{\sZ}\in C_{top}^{2p+2a+b}(\widetilde{\Y^{-a}},\QQ(p+a))$,
$\Omega_{\sZ}\in F^{p+a}\D^{2p+2a+b}(\widetilde{\Y^{-a}})$, $R_{\sZ}\in\D^{2p+2a+b-1}(\widetilde{\Y^{-a}})$
are as in \cite[sec. 5.3-4]{KLM}.) By considering the induced maps
on the cohomologies of the associated simple complexes, we have \emph{$W_{\bb}$-filtered
cycle-class maps} \[
c_{\K}^{p,m}:\, CH^{p}(\X\m\Y,m)\cong H_{\M}^{2p-m}(\X\m\Y,\QQ(p))\to H_{\K}^{2p-m}(\X\m\Y,\QQ(p)),\]
with graded pieces\[
Gr_{b}^{W}c_{\K}^{p,m}:\,\M(p)_{\infty}^{-b-m,b}\to\K(p)_{\infty}^{-b-m,b}.\]

\begin{rem}
If we omit the $0^{\text{{th}}}$ column from the entire construction,
we get cycle-class maps\[
CH^{p-1}(\Y,m-1)\cong H_{\M,\Y}^{2p-m+1}(\X,\QQ(p))\to H_{\K,\Y}^{2p-m+1}(\X,\QQ(p))\]
for cohomology with support on $\Y$.
\end{rem}
Since $\M(p)_{0}^{a,b}=0$ for $b>0$, $H_{\M}^{2p+*}(\X\m\Y,\QQ(p))=W_{0}H_{\M}^{2p+*}(\X\m\Y,\QQ(p))$.
Therefore one consequence of the above is the existence of a commutative
diagram ($\K=\H,\,\D,$ or $\B$)\\
\\
\xymatrix{& CH^p(\X \m \Y ,m) \ar [r]^{c_{\K}} \ar [d]^{Res^m} & W_0 H_{\K}^{2p-m}(\X \m \Y, \QQ (p)) \ar [d]^{Res^m} \\ & CH^{p-m}(\Y^m \m \Y^{m+1}) \ar [r]^{c_{\K}\; \; \; \; \; \; \; \; \; \; \; \;} & {\underline{H_{\K}^{2p-2m}} (\Y^m \m \Y^{m+1},\QQ (p-m))} . } \\
(Note that $CH^{p-m}(\Y^{m}\m\Y^{m+1})=\underline{CH^{p-m}}(\Y^{m}\m\Y^{m+1})$
since we are in the case of algebraic cycles.) 

\begin{example}
We give a simple application of the corresponding diagram for $\K=DR$,
which has upper-r.h. composition \[
CH^{p}(\X\m\Y,m)\rTo^{c_{FDR}}F^{p}W_{2p}H_{DR}^{2p-m}(\X\m\Y,\CC)\rTo^{Res^{m}}F^{p-m}\underline{H^{2p-2m}}(\Y^{m}\m\Y^{m+1},\CC),\]
in the case $p=m\,:=n$.

The construction of \cite{DK} produces families of hypersurfaces\\
\xymatrix{& & & {X_t} \ar @{^(->}[r] \ar [d] & {\X} \ar @{^(->} [r] \ar [d]_{\pi} & {\PP_{\tilde{\Delta}} \times \PP^1} \ar [ld] \\ & & & {\{t\}} \ar @{^(->} [r] & {\PP^1} }\\
in a Toric Fano $n$-fold $\PP_{\tilde{\Delta}}$ ($n=2,3,4$) with
the generic $X_{t}$ nonsingular, and $X_{0}$ a NCD satisfying $\PP_{\tilde{\Delta}}\m X_{0}\cong(\CC^{*})^{n}$
(with coordinates $x_{1},\ldots,x_{n}$). We set $\Y:=X_{0}\subset\X$.
The hypotheses of \cite{DK} Thm. 1 guarantee the existence of $\Xi\in CH^{n}(\X\m\Y,m)$
restricting to $\left\langle \{ x_{1},\ldots,x_{n}\}\right\rangle \in CH^{n}(\X\cap(\CC^{*})^{n},n).$
From this we deduce $\dlog x_{1}\wedge\cdots\wedge\dlog x_{n}\,(=\Omega_{\Xi})$
is a holomorphic form%
\footnote{The assertion is that the pullback of $\bigwedge_{i=1}^{n}\dlog x_{i}$
from $(\CC^{*})^{n}$ to $\X\cap\{(\CC^{*})^{n}\times\PP^{1}\}$ extends
to a holomorphic form on $\X\m X_{0}\times\{0\}$. In fact $\X\m\{\X\cap((\CC^{*})^{n}\times\PP^{1})\}=(D\times\PP^{1})\cup(\Y=X_{0}\times\{0\}),$
but $\Omega_{\Xi}$ can have no residues along $D\times\PP^{1}$.%
} on $\X\m\Y$ representing $c_{DR}(\Xi)$, or equivalently the infinitesimal
invariant of $\{ c_{\H}(\Xi|_{X_{t}})\}\in\Gamma(U\subset\PP^{1},\, R^{n-1}\pi_{*}\CC/\QQ(n))$
as \[
F^{n}H^{n}(\X\m\Y,\CC)\hookrightarrow\Gamma(U,\,\Omega_{U}^{1}\otimes F^{n-1}R^{n-1}\pi_{*}\CC).\]
 Since $Res^{n}(\bigwedge^{n}\dlog x_{i})\in H^{0}(\Y^{n}=\amalg\text{{pts.}},\CC)$
is obviously nonzero, we get generic nontriviality of $\left\langle \Xi|_{X_{t}}\right\rangle \in CH^{n}(X_{t},n).$
\end{example}

\subsection{Generalization of the Griffiths prescription}

The composition of $c_{\H}$ with $\beta$ from (2.9) yields the fundamental
class map\[
cl:\, CH^{p}(\X\m\Y,m)\to\hm(\QQ(0),H^{2p-m}(\X\m\Y,\QQ(p)));\]
we write $cl(\xi)=:[\xi]$. The term {}``Abel-Jacobi map'' refers
either to $c_{\H}$ or its restriction\[
AJ:\, CH_{hom}^{p}(\X\m\Y,m)\to\ext(\QQ(0),H^{2p-m-1}(\X\m\Y,\QQ(p)))=:J^{p,m}(\X\m\Y)\]
to $\ker(cl)$. An interesting consequence of Prop. 2.7 is that $cl$
factors through $\H(p)_{\infty}^{-m,0}$ (via $Gr_{0}^{W}c_{\H}$),
so that its kernel $CH_{hom}^{p}(\X\m\Y,m)\supseteq W_{-1}CH^{p}(\X\m\Y,m)$.

We shall be interested in these maps when $\X\rTo^{\pi}\ms$ is the
$k$-spread of a projective variety $X$ (defined over a field extension
of $k$), and $\Y=\pi^{-1}(D)$ for $D$ a divisor in $\ms$. Relevant
for the application is an extension of the arithmetic Bloch-Beilinson
conjecture to higher Chow groups of smooth projective varieties.

\begin{conjecture}
\emph{{[}BBC{]}} For all smooth projective $\X/\bar{\QQ}$, $\Y=\emptyset$,
the $c_{\H}^{_{p,m}}$ are injective ($\forall p,m$).
\end{conjecture}
For \emph{non}empty $\Y$ (NCD$/\bar{\QQ}$), BBC together with the
Hodge conjecture implies injectivity of the $c_{\H}^{_{p,m}}$ ({}``quasi-projective
BBC''), and even of the $Gr_{b}^{W}c_{\H}^{_{p,m}}$. While injectivity
generally fails for ($\X,\Y$ defined over) $k\subseteq\CC$ of higher
transcendence degree, one might ask whether $c_{\H}^{_{p,m}}$ is
still strictly compatible with $W_{\bb}$; this seems to be a difficult
problem.

Now we describe an approach to $AJ$ to be used in $\S4$. First,
define$\linebreak$$'C_{\H}^{\bb}(W,\QQ(q))\subset C_{\H}^{\bb}(W,\QQ(q))$
($W$ as above) to be the subcomplex consisting of triples of the
form $(0,*,*)$ (no $C_{top}$'s). It maps quasi-isomorphically to\[
C_{\I}^{\bb}(W,\QQ(q)):=\left\{ \begin{array}{cc}
0, & \bb>0\\
\D_{_{d\text{{-cl}}}}^{2q}(W)/F^{q}\D_{_{d\text{{-cl}}}}^{2q}(W), & \bb=0\\
\D^{\bb}(W)/F^{q}\D^{\bb}(W), & \bb<0\end{array}\right..\]
Writing $\widetilde{\I}_{_{\X\m\Y}}(p)_{0}^{a,b}:=C_{\I}^{2p+2a+b-1}(\widetilde{Y^{-a}},\QQ(p+a)),$
we have\[
H^{*}(\tot^{\bb}\widetilde{\I}(p))\cong\frac{W_{2p}H_{DR}^{2p+*-1}(\X\m\Y,\CC)}{F^{p}W_{2p}H_{DR}^{2p+*-1}(\X\m\Y,\CC)}\]
and \begin{equation} {\frac{H^{-m}(\tot^{\bb}\widetilde{\I}(p))}{H^{-m}(\tot \W_{\QQ}(p)[-1]^{\bb})}} \cong J^{p,m}(\X \m \Y). \end{equation}
A class $\xi\in CH_{hom}^{p}(\X\m\Y,m)$ has a $\DD$-closed representative
$\sZ^{\bb}\in\tot^{-m}\M(p)$; this is sent to%
\footnote{here $\bb=i$ indexes entries along the $(-m)^{\text{{th}}}$ diagonal.%
} $(2\pi i)^{p-m}((2\pi i)^{m+\bb}T_{\sZ^{\bb}},\Omega_{\sZ^{\bb}},R_{\sZ^{\bb}})\in\tot^{-m}\H(p)$
by \cite{KLM}, where $(2\pi i)^{p+\bb}T_{\sZ^{\bb}}=\DD\{(2\pi i)^{p+\bb}\Gamma^{\bb}\}$
for $\{(2\pi i)^{p+\bb}\Gamma^{\bb}\}\in\tot^{-m-1}\W_{\QQ}(p)$.
The triple is changed by $\DD$-coboundary to$\linebreak$ $(2\pi i)^{p-m}(0,\Omega_{\sZ^{\bb}},R_{\sZ^{\bb}}+(2\pi i)^{m+\bb}\delta_{\Gamma^{\bb}})$,
then {}``projected'' to $\DD$-closed$\linebreak$ $(2\pi i)^{p-m}\{\overline{R_{\sZ^{\bb}}+(2\pi i)^{m+\bb}\delta_{\Gamma^{\bb}}}\}\in\tot^{-m}\widetilde{\I}(p)$,
which yields the class $AJ(\xi)$ in (3.2) directly. This is used
in $\S4$.

Given an algebraic cycle $\Z\homeq0$ on a smooth projective variety
$X$, the classical \emph{Griffiths prescription} (\cite{Gr}, or
see \cite[12.12]{Le3}) computes $AJ(\Z)$ by integrating suitable
test forms over a topological membrane bounding on $\Z$. With the
last description of $AJ(\xi)$, it is easy to generalize this. Define
a new double complex $\Omega(p)_{i,j}:=\Omega_{_{(\widetilde{\Y^{-i}})^{\infty}}}^{2n-2p-j+1}(\widetilde{\Y^{-i}})$
($n=\dim\X$) with differentials $\mathfrak{I}:\,\Omega(p)_{i,j}\to\Omega(p)_{i-1,j}$
(dual to $Gy$, see $\S8.1$), and $d:\,\Omega(p)_{i,j}\to\Omega(p)_{i,j-1}$.
Letting currents act on $C^{\infty}$ forms (and multiplying the result
by $(2\pi\sqrt{-1})^{-i}$) induces pairings\[
\int:\, F^{n-p+1}\Omega(p)_{i,j}\otimes\widetilde{\I}(p)^{i,j}\to\CC,\]
well-defined since $\dim(\widetilde{\Y^{-i}})=n+i$. These extend
by \emph{summing along the diagonal} to \begin{equation} \sumint : \, \tot_{-m}F^{n-p+1}\Omega(p)_{\bb,\bb} \otimes \tot^{-m}\widetilde{\I}(p)^{\bb,\bb} \rTo \CC ,  \end{equation} 
then by taking cohomology to the perfect pairing\[
F^{n-p+1}H^{2n-2p+m+1}((\X,\Y),\CC)\otimes\frac{H^{2p-m-1}(\X\m\Y,\CC)}{F^{p}}\rTo\CC.\]
By forgetting weights, $J^{p,m}(\X\m\Y)$ maps to the Deligne Jacobian
\SMALL\[
J_{\D}^{p,m}(\X\m\Y):=\frac{H^{2p-m-1}(\X\m\Y,\CC)}{F^{p}H^{2p-m-1}(\X\m\Y,\CC)+H^{2p-m-1}(\X\m\Y,\QQ(p))}\cong\frac{\{ H_{-m}(\tot_{\bb}F^{n-p+1}\Omega(p))\}^{\vee}}{\text{{im}}H^{-m}(\tot^{\bb}\widetilde{\B}(p))},\]
\normalsize and the image of $AJ(\xi)$ is just the functional on
{}``test forms'' $\omega_{\bb}\in\tot_{-m}F^{n-p+1}\Omega(p)$ given
by $\sumintt$ against the {}``generalized membrane'' $\linebreak$$(2\pi i)^{p-m}\{ R_{\sZ^{\bb}}+(2\pi i)^{m+\bb}\delta_{\Gamma^{\bb}}\}$.

\begin{example}
The case of algebraic cycles ($m=0$) on $\X\m\Y$, where $\Y\rInto^{\iota}\X$
is a smooth divisor, is already entertaining. Fix $n(=\dim\X)=5$,
and let $\left\langle \Z\right\rangle \in CH_{hom}^{3}(\X\m\Y)$.
Clearly $\overline{\Z}\in Z_{(\RR)}^{3}(\X)$ gives an element of
$\tot^{0}\M_{_{\X\m\Y}}(3)$ representing $\left\langle \Z\right\rangle $.
Since $[\Z]=0$ in $\hm(\QQ(0),H^{6}(\X\m\Y,\QQ(3))),$ there exist
\[
\Gamma\in C_{top}^{5}(\X,\QQ(3)),\,\gamma\in C_{top}^{4}(\Y,\QQ(2));\,\,\Xi\in F^{3}\D^{5}(\X),\,\xi\in F^{2}\D^{4}(\Y)\]
s.t. $\d\Gamma+(2\pi i)\iota_{*}\gamma=(2\pi i)^{3}\overline{\Z}$,
$d\Xi+(2\pi i)\iota_{*}\xi=(2\pi i)^{3}\delta_{\overline{\Z}}.$ (In
fact, $\Gamma,\,\Xi$ determine $\gamma,\,\xi$.) This yields an element
$\{(\gamma,\xi,0),(\Gamma,\Xi,0)\}\in\tot^{-1}\H_{_{\X\m\Y}}(3),$
while the double-complex $AJ$ applied to $\overline{\Z}$ gives $\linebreak$$(2\pi i)^{3}(\overline{\Z},\delta_{\overline{\Z}},0)\in\tot^{0}\H_{_{\X\m\Y}}(3).$
Adding $\DD$ of the former to the latter gives $\{(0,0,\delta_{\gamma}-\xi),(0,0,\delta_{\Gamma}-\Xi)\}$,
so we have $\{\delta_{\gamma}-\xi,\delta_{\Gamma}-\Xi\}\mapsto\overline{\{\delta_{\gamma},\delta_{\Gamma}\}}\in\tot^{0}\widetilde{\I}_{_{\X\m\Y}}(3).$

Now set\[
\VV:=\frac{\{(\omega,\eta)\,\,|\,\,\omega\in F^{3}\Omega_{\X^{\infty}}^{5}(\X)_{d\text{{-cl}}},\,\eta\in F^{3}\Omega_{\Y^{\infty}}^{4}(\Y),\,\iota^{*}\omega=d\eta\}}{\{(-d\alpha,\iota^{*}\alpha-d\beta)\,\,|\,\,\alpha\in F^{3}\Omega_{\X^{\infty}}^{4}(\X),\,\beta\in F^{3}\Omega_{\Y^{\infty}}^{3}(\Y)\}}\]
\[
\cong H_{0}(\tot_{\bb}F^{3}\Omega(3))\cong F^{3}H^{5}((\X,\Y),\CC)\,,\]
\[
\WW_{\QQ}:=\frac{\{(\tau,\sigma)\,\,|\,\,\tau\in C_{top}^{5}(\X,\QQ(3)),\,\sigma\in Z_{top}^{4}(\Y,\QQ(2)),\,(2\pi i)\iota_{*}\beta=\d\tau\}}{\{(-\d\varepsilon+2\pi i\iota_{*}\varphi,-\d\varphi)\,\,|\,\,\varepsilon\in C_{top}^{4}(\X,\QQ(3)),\,\varphi\in C_{top}^{3}(\Y,\QQ(2))\}}\]
\[
\cong H^{0}(\tot^{\bb}\widetilde{\B}(3))\cong H_{5}((\X,\Y),\QQ(3)).\]
Using (3.3), $\overline{\{\delta_{\gamma},\delta_{\Gamma}\}}$ gives
a functional on $\VV$ by \[
(\omega,\eta)\,\,\,\longmapsto\,\,\,\int_{\Gamma}\omega+(2\pi i)\int_{\gamma}\eta.\]
(Note that if $\gamma$ is algebraic, then $\int_{\gamma}\eta=0$;
but to arrange this in general requires the Hodge conjecture.) This
functional yields an element of $J_{\D}^{3,0}(\X\m\Y)\cong\frac{\VV^{\vee}}{\text{{im}}\WW_{\QQ}}$
which depends only on $\left\langle \Z\right\rangle $ and not on
the choices of $\Gamma,\,\Xi$, etc.
\end{example}

\subsection{Higher Hodge conjecture}

For $X$ smooth projective and defined over a subfield $k$ of $\CC$,
the HC (Hodge Conjecture) states that $cl^{p,0}:\, CH^{p}(X_{(k)})\to\hm(\QQ(0),H^{2p}(X_{\CC}^{an},\QQ(p)))$
should be surjective. Beilinson \cite{Be} and Jannsen \cite{Ja2}
made a prediction extending this to all $p,m$ and beyond the projective
setting, provided the variety is arithmetic ($k\subseteq\bar{\QQ}$). 

\begin{conjecture}
For $U/\bar{\QQ}$ smooth quasiprojective,\[
cl_{U}^{p,m}:\, CH^{p}(U,m)\to\hm(\QQ(0),H^{2p-m}(U_{\CC}^{an},\QQ(p)))\]
is surjective.
\end{conjecture}
To show that this is {}``reasonable'', we relate it to the extended
Bloch-Beilinson conjecture BBC (Conj. 3.4). The next result, due to
M. Saito \cite{SaM4} (by a somewhat less direct inductive proof),
pops right out of the machinery of $\S\S2,3.1$.

\begin{prop}
Conjecture 3.6 follows from HC and BBC.
\end{prop}
\begin{proof}
Let $\X/\bar{\QQ}$ be a good compactification of $U$, i.e. $\Y:=\X\m U$
is a NCD (also defined $/\bar{\QQ}$). Let\[
\widetilde{Res}^{k}:\, W_{b}H_{\K}^{2p-k+b}(\X\m\Y,\QQ(p))\rOnto\K(p)_{\infty}^{-k,b}\]
be the natural map. Consider the diagram \\
\\
\xymatrix{CH^p(\X \m \Y,m) \ar [r]^{c_{\H}} \ar @{->>} [d]^{\widetilde{Res}^m} \ar @/_8pc/ [dd]_{cl} & W_0H^{2p-m}_{\H}(\X \m \Y ,\QQ(p)) \ar @{->>} [d]^{\widetilde{Res}^m} \\ \left\{ \bigcap_{i\geq 2} \ker(d_i) \subseteq \M_{_{\X\m\Y}}(p)_2^{-m,0} \right\} \ar [r]^{\text{KLM}}_{\text{formula}} & \left\{ \bigcap_{i\geq 2} \ker(d_i) \subseteq \H_{_{\X \m \Y}}(p)_2^{-m,0} \right\} \ar [d] \\ \hm(\QQ(0),H^{2p-m}(\X \m \Y,\QQ(p))) \ar @{=} [r] & \left\{ \bigcap_{i\geq 1} \ker(\delta^{[-m]}_{1-i})\subseteq (*) \right\} \\ } 
\\
\\
\\
where $\delta_{j}^{[-m]}$ is defined in $\S2$ and $(*)=$\small \[
\hm(\QQ(0),\B_{_{\X\m\Y}{}}(p)_{2}^{-m,0})\cong\frac{\ker\{ Gy:\, Hg^{p-m}(\widetilde{\Y^{m}})\to Hg^{p-m+1}(\widetilde{\Y^{m-1}})\}}{\text{{im}}\{ Gy:\, Hg^{p-m-1}(\widetilde{\Y^{m+1}})\to Hg^{p-m}(\widetilde{\Y^{m}})\}}.\]
\normalsize (We have made use of the proof of Prop. 2.7. Note that
$Hg^{q}(Y):=\hm(\QQ(0),H^{2q}(Y,\QQ(q)))$.) Also consider \\
\\
\SMALL \xymatrix{CH^{p-m}(\widetilde{\Y^m}) \ar [r]^{Gy=d_1} \ar [d]^{cl^{p-m}} \ar @/_5pc/  [dd]_{c_{\H}^{p-m}} & CH^{p-m+1}(\widetilde{\Y^{m-1}}) \ar [d]_{cl^{p-m+1}} \ar @/^5pc/ [dd]^{c_{\H}^{p-m+1}} \\ Hg^{p-m}(\widetilde{\Y^m}) \ar [r]^{Gy}  & Hg^{p-m+1}(\widetilde{\Y^{m-1}})  \\ H^{2(p-m)}_{\H}( \widetilde{Y^m},\QQ(p-m)) \ar [r]^{Gy=d_1} \ar @{->>} [u] & H_{\H}^{2(p-m+1)}(\widetilde{\Y^{m-1}},\QQ(p-m+1)) \ar @{->>} [u] \\ \ext(\QQ(0),H^{2p-2m-1}(\widetilde{\Y^m},\QQ(p\text{-}m))) \ar [r]^{Gy} \ar @{^(->} [u] & \ext(\QQ(0),H^{2p-2m+1}(\widetilde{\Y^{m-1}},\QQ(p\text{-}m\text{+}1))) \ar @{^(->} [u] \\ } \normalsize 
\\
\\
\\
By the HC, there exists (with $n=\dim(X)$) $\Gamma\in CH^{n-m}(\widetilde{\Y^{m-1}}\times\widetilde{\Y^{m}})$
inducing an inverse of $Gy:\, H^{2p-2m-1}(\widetilde{\Y^{m}})\to H^{2p-2m+1}(\widetilde{\Y^{m-1}})$
on its image, and the $0$-map in other degrees.

Given $\xi\in\hm(\QQ(0),H^{2p-m}(\X\m\Y,\QQ(p)))$, let $\tilde{\xi}$
be a {}``lift'' to $\ker(Gy)\subseteq Hg^{p-m}(\widetilde{\Y^{m}})$.
By HC for $cl^{p-m}$, we have $\sW\in CH^{p-m}(\widetilde{\Y^{m}})$
with $cl^{p-m}(\sW)=\tilde{\xi}$. One easily checks that $c_{\H}^{p-m+1}$
of $\tilde{\sW}:=\sW-\Gamma_{_{*}}Gy(\sW)$ lies in $\ker(Gy)$. By
BBC for $c_{\H}^{p-m+1}$, so does $\tilde{\sW}$ itself. (Of course,
still $cl^{p-m}(\tilde{\sW})=\tilde{\xi}$.)

Now $\overline{\tilde{\xi}}\in\bigcap_{i\geq1}\ker(\delta_{1-i}^{[-m]})\subset(*)$
$\implies$ $c_{\H}^{p-m}(\tilde{\sW})\in\bigcap_{i\geq1}\ker(d_{i})$,
again using the proof of Prop. 2.7. So for $i\geq2$, we have \SMALL \[
c_{\H}^{p-m+i,\, i-1}(d_{i}(\tilde{\sW}))\in\text{{im}}\left\{ \begin{array}{c}
\ext(\QQ(0),H^{2p-2m+i-2}(\widetilde{\Y^{m-i+1}},\QQ(p-m+i-1)))\\
\rTo^{Gy}\ext(\QQ(0),H^{2p-2m+i}(\widetilde{\Y^{m-i}},\QQ(p-m+i)))\end{array}\right\} .\]
\normalsize Let $\Gamma^{i}$ be an algebraic cycle inverse to $Gy$
as above; then writing $\widetilde{d_{i}(\tilde{\sW})}$ for a representative
(of $d_{i}(\tilde{\sW})$) in $CH^{p-m+i}(\widetilde{\Y^{m-i}},i-1)$,
$c_{\H}\{\widetilde{d_{i}(\sW)}-Gy(\Gamma_{_{*}}^{i}\widetilde{d_{i}(\tilde{\sW})})\}=0.$
Using BBC (for $c_{\H,\widetilde{\Y^{m-i}}}^{p-m+i,\, i-1}$) gives
$\widetilde{d_{i}(\tilde{\sW})}\in\text{{im}}(Gy=d_{1})$ $\implies$
$d_{i}(\tilde{\sW})=0$. Hence $\tilde{\sW}$ survives to $\infty$,
and is $Gr_{0}^{W}$ of a class in $CH^{p}(\X\m\Y,m)$.
\end{proof}
\begin{rem}
One case where HC and BBC are actually true in all places they are
used (in the above proof), is where $\X$ and all irreducible components
of each $\Y^{j}$ are rational.
\end{rem}
If (say) $\Y$ is not defined over $\bar{\QQ}$, it is easy to construct
counterexamples to the Conjecture.

\begin{example}
Let $V_{/\bar{\QQ}}\subset\PP^{3}$ be a smooth hypersurface of degree
at least $4$, with $P\in V(\bar{\QQ})$, and $Q\in V(\CC)$ very
general. Let $L_{P}$ and $L_{Q}$ be disjoint lines through $P$
and $Q$ (resp.) intersecting $V$ transversely. (Note that $L_{Q}$
cannot be defined $/\bar{\QQ}$.) We claim the cycle-class maps \begin{equation} CH^2(V\m \{ P,Q \} , 1) \to \hm (\QQ(0),H^3(V \m \{ P,Q \} , \QQ(2))) \end{equation} 
\begin{equation} CH^3(\PP^3 \m V\cup L_P \cup L_Q ,2 ) \to \hm (\QQ(0), H^4(\PP^3 \m V\cup L_P \cup L_Q ,\QQ(3))) \end{equation} 
are not surjective. Consider the blow-ups \[
\beta_{V}:\, B_{P\amalg Q}(V)=:\tilde{V}\rOnto V,\,\,\,\,\,\,\,\,\,\,\beta_{\PP}:\, B_{L_{Q}\amalg L_{P}}(\PP^{3})=:\tilde{\PP^{3}}\rOnto\PP^{3},\]
with (resp.) exceptional divisors $\ell_{P}\amalg\ell_{Q}$ and $S_{P}\amalg S_{Q}$;
set $'\tilde{V}:=\beta_{\PP}^{-1}(V)$.

For (3.4), let $\X=\tilde{V}$, $\Y=\ell_{P}\amalg\ell_{Q}(=\widetilde{\Y^{1}})$;
take $\tilde{\xi}\in Hg^{1}(\Y)$ to be the class of a difference
of points (-pt. on $\ell_{P}$, +pt. on $\ell_{Q}$). Then in fact
$\tilde{\xi}\in\ker\{ Gy=d_{1}:\, H_{\H}^{2}(\Y,\QQ(1))\to H_{\H}^{4}(\X,\QQ(2))\}$,
so $\tilde{\xi}$ determines a class in the r.h.s. of (3.4). But any
$\Z\in CH^{1}(\Y)$ with $cl^{1}(\Z)=\tilde{\xi}$ has $(\beta_{V})_{_{*}}Gy(\Z)=Q-P\in CH^{2}(V_{\CC})$
which is $\nrateq0$; so $d_{1}(\Z)\nequiv0$ and $\Z$ cannot be
$Gr_{0}^{W}$ of a class in the l.h.s. of (3.4).

For (3.5), take $'\X=\tilde{\PP^{3}}$, $'\Y={}'\tilde{V}\cup S_{p}\cup S_{Q}$
$\implies{}'\widetilde{\Y^{1}}={}'\tilde{V}\amalg S_{P}\amalg S_{Q}$,
$'\widetilde{\Y^{2}}=$ disjoint union of lines on $'\tilde{V}$ (including
$\ell_{P},\,\ell_{Q}$). Our $\tilde{\xi}$ from above maps to a $'\tilde{\xi}\in Hg^{1}({}'\widetilde{\Y^{2}})\cong H_{\H}^{2}({}'\widetilde{\Y^{2}},\QQ(1)),$
$Gy\,(=d_{1})$ of which is $0$ in $H_{\H}^{4}({}'\widetilde{\Y^{1}},\QQ(2))\cong Hg^{2}({}'\widetilde{\Y^{1}})$.
Moreover, $d_{2}$ maps $'\tilde{\xi}$ to \[
\frac{H_{\H}^{5}({}'\X,\QQ(3))}{Gy\{ H_{\H}^{3}({}'\widetilde{\Y^{1}},\QQ(2))\}}\cong\frac{H^{4}({}'\X,\CC/\QQ(3))}{Gy\{ H^{2}({}'\widetilde{\Y^{1}},\CC)\}}=0;\]
 so we get a class in the r.h.s. of (3.5). But $d_{1}=Gy$ of any
$'\Z\in CH^{1}({}'\widetilde{\Y^{2}})$ with $cl^{1}({}'\Z)={}'\tilde{\xi}$,
maps to $Q-P\,(\nrateq0)$ under \[
CH^{2}({}'\widetilde{\Y^{1}})\rTo^{\text{{restrict}}}CH^{2}({}'\tilde{V})\rTo^{(\beta_{\PP}|_{_{{}'\tilde{V}}})_{_{*}}}CH^{2}(V_{\CC});\]
so once again we are done.
\end{example}

\subsection{Coniveau filtration}

To conclude this section we mention an alternative approach to sequential
residue maps on $H_{\K}^{2p+*}(\X\m\Y,\QQ(p))$. While the Gysin/weight
setup above has the highest-codimension residue defined first, the
approach via coniveau reverses this (so that codim.-1 residues are
defined first). The drawback is that what follows is valid only for
$\K=\M,\D,$ and $\B$ (not $\H$). 

Define\[
N^{k}C_{\K}^{2p+*}(\X,\QQ(p)):=\]
\[
im\left\{ (2\pi i)^{k}\iota_{_{*}}^{\widetilde{\Y^{k}}}:\, C_{\K}^{2p-2i+*}(\widetilde{\Y^{i}},\QQ(p-i))\to C_{\K}^{2p+*}(\X,\QQ(p))\right\} ;\]
taking \[
\textrm{{}}'\K(p)_{0}^{a,b}:=Gr_{N}^{a}C_{\K}^{2p+a+b}(\X,\QQ(p))\,,\,\, d_{0}:=(Gr_{N}^{a})d_{\K}\]
then yields (4th quadrant) spectral sequences with \[
{}'\K(p)_{1}^{a,b}\cong H_{\K}^{2p-a+b}(\Y^{a}\m\Y^{a+1},\QQ(p-a))\,,\]
\[
{}'\K(p)_{\infty}^{a,b}\cong Gr_{N}^{a}H_{\K}^{2p+a+b}(\X,\QQ(p)).\]
The higher differentials $d_{r}$ on ${}'\K(p)_{r}^{a,b}$ then yield
the desired residue maps\[
{}'Res^{r}:\,\left\{ \ker({}'Res^{r-1})\subseteq H_{\K}^{2p+b}(\X\m\Y,\QQ(p))\right\} \mspace{50mu}\]
\[
\mspace{50mu}\to\left\{ \begin{array}{c}
\text{{subquotient\, of}}\\
H_{\K}^{2p+b-2k+1}(\Y^{r}\m\Y^{r+1},\QQ(p-r))\end{array}\right\} \]
with $\bigcap\ker('Res^{r})=\underline{H_{\K}^{2p+b}}(\X\m\Y,\QQ(p)).$
The $'Res^{r}$ have been studied for $\K=\M,\D$ (in a special case)
in \cite{Ke4}. 

It is intriguing to compare the length of the filtration by their
kernels to that of the {}``weight'' filtration of $\S3.1$. To determine
that an element of $CH^{p}(\X\m\Y,m)$ comes from $CH^{p}(\X,m)$,
one must check the vanishing of $m$ higher $\widetilde{Res}^{k}$
maps (see $\S3.3$) or $[p-\frac{m}{2}]$ (where $[\cdot]=$greatest-integer
function) higher $'Res^{r}$ maps.

\section{\textbf{Higher normal functions}}

The next two sections describe our take on higher Abel-Jacobi maps
for higher Chow groups --- invariants that detect cycles in the regulator
(or $AJ$) kernel. Here we will adopt the formalism of M. Saito's
mixed Hodge modules \cite{SaM2}, while $\S5$ (which is more {}``computational'')
reverts to absolute Hodge cohomology in a special (but broadly applicable)
case.

\subsection{Spreads and the Leray filtration}

We first recall the relationship between these, for $W$ (defined
over $k\subset\CC$) a smooth quasi-projective variety; for us MHM($W$)
will always mean MHM($W_{\CC}$). There is an equivalence of categories
between MHM($pt.$) and graded-polarizable mixed Hodge structures
(PMHS), so that push-forward along the structure morphism $a_{_{W}}:\, W_{\CC}\to\text{{Spec}}\CC$
yields $(a_{_{W}})_{_{*}}:\,\MHM(W)\to D^{b}\PMHS$ satisfying $H^{i}(a_{_{W}})_{_{*}}\QQ_{_{W}}(p)\cong H^{i}(W,\QQ(p))$.
Writing $H^{0}:=H^{2p-m}(W,\QQ(p))$, $H^{-1}:=H^{2p-m-1}(W,\QQ(p))$,
we have a diagram with short-exact rows: \begin{equation} \tiny \xymatrix{{\ext} (\QQ(0),H^{-1}) \ar @{=} [d] & H^{2p-m}_{\H}(W,\QQ(p)) \ar @{=} [d] \\ {\frac{{W_0} H_{\CC}^{-1}}{{F^0} {W_0} H^{-1}_{\CC}+{W_0} H^{-1} } } \ar [r] & Ext^{2p-m}_{_{D^b\MHS}}(\QQ(0),(a_{_{W}})_{_*} \QQ_{_{W}}(p)) \ar [r] & {\hm}(\QQ(0),H^0) \ar @{=} [d]  \\ {\frac{W_{-1} H_{\CC}^{-1}}{\{ F^0 W_0 H^{-1}_{\CC} + W_0 H^{-1}\} \cap W_{-1}H^{-1}_{\CC}}} \ar [r] \ar @{^(->} [u] & Ext^{2p-m}_{_{D^b\PMHS}}(\QQ(0),(a_{_W})_{_*}\QQ_{_{W}}(p)) \ar [r] \ar @{^(->} [u] & Hom_{_{\PMHS}}(\QQ(0),H^0) \\ Ext^1_{_{\PMHS}}(\QQ(0),H^{-1}) \ar @{=} [u] & Ext^{2p-m}_{_{\MHM (W)}} (\QQ_{_{W}}(0), \QQ_{_{W}}(p)) \ar @{=} [u]^{(*)} \\ } \normalsize \end{equation}  
where ({*}) follows from $\QQ_{_{W}}(0)=(a_{_{W}})^{^{*}}\QQ(0)$
and adjointness of $(a_{_{W}})^{^{*}}$, $(a_{_{W}})_{_{*}}$. There
is then a cycle-class map (described in \cite{As})\[
c_{_{\MHM}}^{_{W}}:\, CH^{p}(W,m)\to Ext_{_{\MHM(W)}}^{2p-m}(\QQ_{_{W}}(0),\QQ_{_{W}}(p))\]
factoring $c_{\H}$ (as defined in $\S3$, thinking of $W$ as having
a good compactification).

Next, we briefly revisit the approach to spreads discussed in \cite{Ke2}.
Begin with the following data:

\begin{itemize}
\item fields $k\subset K\subset\CC$, $k$ algebraically closed and $K/k$
finitely generated of transcendence degree $t$;
\item a smooth projective variety $X$ defined $/K$, of dimension $d$;
and
\item a higher Chow (pre)cycle $\Z\in Z^{p}(X_{(K)},m)$ with $\db\Z=0$.
\end{itemize}
The theory of spreads then provides:

\begin{itemize}
\item a smooth projective variety $\ms$ $/k$, with an isomorphism $k(\ms)\to K$;
and $p_{g}:\,\text{{Spec}}K\to\ms$ the corresponding generic point.
(Clearly $\dim_{k}\ms=t$.) Denote by $U\subset\ms$ any affine Zariski-open
subvariety defined $/k$.
\item $\X$ smooth projective and $\pi:\,\X\to\ms$ projective, both defined
$/k$, s.t. $X=\X\times_{p_{g}}\text{{Spec}}K$. Write $\pi^{-1}(U)=:\X_{U}$.
\item $\bar{\sZ}\in Z^{p}(\X,m)$ with $\Z=\bar{\sZ}\times_{p_{g}}\text{{Spec}}K$,
and some $U_{0}\subset\ms$ (as above) s.t. $\sZ_{U_{0}}:=\bar{\sZ}|_{\X_{U_{0}}}$
is $\db$-closed. {[}Note: by $\bar{\sZ}$ we do not necessarily mean
the Zariski closure, though that would be one possible choice.{]}
\end{itemize}
\begin{rem}
It is possible that none of the many choices of $\bar{\sZ}$ is $\db$-closed;
on the other hand, if such a choice (with $\db\bar{\sZ}=0$) does
exist, we call it a \emph{complete $k$-spread} of $\Z$. Existence
or nonexistence of a complete spread depends only on the class of
$\Z$ in $CH^{p}(X_{(K)},m)$.
\end{rem}
To get rid of ambiguities (i.e., the nonuniqueness of $\X$, $\ms$,
$\bar{\sZ}$, etc.), we take the limit over $U\subset\ms$ ($U$ as
above): $\eta_{\ms}:=\varprojlim U$, $\X_{\eta}:=\varprojlim\X_{U}$.
These have precise meaning only under a contravariant functor to abelian
groups,%
\footnote{the resulting direct limits are exact in $\underline{\text{{Ab}}}$.%
} e.g.\[
\ext(\QQ(0),H^{a}(\gs,R^{b}\pi_{_{*}}\QQ(p))):=\varinjlim_{U}\ext(\QQ(0),H^{a}(U,R^{b}\pi_{_{*}}\QQ(p))),\]
\[
H_{_{DR}}^{*}(\X_{\eta},\CC):=\varinjlim_{U}H^{*}(\D^{\bb}((\X_{U})_{\CC}^{an})),\]
\[
CH^{p}(\X_{\eta},m):=\varinjlim_{U}CH^{p}(\X_{U},m)\rTo_{p_{g}^{^{*}}}^{\cong}CH^{p}(X,m).\]
The \emph{$k$-spread} $\sZ$ of $\Z$ is just the restriction of
$\sZ_{U_{0}}$ to $Z^{p}(\X_{\eta},m);$ writing $\mathfrak{{s}}:=(p_{g}^{^{*}})^{-1}$,
the class $\left\langle \sZ\right\rangle \in CH^{p}(\X_{\eta},m)$
(frequently denoted just $\sZ$) is then $\mathfrak{{s}}\left\langle \Z\right\rangle $.
Note that the complex points of $\ms$ that survive in the limit over
$U$ are the very general points (of maximal transcendence degree),
corresponding to the various embeddings of $k(\ms)$ into $\CC$.

The spectral sequence \begin{equation} E^{i,j}_2 := Ext^i_{_{\MHM (\gs)}} \left( \QQ_{_{\eta}}(0), H^j \pi_{_*}\QQ_{_{\X_{\eta}}} (p) \right) \end{equation} 
computing \begin{equation} Ext^{2p-m}_{_{\MHM (\X_{\eta}) }}  \left( \QQ_{_{\X_{\eta}}}(0), \QQ_{_{\X_{\eta}}}(p) \right) \cong Ext^{2p-m}_{_{D^b \MHM (\gs )}} \left( \QQ_{_{\eta}}(0), \pi_{_*} \QQ_{_{\X_{\eta}}}(p) \right) \end{equation} 
degenerates at $E_{2}$ since $\pi$ is proper; so the resulting Leray
filtration $\L^{\bb}$ on (4.3) has $Gr_{\L}^{i}$(4.3)$\cong E_{2(=\infty)}^{i,2p-m-i}$.
Set $\L^{i}CH^{p}(X,m):=(c_{_{\MHM}}^{_{\X_{\eta}}}\circ\mathfrak{s})^{-1}\L^{i}$,
so that $c_{_{\MHM}}^{_{\X_{\eta}}}$ induces\[
c_{_{\MHM}}^{_{X,i}}:\, Gr_{\L}^{i}CH^{p}(X,m)\rInto Ext_{_{\MHM(\gs)}}^{i}(\QQ_{_{\eta}}(0),H^{2p-m-i}\pi_{_{*}}\QQ_{_{\X_{\eta}}}(p)).\]
The filtration on $CH^{p}(X,m)$ stabilizes $\L^{i}=\L^{i+1}$ $\forall i\geq\text{{min}}\{ t+2,p+1\}$
by a standard argument (using relative hard Lefschetz for $\X_{\eta}\to\gs$,
see \cite{Ke1} or \cite{Le1}). If $k=\bar{\QQ}$ then the {}``quasi-projective
BBC'' of $\S3.2$ implies that $c_{\H}^{_{\X_{\eta}}}$, hence $c_{_{\MHM}}^{_{\X_{\eta}}}$,
is injective; so $\L^{\bb}$ (conjecturally) stabilizes at $\{0\}$,
hence should be viewed as a candidate Bloch-Beilinson filtration.

\begin{rem}
Note that this is essentially the filtration described in \cite{As},
except without the limit over all $K\subseteq\CC$ finitely generated
over $k$ (and containing the minimal field of definition of $X,\,\Z$).
\end{rem}
A second spectral sequence (for the structure morphism $a_{_{\eta}}:\gs\to\text{{Spec}}k$)
yields short-exact sequences $\underline{E}_{\infty}^{i,j}\rInto E_{\infty}^{i,j}\rOnto^{p^{i,j}}\underline{\underline{E}}_{\infty}^{i,j}$
where \[
\underline{E}_{\infty}^{i,j}=Ext_{_{\PMHS}}^{1}(\QQ(0),H^{i-1}(\gs,R^{j}\pi_{_{*}}\QQ(p))),\]
\[
\underline{\underline{E}}_{\infty}^{i,j}=Hom_{_{\text{{(P)MHS}}}}(\QQ(0),H^{i}(\gs,R^{j}\pi_{_{*}}\QQ(p))).\]
(Arapura \cite{Ar2} has proved that the MHS on $H^{i}(U,R^{j}\pi_{_{*}}\QQ(p))$
arising from Saito's theory agrees with the {}``natural'' MHS coming
from his geometric definition \cite{Ar1} of the Leray filtration
on $H^{*}(\X_{U},\QQ(p))$.) In particular, the $[\sZ]_{i}:=p^{i,2p-m-i}(c_{_{\MHM}}^{_{X,i}}(\Z))$
are just Leray graded pieces of the fundamental class $[\sZ]\in\hm(\QQ(0),H^{2p-m}(\X_{\eta},\QQ(p)))$.
Both are obviously (well-)defined (given $\left\langle \Z\right\rangle $)
without recourse to any of the theory just discussed (except spreads).

\subsection{A geometric construction: the $\Lambda$-filtration}

We now propose a notion of \emph{higher normal functions} having the
$[\sZ]_{i}$ as {}``topological invariants'', which arose from an
attempt to extend Arapura's geometric approach to Leray to absolute
Hodge cohomology. The definition does not use MHM, and leads to a
more concrete description of $\L^{\bb}$.

Recall that a complex subvariety $\T_{\CC}\subset\ms_{\CC}$ is called
\emph{very general} if no rational function $f\in\bar{\QQ}(\ms)^{*}$
has $f|_{\T_{\CC}}\equiv0$. Equivalently, the minimal field of definition
$L(\subset\CC)$ of $\T$ has $\text{{trdeg}}(L/k)=\text{{codim}}_{\ms}(\T)$;
$\T$ always means $\T_{L}$, and $\eta_{\T}:=\varprojlim\V$ over
$\V\subset\T$ affine Zariski open $/L$. Write $\pi$ (by abuse of
notation) for the proper morphism $\X_{\eta_{\T}}\to\eta_{\T}$, and
$\ms[i]$ for the set of $(i-1)$-dimensional very general subvarieties
of $\ms_{(\CC)}$.

\begin{defn}
If $[\sZ]_{0}=\cdots=[\sZ]_{i-1}=0$, the \emph{$i^{\text{{th}}}$
higher normal function}\[
\nu_{\sZ}^{i}:\,\ms[i]\to\coprod_{\T\in\ms[i]}\ext(\QQ(0),H^{2p-m-1}(\X_{\eta_{\T}},\QQ(p)))\]
\emph{associated to $\sZ$}, is given by $\nu_{\sZ}^{i}(\T):=AJ(\sZ|_{\X_{\eta_{\T}}})$.
\end{defn}
\begin{rem}
This can be {}``higher'' in two senses ($i>1$ and $m>0$). If $i=1$
and $m=0$, $\nu_{\sZ}^{1}$ is essentially a classical normal function,
restricted to $\eta_{\ms}(\CC)=\varprojlim U(\CC)=\ms[1]$.
\end{rem}
To see that this is well-defined, note that $\pi:\,\X_{\eta}\to\eta(=\eta_{\T}\text{{\, or\,\,}}\eta_{\ms})$
is proper and smooth, so that $E_{2}^{a,b}=H^{a}(\eta,R^{b}\pi_{_{*}}\QQ(p))$
converges at $E_{2}$ and the (canonical) Leray filtration has $Gr_{\L}^{a}H^{a+b}(\X_{\eta},\QQ(p))=E_{2,\eta}^{a,b}$.
If $[\sZ]_{a}\in E_{2,\gs}^{a,2p-m-a}$ vanishes for $0\leq a\leq i-1$,
then so do $[\sZ|_{\X_{\eta_{\T}}}]\in E_{2,\eta_{\T}}^{a,2p-m-a}$,
while $E_{2,\eta_{\T}}^{a,2p-m-a}=0$ for $a\geq i$ (as $\dim(\eta_{\T})=i-1$
and $\eta_{\T}$ is a limit of affines); hence $[\sZ|_{\X_{\eta_{\T}}}]=0$
and $AJ$ is defined.

Moreover, we have (noncanonically) $H^{*}(\X_{\eta},\QQ(p))\cong\oplus_{a+b=*}E_{2,\eta}^{a,b}$
as MHs (with $\L^{q}=\oplus_{a\geq q}$), so that all $\L^{q+1}\to\L^{q}\to Gr_{\L}^{q}$
split and $L^{a}\ext(\QQ(0),H^{a+b}(\X_{\eta},\QQ(p))):=\ext(\QQ(0),\L^{a}H^{a+b})$
$\linebreak$satisfies $Gr_{L}^{a}\cong\ext(\QQ(0),E_{2,\eta}^{a,b})=:\E_{2,\eta}^{a,b}.$
If $\T'\subset\T$ is a (very general) hyperplane section, consider
$\theta^{n}:\, H^{n}(\X_{\eta_{\T}},\QQ(p))\to H^{n}(\X_{\eta_{\T'}},\QQ(p))$,
with $\theta^{a,b}:=Gr_{\L}^{a}\theta^{a+b}:\, E_{2,\eta_{\T}}^{a,b}\to E_{2,\eta_{\T'}}^{a,b}$
injective for $a\leq i-2\,(=\dim(\T')).$ If also $2a+b<2p$, then
weights of $E_{2,\eta_{\T}}^{a,b}$ are $<0$, so that $\hm(\QQ(0),\text{{coker}}(\theta))=0$
and $\E_{2,\eta_{\T}}^{a,b}\rInto_{Ext^{1}(\theta^{a,b})}\E_{2,\eta_{\T'}}^{a,b}.$
In particular, if $i\leq m+2$, then $Ext^{1}(\theta^{a,2p-m-a-1})$
is injective for $0\leq a\leq i-2$, while $\E_{2,\eta_{\T}}^{a,2p-m-a-1}$
(resp. $\E_{2,\eta_{\T'}}^{a,2p-m-a-1}$) is $0$ for $a\geq i$(resp.
$i-1$); thus $\ker(Ext^{1}(\theta^{2p-m-1}))\cong\E_{2,\eta_{\T}}^{i-1,2p-m-i}.$
This proves (a) of

\begin{prop}
Assume $[\sZ]_{0}=\cdots=[\sZ]_{i-1}=0$.\\
\emph{(a)} If $\nu_{\sZ}^{i-1}=0$ and $i\leq m+2$, then $\nu_{\sZ}^{i}$
factors through \[
\coprod_{\T\in\ms[i]}\ext(\QQ(0),H^{i-1}(\eta_{\T},R^{2p-m-1}\pi_{_{*}}\QQ(p)));\]
 \emph{(b)} If $\nu_{\sZ}^{i}=0$, then $[\sZ]_{i}=0$. In fact, $[\sZ]_{i}$
is the {}``topological invariant'' of $\nu_{\sZ}^{i}$ (to be defined
below).
\end{prop}
Assuming (b), we can use the normal functions to define a filtration
$\Lambda^{\bb}$ without reference to vanishing of the $\{[\sZ]_{i}\}$.

\begin{defn}
Set $\Lambda^{0}CH^{p}(X,m)=CH^{p}(X,m)$, $\Lambda^{1}CH^{p}(X,m)=CH_{hom}^{p}(X,m)$
($\Lambda^{0}=\Lambda^{1}$ if $m>0$), and\[
\Lambda^{i}CH^{p}(X,m):=\left\{ \Z\in CH_{hom}^{p}(X,m)\,\,\,|\,\,\,\nu_{\mathfrak{{s}}(\Z)}^{1}=\cdots=\nu_{\mathfrak{{s}}(\Z)}^{i-1}=0\right\} .\]

\end{defn}
To facilitate the proof of (b), we restrict the domains of the $\{\nu_{\sZ}^{i}\}$.
\emph{Choose} nested very general hyperplane sections $(\ms_{\CC}\supsetneq)\T_{\CC}^{1}\supsetneq\T_{\CC}^{2}\supsetneq\cdots\supsetneq\T_{\CC}^{t-1},$
with $\T^{j}$ of codimension $j$ and minimal field of definition
$L_{j}$. These satisfy $k\subsetneq L_{1}\subsetneq L_{2}\subsetneq\cdots\subsetneq L_{t-1}\subsetneq K$
with $\text{{trdeg}}(L_{j}/k)=j$. One can show that there exist smooth
affine varieties and morphisms\\
\xymatrix{& {\M_0} \ar @{=} [d] & & & & & {\M_t} \ar @{=} [d] \\ (\ms & )U \ar @{_(->} [l] \ar [r]^{\rho_1} \ar @/^2pc/ [rr]^{\rho_2} \ar @/_1pc/ [rrrr]_{\rho_{_{t-1}}} & \M_1 \ar [r]^{\rho_{12}} & \M_2 \ar [r] & \cdots \ar [r] & \M_{t-1} \ar [r] & \text{Spec}(k) \\ }\\
defined $/k$ (with $\dim(\M_{j})=t-j$, $k(\M_{j})\cong L_{j}$),
such that fibers of each $\rho_{j}$ are affine $j$-dimensional (multiple-)hyperplane
sections of $U$. (In particular, $(\T^{j}\cap U)_{\CC}=\rho_{j}^{-1}(q_{j})$
for some very general $q_{j}\in\M_{j}(\CC)$.) Denote the restriction
$\eta_{\ms}\to\eta_{\M_{j}}$ of $\rho_{j}$ by $\rho_{j}^{\eta}$.
Given $\mu\in\M_{j}$, write $\T_{\mu}:=\rho_{j}^{-1}(\mu)$, and
set $H_{\mu}:=H^{2p-m-1}(\X_{\eta_{\T_{\mu}}},\QQ(p)).$

The restricted normal function \[
\bar{\nu}_{\sZ}^{i}:\,\eta_{\M_{i-1}}(\CC)\to\coprod_{\mu\in\eta_{\M_{i-1}}(\CC)}\ext(\QQ(0),H_{\mu})\]
still clearly satisfies Prop. 4.5(a) (this will be important later).
Define a local system on $\eta_{\M_{i-1}}$ with stalks $W_{0}H_{\mu}$
by $\HH_{i}:=W_{0}R^{2p-m-1}(\rho_{i-1}^{\eta}\circ\pi)_{_{*}}\QQ(p),$
and set $\H_{i}:=\HH_{i}\otimes\O_{\eta_{\M_{i-1}}}$; take $\F_{i}\subset\H_{i}$
to be the subsheaf with stalks $F^{p}W_{2p}H^{2p-m-1}(\X_{\eta_{\T_{\mu}}},\CC),$
and\[
\FF_{i}:=\ker\left\{ \HH_{i}\to\frac{\H_{i}}{\F_{i}}\right\} \,,\,\,\,\,\,\J_{i}:=\text{{coker}}\left\{ \frac{\HH_{i}}{\FF_{i}}\rInto\frac{\H_{i}}{\F_{i}}\right\} .\]
Here $\J_{i}$ has stalks $\J^{p,m}(\X_{\eta_{\T_{\mu}}})$, $\FF$
is a VMHS of pure type $(0,0)$ (since its stalks are $\hm(\QQ(0),H_{\mu})$),
and $\bar{\nu}_{\sZ}^{i}$ may be viewed as a holomorphic section
of $\J_{i}$. Applying the connecting homomorphism\[
H^{0}(\eta_{\M_{i-1}},\J_{i})\rTo^{\delta}H^{1}\left(\eta_{\M_{i-1}},\frac{\HH_{i}}{\FF_{i}}\right)\]
gives the \emph{topological invariant} $\delta\bar{\nu}_{\sZ}^{i}$.

Thinking of $\eta_{\M}$ as some sufficiently small affine $\V\subset\M_{i-1}$
(in the limit), cover it with analytic balls $V_{\alpha}$. Notice
that $\frac{\H_{i}}{\F_{i}}\cong\RR^{-m}(\rho_{i-1}^{\eta}\circ\pi)_{_{*}}\tot^{\bb}\widetilde{\I}_{_{\left(\X_{\gs}/\eta_{\M_{i-1}}\right)}}(p)$,
so that we can apply the version of $AJ$ via (3.2) continuously fiberwise
over each $V_{\alpha}$ to get a section (omitting the {}``$\bb$'')
$(2\pi i)^{p-m}\left\{ R_{\alpha,\mu}+(2\pi i)^{m}\delta_{\Gamma_{\alpha,\mu}}\right\} $,
where $R_{\alpha,\mu}=R_{\sZ|_{\X_{\eta_{\T_{\mu}}}}}$, $\d\Gamma_{\alpha,\mu}=T_{\sZ|_{\X_{\eta_{\T_{\mu}}}}}$.
Over $V_{\alpha}\cap V_{\beta}$, these differ by only $(2\pi i)^{p}\delta_{\Gamma_{\alpha,\mu}-\Gamma_{\beta,\mu}}$
(which are rational cycles); hence we get a Cech cocycle in $H^{1}(\eta_{\M},\HH/\FF),$
and this is $\delta\bar{\nu}_{\sZ}^{i}$.

Now the fundamental class $[\sZ]$ is given by $(2\pi i)^{p}T_{\sZ}$,
which (fibered by $\rho_{i-1}^{\eta}\circ\pi$) gives a trivial section%
\footnote{we are still assuming $[\sZ]_{0}=\cdots=[\sZ]_{i-1}=0.$%
} of $\HH_{i}\cong\RR^{-m}(\rho_{i-1}^{\eta}\circ\pi)_{_{*}}\tot^{\bb}\W_{_{\QQ,(\X_{\gs}/\eta_{\M_{i-1}})}}(p)$
and thus maps to some $\overline{[\sZ]}\in W_{0}H^{1}(\eta_{\M_{i-1}},\HH_{i})$.
(This is by Leray for $\rho_{i-1}^{\eta}\circ\pi$; while the relevant
s.s. doesn't degenerate at $E_{2}$, the $E_{\infty}^{1,*}$ terms
inject into $E_{2}^{1,*}$.) This $\overline{[\sZ]}$ is calculated
by bounding the fibers $T_{\sZ}|_{\X_{\eta_{\T_{\mu}}}}$ over the
$V_{\alpha}$ and taking differences over $V_{\alpha}\cap V_{\beta}$;
reusing the $\{\Gamma_{\alpha,\mu}\}$ shows immediately that the
image of $\overline{[\sZ]}$ under $\jmath:\, W_{0}H^{1}(\eta_{\M_{i-1}},\HH_{i})\to H^{1}(\eta_{\M_{i-1}},\HH_{i}/\FF_{i})$
is $\delta\bar{\nu}_{\sZ}^{i}$. Moreover, since $\FF_{i}$ is a locally
constant system of (limits of) pure HS of weight $0$, $W_{0}H^{1}(\eta_{\M},\FF)=0$
and $\jmath$ is injective. So from $\delta\bar{\nu}_{\sZ}^{i}$ we
recover $\overline{[\sZ]}$, and map this to $H^{1}(\eta_{\M_{i-1}},R^{2p-m-1}(\rho_{i-1}^{\eta}\circ\pi)_{_{*}}\QQ(p)),$
which has a filtration $\LL$ (using Leray for $\pi$ proper) with
graded pieces $Gr_{\LL}^{a}\cong H^{1}(\eta_{\M_{i-1}},R^{a-1}(\rho_{i-1}^{\eta})_{_{*}}R^{2p-m-a}\pi_{_{*}}\QQ(p)).$
By Leray for $\rho_{i-1}^{\eta}$ and vanishing of $[\sZ]_{0}$ thru
$[\sZ]_{i-1}$, $\overline{[\sZ]}$ lives in $\LL^{i}$; its image
in $Gr_{\LL}^{i}$ identifies with that of $[\sZ]_{i}$ under\[
\kappa_{i}:\, H^{i}(\gs,R^{2p-m-i}\pi_{_{*}}\QQ(p))\to H^{1}(\eta_{\M_{i-1}},R^{i-1}(\rho_{i-1}^{\eta})_{_{*}}R^{2p-m-i}\pi_{_{*}}\QQ(p)).\]
(This exists by Leray for $\rho_{i-1}^{\eta}$ and the vanishing of
$R^{i}(\rho_{i-1}^{\eta})_{_{*}}(\cdots)$ by fiber dimension.) Fix
$q_{i}\in\eta_{\M_{i-1}}(\CC)$, and let $\C_{i}:=\rho_{i-1,i}(q_{i})\subset\M_{i-1}$,
$W_{i}:=\rho_{i}^{-1}(q_{i})$ (multi-hyperplane section of $U\subset\ms$).
By affine weak Lefschetz and dimension, we have\\
\SMALL \xymatrix{H^i(\gs,R^{2p-m-i}\pi_{_*}\QQ(p)) \ar @{^(->} [r]_{\gamma_i} \ar [d]^{\kappa_i} & H^i(\eta_{W_i},R^{2p-m-i}\pi_{_*}\QQ(p)) \ar [d]^{\cong} \\ H^1(\eta_{\M_{i-1}},R^{i-1}(\rho^{\eta}_{i-1})_{_*}R^{2p-m-i}\pi_{_*}\QQ(p)) \ar [r] & H^1(\eta_{\C_i},R^{i-1}(\rho^{\eta}_{i-1})_{_*}R^{2p-m-i}\pi_{_*}\QQ(p))} \normalsize \\
which shows $\kappa_{i}$ injective. This completes the proof of Prop.
4.5(b).

\subsection{Comparing filtrations}

In fact, we have shown $\bar{\nu}_{\sZ}^{i}=0\,\,\implies\,\,[\sZ]_{i}=0$,
and so we can define a third filtration on $CH^{p}(X,m)$ by $'\Lambda^{0}:=\Lambda^{0}$,
$'\Lambda^{1}:=\Lambda^{1}$, and\[
'\Lambda^{i}CH^{p}(X,m):=\left\{ \Z\in CH_{hom}^{p}(X,m)\,\,\,|\,\,\,\bar{\nu}_{\mathfrak{{s}}(\Z)}^{1}=\cdots=\bar{\nu}_{\mathfrak{s}(\Z)}^{i-1}=0\right\} .\]
Evidently $\Lambda^{i}\subset\,'\Lambda^{i};$ obviously $'\Lambda^{\bb}$
depends \emph{a priori} on our choices of hyperplane sections. 

\begin{thm}
$\L^{r}\subseteq\Lambda^{r}\subseteq{}'\Lambda^{r}$ \emph{(}$\forall r$\emph{),}
with equality for $r=0,\ldots,m+2$.
\end{thm}
\begin{proof}
First some notation. For $\mu\in\eta_{\M_{i-1}}(\CC)$, let \[
\tilde{H}_{\mu}:=H^{i-1}(\eta_{\T_{\mu}},R^{2p-m-i}\pi_{_{*}}\QQ(p));\]
define a local system $\tilde{\HH}_{i}:=W_{0}R^{i-1}(\rho_{i-1}^{\eta})_{_{*}}R^{2p-m-i}\pi_{_{*}}\QQ(p)$
on $\eta_{\M_{I-1}}$ (with stalks $W_{0}\tilde{H}_{\mu}$), and $\tilde{\H}_{i}$,
$\tilde{\F}_{i}$, $\tilde{\FF}_{i}$, $\tilde{\J}_{i}$ by analogy
to $\H_{i}$, $\F_{i}$, $\FF_{i}$, $\J_{i}$ above. Set $\tilde{\HH}_{i}^{pol}:=W_{-1}\tilde{\HH}_{i}$
and $\tilde{\J}_{i}^{pol}:=\text{{im}}\{\tilde{\H}_{i}^{pol}\to\tilde{\J}_{i}\}$;
the Jacobian sheaves have stalks $(\tilde{\J}_{i})_{\mu}=\ext(\QQ(0),\tilde{H}_{\mu})$
and $(\tilde{\J}_{i}^{pol})_{\mu}=$$\linebreak$ $Ext_{_{\PMHS}}^{1}(\QQ(0),\tilde{H}_{\mu})$.
Using the analysis preceding Prop. 4.5, $\tilde{\J}_{i}^{pol}\rInto$$\linebreak$$\tilde{\J}_{i}\rInto\J_{i}$.

Now we explain some details in the big commutative diagram (4.4) on
the next page. First, $Ext(\gamma_{i})$ is injective for $i\leq m+1$,
because $\gamma_{i}$ is injective ($\forall i$) and $i\leq m+1$
$\implies$ weights of $H^{i-1}(\eta_{W_{i-1}},R^{2p-m-i}\pi_{_{*}}\QQ(p))$
are $\leq2(i-1)+(2p-m-i)-2p=i-(m+2)<0.$ Also, $Hom(\kappa_{i})$
is injective because $\kappa_{i}$ is (see above). Next, the local
system underlying $H^{i-1}(\rho_{i-1}^{\eta})_{_{*}}H^{2p-m-i}\pi_{_{*}}\QQ_{\X_{\gs}}(p)$
is smooth (as the limit over $\V\subset\M_{i-1}$ allows for removal
of jumps), so the $Ext_{_{\MHM(\eta_{\M_{i-1}})}}^{1}$ can be replaced
by extensions $Ext_{_{\VMHS(\eta_{\M})_{ad}}}^{1}$ of admissible
variations of MHS. By M. Saito's description \cite[Rmk. 3.7(ii)]{SaM1}
of $\Delta$, it is clear that $\jmath'\circ\Delta=\tilde{\delta}\circ j\circ\varepsilon_{i}$
(the {}``pentagon'' in (4.4) commutes).

More precisely, consider an admissible VMHS on an affine (Zar. op.)
$\V\subset\M_{i-1}$ with (smooth) underlying locals system $\tilde{\HH}$,
and an extension $\tilde{\HH}\to\EE\to\QQ_{_{\V}}(0)$ with class
$e\in Ext_{_{\VMHS(\V)_{ad}}}^{1}(\QQ_{_{\V}}(0),\tilde{\HH}).$ Associate
to this a {}``normal function'' $\epsilon(e)\in H^{0}(\V,\tilde{\J})$
whose value at $v\in\V$ is the class of $\tilde{\HH}_{v}\to\EE_{v}\to\QQ(0)$
in $\tilde{\J}_{v}^{pol}=Ext_{_{\PMHS}}^{1}(\QQ(0),\tilde{\HH}_{v}).$
To produce $\Delta(e)\in\hm(\QQ(0),H^{1}(\V,\tilde{\HH}))$, consider
a cover $\{ V_{\alpha}\}$ of $\V$ by acyclic analytic balls, and
local {}``liftings'' (morphisms of local systems) $\sigma_{\alpha}:\QQ_{_{\V}}(0)\to\EE$;
and take differences over the intersections. If $\epsilon(e)$ is
zero, i.e. if the pointwise extensions of MHS split, then $e$ splits
locally and we may take the $\sigma_{\alpha}$ to be morphisms of
VMHS; it follows that $\Delta(e)$ factors through $\hm(\QQ(0),H^{1}(\V,\tilde{\FF})),$
which is zero since $\tilde{\FF}$($=${}``$Hg(\H)$'') has pure
weight $0$ ($\implies$ $H^{1}(\V,\tilde{\FF})$ has weights $>1$).
In the limit over $\V$ and $U$ (with $\tilde{\HH}$ essentially
$\tilde{\HH}_{i}$), $\epsilon$ gives $\varepsilon_{i}$; so we have
established $\varepsilon_{i}$'s injectivity.

We easily show $\L^{r}\subseteq{}'\Lambda$ $(\forall r)$; since
the proof doesn't depend on fibers of $\rho_{i}$ being hyperplane
sections, $\L^{r}\subseteq\Lambda^{r}$ follows as well. Assume $\L^{i}\subseteq{}'\Lambda^{i}$,
and let $\Z\in\L^{i}CH^{p}(X,m)$ (and $\sZ:=\mathfrak{s}(\Z)$).
It is easy to see that $(j\circ\varepsilon_{i}\circ\alpha_{i})(c_{_{\MHM}}^{_{X,i}}(\Z))=:\tilde{\bar{\nu}}_{\sZ}^{i}$
maps to $\bar{\nu}_{\sZ}^{i}$ under $H^{0}(\tilde{\J}_{i})\rInto H^{0}(\J_{i})$.
If $\Z\in\L^{i+1}$, then $\tilde{\bar{\nu}}_{\sZ}^{i}=0$ $\implies$
$\bar{\nu}_{\sZ}^{i}=0$ $\implies$ $\Z\in{}'\Lambda^{i+1}$, done.\begin{equation} $$

\includegraphics[%
  width=1.20\textheight,
  height=8in,
  angle=180]{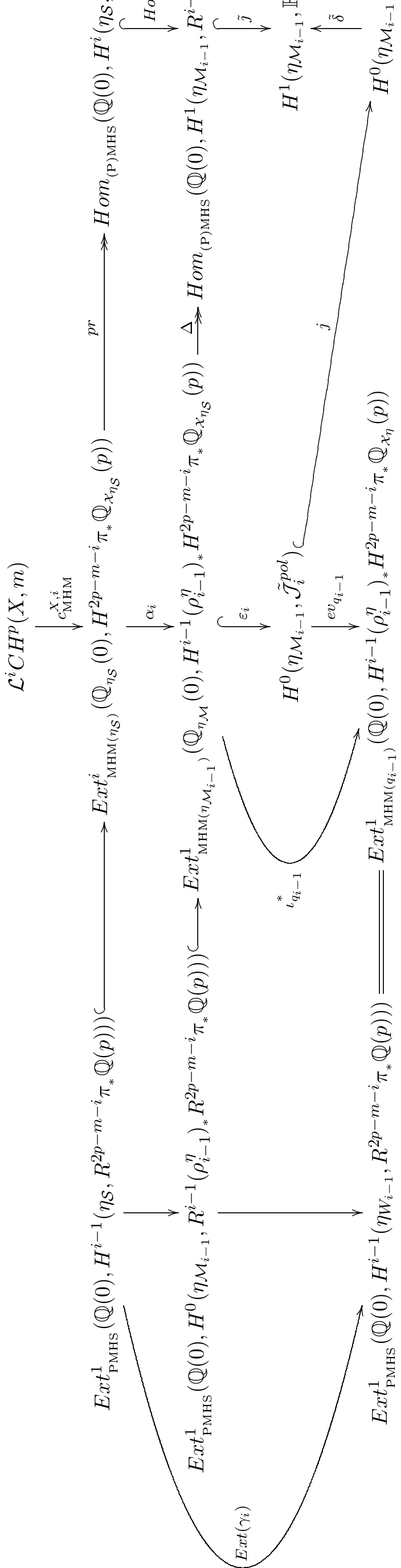}$$ \end{equation}

We remark that whenever $\bar{\nu}_{\sZ}^{i}$ lifts to (or {}``factors
through'') some $\tilde{\bar{\nu}}_{\sZ}^{i}\in H^{0}(\eta_{\M_{i-1}},\tilde{\J}_{i})$,
the topological-invariant map $\delta$ lifts to $\tilde{\delta}$
in (4.4). Such a lifting occurs (for any $i$) when $\Z\in\L^{i}$
as above, but also (for $i\leq m+2$) if $\Z\in{}'\Lambda^{i}$ (since
Prop. 4.5 is valid for \emph{restricted} higher normal functions as
well).

Assume inductively $\L^{r}={}'\Lambda^{r}$ $\forall r\leq i$, where
$i\leq m+1$. To prove $\L^{i+1}={}'\Lambda^{i+1}$ it suffices to
show $'\Lambda^{i+1}\cap\L^{i}\subseteq\L^{i+1}$. Given $\Z\in{}'\Lambda^{i+1}\cap\L^{i}$,
$(j\circ\varepsilon_{i}\circ\alpha_{i})(c_{_{\MHM}}^{_{X,i}}(\Z))$
is zero since it goes to $\bar{\nu}_{\sZ}^{i}$ under an injective
map. Hence $(\Delta\circ\alpha_{i})(c_{_{\MHM}}^{_{X,i}}(\Z))=0$
$\implies$ $pr(c_{_{\MHM}}^{_{X,i}}(\Z))=0$. Injectivity of $Ext(\gamma_{i})$
and an obvious diagram chase show that in fact $c_{_{\MHM}}^{_{X,i}}(\Z)=0$
and so $\Z\in\L^{i+1}$.
\end{proof}
\begin{rem}
(i) The (weak) analogue here of Zucker's theorem on normal functions
\cite[sec. 9]{Z} is the following trivial consequence of the diagram
(4.4): the image of $\tilde{\delta}:\,\{\text{{im}}(j\circ\varepsilon_{i})\subseteq H^{0}(\eta_{\M_{i-1}},\tilde{\J}_{i})\}\to H^{1}(\eta_{\M_{i-1}},\tilde{\HH}_{i}/\tilde{\FF}_{i})$
is precisely the {[}rational{]} $(p,p)$ classes in the target. (We
stress that this analogy is imperfect, as we do not have a concept
of horizontality for sections of $\tilde{\J}_{i}$.)

(ii) For $\Z\in\L^{i}$ ($\forall i$) or $\Z\in{}'\Lambda^{i}$ ($i\leq m+2$),
the diagram actually gives another proof%
\footnote{we have proved this for $\Z\in{}'\Lambda^{i}$ ($\forall i$).%
} that $\bar{\nu}_{\sZ}^{i}$ has $[\sZ]_{i}$ as its topological invariant,
since $\tilde{\delta}\tilde{\bar{\nu}}_{\sZ}^{i}=(\tilde{\delta}\circ j\circ\varepsilon_{i}\circ\alpha_{i})(c_{_{\MHM}}^{_{X,i}}(\Z))=(\tilde{\jmath}\circ Hom(\kappa_{i})\circ pr)(c_{_{\MHM}}^{_{X,i}}(\Z))=(\tilde{\jmath}\circ Hom(\kappa_{i}))[\sZ]_{i}.$
This proves the conjecture ({}``Expectation'') of \cite[sec. 7.1]{Ke2}.
\end{rem}
We finish with some remarks on infinitesimal invariants. Assume our
affine $U\subset\ms$ is sufficiently small that $\pi:\,\X_{U}\to U$
is smooth. Define $\nabla J^{p,i,m}(\X_{U}/U)$ to be the cohomology
of \[
\Gamma(\Omega_{U}^{i-1}\otimes\RR^{2p-m-i}\pi_{_{*}}\Omega_{\X_{U}/U}^{\bb\geq p-i+1})\rTo^{\nabla}\Gamma(\Omega_{U}^{i}\otimes\RR^{2p-m-i}\pi_{_{*}}\Omega_{\X_{U}/U}^{\bb\geq p-i})\]
\[
\rTo^{\nabla}\Gamma(\Omega_{U}^{i+1}\otimes\RR^{2p-m-i}\pi_{_{*}}\Omega_{\X_{U}/U}^{\bb\geq p-i-1})\]
at the middle term; M. Saito \cite{SaM3} constructed a map\[
\Phi_{p,i,m}:\,\hm\left(\QQ(0),H^{i}(U,R^{2p-m-i}\pi_{_{*}}\QQ(p))\right)\to\nabla J^{p,i,m}(\X_{U}/U).\]
Passing to the $\varinjlim_{U}$, the \emph{infinitesimal invariant}
of $\nu_{\sZ}^{i}$ (or $\bar{\nu}_{\sZ}^{i}$) is given by\[
d\nu_{\sZ}^{i}:=\Phi_{p,i,m}(\delta\bar{\nu}_{\sZ}^{i})=\Phi_{p,i,m}([\sZ]_{i})\in\nabla J^{p,i,m}(\X_{\eta}/\gs).\]

\begin{rem}
For $m=0$, {}``higher normal functions'' $\rho_{\X_{\eta}/\gs}^{p,i}(\sZ)$
were defined by S. Saito \cite{SaS}. These, too, formally take $\Phi_{p,i,m}([\sZ]_{i})$
as their infinitesimal invariants. However, for $i\geq2$ they are
missing the {}``higher $AJ$'' data contained in $\nu_{\sZ}^{i}$
(and indeed contain no more information than $\Phi_{p,i,m}([\sZ]_{i})$
itself).
\end{rem}

\section{\textbf{Leray filtration in the product case}}

Imagine that the smooth projective {}``$X$'' treated in the discussion
of spreads in $\S4$, was actually $X_{K}=X\otimes_{k}K$ for some
$X$ defined $/k$. The $k$-spread geometry is then $\pi:\,\X=X\times\ms\to\ms$,
where $K\cong k(\ms)$; if also $\Z\in Z^{p}(X_{K},m)$, then $\sZ\in Z^{p}((X\times\gs)_{k},m)$.
This is the setup for $\S\S5,6,7$.

\subsection{Higher Abel-Jacobi maps}

For $W$ smooth projective we let (as in \cite{RS}) $\K_{\H}^{\bb}(W)$
denote a (canonically chosen) complex of MHS with cohomology computing
$H^{*}(W,\QQ)$. We have a noncanonical isomorphism \[
\K_{\H}^{\bb}(W)\simeq\bigoplus_{j}H^{j}(W,\QQ)[-j]\]
in the bounded derived category $D^{b}MHS$. In particular, \[
\K_{\H}^{\bb}(X\times\gs)\,\simeq\,\K_{\H}^{\bb}(\gs)\otimes\K_{\H}^{\bb}(X)\,\simeq_{n.c.}\]
\begin{equation}  \begin{matrix} \bigoplus_j  \K_{\H}^{\bb}(\gs) \otimes H^{2p-m-j}(X,\QQ) [-2p+m+j] \\ \end{matrix} \end{equation} 
\[
\implies\, H_{\H}^{2p-m}(X\times\gs,\QQ(p))\cong Ext_{_{D^{b}MHS}}^{2p-m}(\QQ(-p),\K_{\H}^{\bb}(X\times\gs))\,\simeq\]
\begin{equation} \begin{matrix}  \bigoplus_j Ext^j_{_{D^bMHS}}(\QQ(-p),\K^{\bb}_{\H}(\gs)\otimes H^{2p-m-j}(X)) \\ \end{matrix} \end{equation}
yields an obvious Leray filtration on $H_{\H}^{2p-m}(X\times\gs,\QQ(p))$,
via $\L^{i}=\oplus_{j\geq i}$. This is canonical, and is easily shown
to be strictly compatible with $\L^{\bb}$ of $\S4$, under the map\\
\xymatrix{Ext^{2p-m}_{_{D^bPMHS}}(\QQ(0),(a_U)_*\pi_*\QQ_{_{X\times U}}(p)) \ar @{^(->} [r] \ar @{=} [d] & Ext^{2p-m}_{_{D^bMHS}}(\QQ(0),(a_U)_*\pi_*\QQ_{_{X\times U}}(p)) \ar @{=} [d] \\ Ext^{2p-m}_{_{MHM(X\times U)}}(\QQ_{_{X\times U}}(0),\QQ_{_{X\times U}}(p)) & H^{2p-m}_{\H}(X\times U,\QQ(p))} \\
(after taking $\varinjlim_{U}$). Simply identify $\K_{\H}^{\bb}(X\times U)(p)$
with $(a_{_{X\times U}})_{*}\QQ_{_{X\times U}}(p)\simeq(a_{U})_{*}\pi_{*}\QQ_{_{X\times U}}(p)$,
and use the filtration of (5.1) by $\oplus_{\geq i}$ to produce (compatibly)
both (5.2) and (4.2).

Consequently we may reproduce the $\L^{\bb}$ of $\S4$ on $CH^{p}(X_{K},m)$
(in our present special case) by pulling $\L^{\bb}$ back along the
composite\[
CH^{p}(X_{K},m)\begin{array}{c}
_{\cong}\\
\longrightarrow\\
^{\mathfrak{{s}}:=k\text{{-spread}}}\end{array}CH^{p}\left((X\times\gs)_{k},m\right)\begin{array}[t]{c}
\to\\
^{c_{\H}}\end{array}H_{\H}^{2p-m}(X\times\gs,\QQ(p)).\]
We always denote $\mathfrak{{s}}\left\langle \Z\right\rangle $ by
$\left\langle \sZ\right\rangle $ (or just $\sZ$). We will be particularly
interested in the restriction of $\L^{\bb}$ (by intersection) to
\[
\underline{CH^{p}}(X_{K},m):=\mathfrak{{s}}^{-1}\left(im\left\{ CH^{p}\left((X\times S)_{k},m\right)\to CH^{p}\left((X\times\gs)_{k},m\right)\right\} \right),\]
 the {}``higher Chow cycles with complete $k$-spread'' $\left\langle \bar{\sZ}\right\rangle \in CH^{p}(X\times S,m)$. 

Recall that for smooth quasi-projective $Y$, $\underline{H^{n}}(Y)$
for $Y$ denotes $im\{ H^{n}(\bar{Y})\to H^{n}(Y)\}$ for any good
compactification $\bar{Y}$, or equivalently the lowest weight filtrand
$W_{n}H^{n}(Y)$. (Dually, $H_{n}(Y):=coim\{ H_{n}(Y)\to H_{n}(\bar{Y})\}$.)
Calculating for $m>0$\\
\[
im\left\{ \begin{array}{c}
Ext_{^{D^{b}MHS}}^{i}\left(\QQ(-p),\K_{\H}^{\bb}(\ms)\otimes H^{2p-m-i}(X)\right)\\
\to Ext_{^{D^{b}MHS}}^{i}\left(\QQ(-p),\K_{\H}^{\bb}(\gs)\otimes H^{2p-m-i}(X)\right)\end{array}\right\} \]
\\
\[
\cong im\left\{ \begin{array}{c}
Ext_{^{MHS}}^{1}\left(\QQ(-p),H^{i-1}(\ms)\otimes H^{2p-m-i}(X)\right)\\
\to Ext_{^{MHS}}^{1}\left(\QQ(-p),H^{i-1}(\gs)\otimes H^{2p-m-i}(X)\right)\end{array}\right\} \]
\\
\[
\cong im\left\{ \begin{array}{c}
Ext_{^{MHS}}^{1}\left(\QQ(-p),\underline{H^{i-1}}(\gs)\otimes H^{2p-m-i}(X)\right)\\
\to Ext_{^{MHS}}^{1}\left(\QQ(-p),H^{i-1}(\gs)\otimes H^{2p-m-i}(X)\right)\end{array}\right\} \]
 \\
\begin{equation} \cong \frac{\ext\left( \QQ(-p), \underline{H^{i-1}}(\gs) \otimes H^{2p-m-i}(X) \right)}{im \left\{ \hm \left( \QQ(-p), \left\{ \frac{W_{m+i}}{W_{i-1}} H^{i-1}(\gs) \right\} \otimes H^{2p-m-i}(X) \right) \right\}} \end{equation} \\
we have our basic invariants (for $m>0$)\\
\\
\SMALL \xymatrix{ {Gr_{\L}^i \underline{CH^p}(X_K,m)} \ar @{^(-->} [r] \ar @{^(~>}  [dd] |{Gr^i_{\L} \underline{c_{\H}}} & {Gr_{\L}^i CH^p(X_K,m)} \, \, \ni \, \Z  \ar @{^(~>} @/_15pc/ [dd] |{Gr^i_{\L}c_{\H}} \\ & { Hom_{MHS} \left( \QQ(-p), H^i(\gs) \otimes H^{2p-m-i}(X) \right) } \\  {Gr^i_{\L} \underline{H^{2p-m}_{\H}}(X\times \gs,\QQ(p))}  \ar @{^(->} [r] & {Gr^i_{\L}H^{2p-m}_{\H}(X\times \gs,\QQ(p)) } \ar @{>>} [u]^{\beta} \\  {\frac{Ext^1_{MHS} \left( \QQ(-p),\underline{H^{i-1}}(\eta_S)\otimes H^{2p-m-i}(X) \right) }{im \left[ Hom_{MHS} \left( \QQ(-p), \frac{W_{m+i}}{W_{i-1}}H^{i-1}(\gs) \otimes H^{2p-m-i}(X) \right) \right] }}  \ar @{^(->} [r] \ar [u]^{\cong} & {Ext^1_{MHS} \left( \QQ(-p),H^{i-1}(\gs)\otimes H^{2p-m-i}(X) \right) } \ar @{^(->} [u]^{\alpha} \\} \normalsize \\
\\
\\
sending $\Z\longmapsto[c_{\H}(\sZ)]_{i}\begin{array}[t]{c}
\longmapsto\\
^{\beta}\end{array}[\sZ]_{i}=cl_{i}(\Z)$, and if this is $0$ (pulling back along $\alpha$) to $[AJ(\sZ)]_{i-1}=AJ_{i-1}(\Z)$.
As in $\S4$, these invariants vanish for $i>min\{ t+1,p\}$. They
are called (resp.) $i^{\text{{th}}}$ higher cycle-class and $(i-1)^{\text{{st}}}$
higher Abel-Jacobi class (of $\Z$).

We want to say something about computing these; a different approach
will be taken in $\S7$. The trouble is that the only datum we have
from spreading out is $\sZ_{U_{0}}\in Z^{p}(Z\times U_{0},m)_{_{\db-\text{{cl.}}}}$.
To produce a double-complex cocycle $\sZ^{\bb}\in\tot^{-m}\M_{_{X\times U_{0}}}(p)$
from this might require an explicit (possibly difficult) application
of Bloch's moving lemma. (The one exception to this is the case $m=0$
of algebraic cycles.) So we will make do without the double-complex
$AJ$-formulas of $\S3$.

To compute $[\sZ]_{i}$, we use the $C^{\infty}$ chain $\overline{T_{\sZ_{U_{0}}}}$
(viewed as a topological correspondence on $X\times(\ms,D_{0})$,
where $D_{0}:=\ms\m U_{0}$) to induce a map $H_{i}(\gs,\QQ)\to H_{2d-2p+m+i}(X,\QQ)$.
This map is nonzero if and only if $[\sZ]_{i}\in Hom_{^{MHS}}(\QQ(-p),H^{i}(\gs)\otimes H^{2p-m-i}(X))$
is.

When $[\sZ]_{i}$ vanishes, one would like an explicit description
of $[AJ(\sZ)]_{i-1}$ as a collection of {}``generalized membrane
integrals'' --- i.e., as functionals on test $C^{\infty}$ forms.
Even if $[\sZ]\in Hom_{^{MHS}}(\QQ(-p),H^{2p-m}(X\times\gs,\QQ))$
vanishes, this can be computationally horrendous. The simplest description
(in general) of $F^{d+t-p+1}\{ H^{2t-i+1}((\ms,D_{0}),\CC)\otimes H^{2d-2p+m+i}(X,\CC)\}$
in terms of $C^{\infty}$ forms uses Hodge-filtered K\"unneth components
of $\tot^{\bb}\Omega(p)^{\bb,\bb}$ from $\S3$ (with $\X=X\times\ms$,
$\Y=X\times D_{0}$). (It is well-known that $C^{\infty}$ forms with
compact support on $U_{0}$ are not amenable to filtration by {[}Hodge{]}
type.)

The one really nice case is when $\Z\in\L^{i}\underline{CH^{p}}(X_{K},m)$,
so that there exists $\bar{\sZ}\in CH^{p}(X\times\ms,m)$ --- automatically
$\homeq0$ for $m>0$ (and we $assume$ for the moment $0=[\bar{\sZ}]\in Hom_{^{MHS}}(\QQ(-p),H^{2p}(X\times\ms))$
if $m=0$). The \cite{KLM} current $(2\pi i)^{p-m}R_{\bar{\sZ}}\pm(2\pi i)^{p}\delta_{\partial^{-1}T_{\bar{\sZ}}}=:(2\pi i)^{p-m}R_{\bar{\sZ}}'$
yields a functional on ($C^{\infty}$ forms representing) $F^{d+t-p+1}\{ H^{2t-i+1}(\ms,\CC)\otimes H^{2d-2p+m+i}(X,\CC)\},$
modulo periods. Now since\[
Hom_{MHS}\left(\QQ(-p),\frac{W_{m+i}}{W_{i-1}}H^{i-1}(\gs)\otimes H^{2p-m-i}(X)\right)\]
\[
\lTo^{\cong}Hom_{MHS}\left(\QQ(-p),\frac{W_{m+i}}{W_{i-1}}H^{i-1}(\gs)\otimes F_{h}^{p-i+1}H^{2p-m-i}(X)\right),\]
 we can kill the {}``$im\{ Hom\}$'' in (5.3) by projecting\small \[
Gr_{\L}^{i}\underline{H_{\H}^{2p-m}}(X\times\gs,\QQ(p))\twoheadrightarrow Ext_{^{MHS}}^{1}\left(\QQ(-p),\underline{H^{i-1}}(\gs)\otimes\frac{H^{2p-m-i}(X)}{F_{h}^{p-i+1}H^{2p-m-i}(X)}\right);\]
\normalsize the image of $[AJ(\sZ)]_{i-1}$ by this is denoted $[AJ(\sZ)]_{i-1}^{tr}=AJ_{i-1}^{tr}(\Z)$.
The functional described above (an element of $Ext_{^{MHS}}^{1}(\QQ(-p),H^{i-1}(\ms)\otimes H^{2p-m-i}(X))$
projects to compute this $[AJ(\sZ)]_{i-1}^{tr}$.

In short, we have constructed a slight {}``quotient'' of $AJ_{i-1}$,\[
AJ_{i-1}^{tr}:\,\L^{i}\underline{CH^{p}}(X_{K},m)\to\ext\left(\QQ(-p),\underline{H^{i-1}}(\gs)\otimes\frac{H^{2p-m-i}(X)}{F_{h}^{p-i+1}}\right),\]
that can be described explicitly in the spirit of the {}``Griffiths
prescription''; it will be used in $\S6$.

\subsection{Weight filtration redux}

Our next aim is to Leray-filter the $\H(p)$-spectral sequence and
the residue maps from $\S3$, when $\X=X\times\ms$. Recall from (5.2)
that for $W$ smooth quasiprojective \[
H_{\H}^{M}(X\times W,\QQ(q))\cong\bigoplus_{\nu}Ext_{^{D^{b}MHS}}^{M-\nu}\left(\QQ(-\nu),\K_{\H}^{\bb}(W)\otimes H^{\nu}(X)\right)\]
\[
=:\bigoplus_{\nu}H_{\H}^{M-\nu}\left(W,H^{\nu}(X,\QQ(q))\right)=:\bigoplus_{\nu}Gr_{\nu}^{\L}H_{\H}^{M}\left(X\times W,\QQ(q)\right)\]
\[
\cong\bigoplus_{\nu}Gr_{\L}^{M-\nu}H_{\H}^{M}\left(X\times W,\QQ(q)\right).\]
If $W$ is projective, there are obvious cochain complexes \[
Gr_{\nu}^{\L}C_{\H}^{\bb}(X\times W,\QQ(q))\]
for computing the graded pieces, with $\oplus_{\nu}$ $\simeq$ to
the usual $C_{\H}^{\bb}$. For $\bb<2q$ the $Gr_{\nu}^{\L}$ complex
is given by\[
C_{top}^{\bb-\nu}\left(W,H^{\nu}(X,\QQ(q))\right)\oplus F^{q}\left(D^{\bb-\nu}(W)\otimes H^{\nu}(X,\CC)\right)\oplus\]
\[
\left(D^{\bb-\nu-1}(W)\otimes H^{\nu}(X,\CC)\right),\]
 for $\bb>2q+1$ it is $0$, and the reader can fill in the rest.

Passing to the spread situation $\X=X\times\ms$, $\Y^{i}=X\times D^{i}$
for $D\subset\ms$ NCD, we form double complexes\[
Gr_{\nu}^{\L}\H_{\X\m\Y}(p)_{0}^{a,b}:=Gr_{\nu}^{\L}C_{\H}^{2p+2a+b}(X\times\widetilde{D^{-a}},\QQ(p+a))\]
with $s^{\bb}\H(p)_{0}\simeq\oplus_{\nu}s^{\bb}Gr_{\nu}^{\L}\H(p)_{0}$;
in fact\[
\H(p)_{1}^{a,b}\cong\bigoplus_{\nu}H_{\H}^{2p+2a+b-\nu}(\widetilde{D^{-a}},H^{\nu}(X,\QQ(p+a))\cong\bigoplus_{\nu}Gr_{\nu}^{\L}\H(p)_{1}^{a,b},\]
and the $\cong$ persists for all pages $r\geq1$. Hence the $Gr_{\nu}^{\L}$
spectral sequences converge to $Gr_{\nu}^{\L}H_{\H}^{2p+*}(X\times(\ms\m D),\QQ(p)),$
and as before we get weight filtrations with maps\[
Gr_{b}^{W}Gr_{\nu}^{\L}H_{\H}^{2p+a+b}(X\times(\ms\m D),\QQ(p))\cong Gr_{\nu}^{\L}\H(p)_{\infty}^{a,b}\]
\[
\rTo^{Res^{-a}}Gr_{\nu}^{\L}\underline{H_{\H}^{2p+2a+b}}\left(X\times(D^{-a}\m D^{-a+1}),\QQ(p+a)\right).\]
 We proceed to the main application.

With the understanding that the limit (over $U$ affine Zariski open
in $\ms$) defining $\gs$ is taken together with a system of good
compactifications (or of $(\tilde{\ms},\tilde{D})$ with $\tilde{D}$
a NCD), we can define%
\footnote{Warning: one can also define $W_{*}$ directly on $CH^{p}(X_{K},m)$,
and $Gr_{\L}$ and $Gr^{W}$ do not {}``commute''. Saying the $Gr_{\L}$
of a class is in $W_{\ell}$ does not mean the whole class is in $W_{\ell}$.%
} \[
W_{\ell}Gr_{\L}^{i}CH^{p}(X_{K},m):=\mathfrak{{s}}^{-1}\left(W_{\ell}Gr_{\L}^{i}CH^{p}((X\times\gs)_{L},m)\right).\]
The following is a simple means for constructing nontrivial $m$-decomposable
higher cycles in the weight-graded pieces.

\begin{thm}
Let $\tau_{1},\ldots,\tau_{m}\in\CC$ be algebraically independent
$/\bar{K}$, and set $L:=K(\tau_{1},\ldots,\tau_{m})$. The exterior
product with $\{\tau_{1},\ldots,\tau_{m}\}\in CH^{m}(\text{{Spec}}L,\, m)$
gives a map\[
Gr_{\L}^{i}\underline{CH^{p}}(X_{K},n)\hookrightarrow Gr_{-n}^{W}Gr_{\L}^{i+m}CH^{p+m}(X_{L},n+m).\]

\end{thm}
\begin{proof}
$L$ is the compositum of $L_{0}:=k(\tau_{1},\ldots,\tau_{m})$ and
$K\cong k(\ms)$; writing $\T:=\PP_{t_{1}}^{1}\times\cdots\times\PP_{t_{m}}^{1}$
we have $L_{0}\cong k(\T)$, $L\cong k(\ms\times\T)$. Taking $D,E$
NCD's in $\ms,\T$ (resp.), with $E\supseteq\cup_{j}|(t_{j})|$, set
$\mathfrak{{D}}:=(\ms\times E)\cup(D\times\T)$. The $k$-spread of
the exterior product $\Z\times\{\underline{\tau}\}$ is the exterior
product of \[
\sZ\in\L^{i}\underline{CH^{p}}(X\times(\ms\m D),\, n)\text{{\,\,\, and\,\,\,}}\{\underline{t}\}\in[\L^{m}]CH^{m}\left([\text{{pt.}}\times]\T\m E,\, m\right).\]
 Assuming $\Z$ is nonzero in $Gr_{\L}^{i}\underline{CH^{p}}(X_{K},n)$,
the class of $\sZ$ is nontrivial (in the limit over all $D$) in 

\begin{equation} \\
Gr^i_{\L} \underline{CH^p}(X\times (\ms \m D), n) \rInto^{Gr^i_{\L}\underline{c_{\H}}} \underline{H_{\H}^i} ( \ms \m D , H^{2p-n-i}(X,\QQ(p)) ) \\
\end{equation}

Evidently the restriction of $Res^{m}(c_{\H}\{\underline{t}\})\in H_{\H}^{0}(E^{m}=\amalg\text{{pts.}},\QQ(0))\cong\oplus\QQ$
to $p\in\cap_{j}|(t_{j})|$ is $\pm1$ ($\neq0$, which puts $\{\underline{t}\}$
in the highest weight $W_{0}$). Since $\sZ$ is Res-free (hence in
$W_{-n}$), the product gives a class in \[
Gr_{-n}^{W}Gr_{\L}^{i+m}CH^{p+m}(X\times\{(\ms\times\T)\m\mathfrak{{D}}\},\, n+m)\]
\[
\hookrightarrow Gr_{-n}^{W}Gr_{\L}^{i+m}H_{\H}^{2p+m-n}\left(X\times\{(\ms\times\T)\m\mathfrak{{D}}\},\QQ(p+m)\right)\]
\[
\begin{array}[b]{c}
_{Res^{m}}\\
\longrightarrow\end{array}Gr_{\L}^{i}\underline{H_{\H}^{2p-n}}(X\times\mathfrak{{D}}^{m}\m\mathfrak{{D}}^{m+1},\QQ(p))\begin{array}[b]{c}
_{|_{(\ms\m D)\times\{ p\}}}\\
\longrightarrow\end{array}\left\{ \text{{r.h.s.\, of\,(5.4)}}\right\} .\]
Clearly $Res_{(\ms\m D)\times\{ p\}}^{m}\left\{ Gr_{\L}^{i+m}c_{\H}(\sZ\times\{\underline{t}\})\right\} =Gr_{\L}^{i}c_{\H}(\sZ)\otimes Res_{\{ p\}}^{m}(c_{\H}\{\underline{t}\})$$\linebreak$$\neq0$;
and it only remains to be said that this does not vanish upon enlarging
$\mathfrak{{D}}$ provided we enlarge $D$ so that $D\times\{ p\}\supseteq(\ms\times\{ p\})\cap\mathfrak{{D}}$.
\end{proof}

\section{\textbf{Detecting indecomposables in the \\ regulator kernel}}

The $j$-decomposable cycles are those in the image of the exterior
product\[
\theta^{j}:\, CH^{p-j}(X_{K},m-j)\otimes CH^{j}(\text{{Spec}}K,\, j)\to\L^{j}CH^{p}(X_{K},m).\]
For any subgroup $\Xi\subseteq CH^{p}(X_{K},m)$, define\[
\Xi_{j\text{{-ind}}}:=\frac{\Xi}{im(\theta^{j})\cap\Xi}.\]
The usual definition of indecomposable is $\Xi_{m\text{{-ind}}}$,
while the {}``strong indecomposable'' notion of \cite{CF} corresponds
to $\Xi_{1\text{{-ind}}}$. We study here the case $\Xi=\L^{i}\underline{CH^{p}}(X_{K},m)$
--- these are cycles with complete spread; that is, their spreads
have no residues along substrata of $\ms$.

\subsection{An invariant for cycles with complete spread}

Our aim is to produce, starting from the higher $\{ AJ_{i-1}^{tr}\}$
maps of $\S5.1$, further quotients which detect indecomposables.
The principal difficulty is that an element of $im(\theta^{j})\cap\Xi$
(decomposables with complete spread) may involve terms whose spreads
have residues (although the sum of these residues must vanish).

\begin{thm}
(i) For $1\leq j\leq m\leq i\leq p$, $Gr_{\L}^{i}\underline{c_{\H}}$
(composed with a projection) descends to a map\[
\L^{i}\underline{CH}_{j\text{{-ind}}}^{p}(X_{K},m)\to Ext_{MHS}^{1}\left(\QQ(-p),\frac{H^{2p-m-i}(X)}{F_{h}^{p-i}}\otimes\underline{H^{i-1}}(\gs)\right).\]
 (ii) If $1\leq j$ and $i\leq p$, but one of the other inequalities
fails, then $F_{h}^{p-i}$ may be replaced by $F_{h}^{p-i+1}$ (i.e.,
detecting indecomposables is a trivial matter). If $j=0$ or $i>p$,
there is nothing to detect.
\end{thm}
\begin{proof}
From $\S1$ we know the map from $\L^{i}\underline{CH^{p}}(X_{K},m)$
to$\linebreak$ $Ext_{MHS}^{1}\left(\QQ(-p),\frac{H^{2p-m-i}(X)}{F_{h}^{p-i+1}}\otimes\underline{H^{i-1}}(\gs)\right)$
is well-defined. So given a class $\xi$ in $\L^{i}\underline{H_{\H}^{2p-m}}(X\times\gs,\QQ(p))$
$and$ in the image of

\begin{equation} \\
\begin{matrix} H_{\H}^{2p-m-j}(X\times \gs,\QQ(p-j))\otimes H^j_{\H}(\gs,\QQ(j)) \longrightarrow \\ \\ H^{2p-m}_{\H}(X\times \gs,\QQ(p)), \end{matrix} \\
\end{equation}\\
we need to show its projection $\bar{\xi}$ to \[
Gr_{\L}^{i}H_{\H}^{2p-m}\cong Ext_{MHS}^{1}\left(\QQ(-p),H^{2p-m-i}(X)\otimes H^{i-1}(\gs)\right)\]
 is contained in the image of 

\begin{equation} \\
\begin{matrix} Ext^1_{MHS}(\QQ(-p),F^{p-i}_hH^{2p-m-i}(X) \otimes \underline{H^{i-1}}(\gs ) \longrightarrow \\ \\

Ext^1_{MHS}(\QQ(-p),H^{2p-m-i}(X)\otimes H^{i-1}(\gs)). 

\end{matrix} \\
\end{equation}\\
Since (6.1) $splits$ into a $\oplus$ of Leray-graded pieces, we
have immediately that $\bar{\xi}\in$\[
\begin{array}{c}
Gr_{\L}^{i}\underline{H_{\H}^{2p-m}}\left(X\times\gs,\QQ(p)\right)\,\,\,\bigcap\\
\\im\left(\begin{array}{c}
Gr_{\L}^{i-j}H_{\H}^{2p-m-j}(X\times\gs,\QQ(p))\otimes H_{\H}^{j}(\gs,\QQ(j)\\
\to Gr_{\L}^{i}H_{\H}^{2p-m}(X\times\gs,\QQ(p))\end{array}\right)\end{array}\,\,\,\textbf{{=}}\]
\\

\small \begin{equation} \\
\begin{matrix} { im \left( \begin{matrix} Ext^i_{^{D^bMHS}} ( \QQ(-p), H^{2p-m-i}(X) \otimes \underline{\K_{\H}^{\bb}} (\gs) ) \rTo^{\epsilon^j_{\H}}  \\ Ext^i_{^{D^bMHS}} ( \QQ(-p), H^{2p-m-i}(X) \otimes \K_{\H}^{\bb}(\gs) ) \end{matrix} \right) \; \; \; \bigcap } \\ \\ { im \left( \begin{matrix} Ext^{i-j}_{^{D^bMHS}}(\QQ(-p+j),H^{2p-m-i}(X)\otimes \K_{\H}^{\bb}(\gs)) \otimes Ext^j_{^{D^bMHS}}(\QQ(-j),\K_{\H}^{\bb}(\gs)) \\ \rTo^{\theta_{\H}^j} Ext^i_{^{D^bMHS}}(\QQ(-p),H^{2p-m-i}(X)\otimes \K_{\H}^{\bb}(\gs)) \end{matrix} \right) } \end{matrix} .
\end{equation} \normalsize\\

Here we simply take \[
\K_{\H}^{\bb}(\gs):=\bigoplus_{\ell}\underline{H^{\ell}}(\gs,\QQ)[-\ell]\,,\,\,\,\underline{\K_{\H}^{\bb}}(\gs):=\bigoplus_{\ell}\underline{H^{\ell}}(\gs,\QQ)[-\ell].\]
 Moreover, using an idea of \cite{RS}, $Hom_{MHS}(\QQ(-j),H^{j}(\gs))$
gives a subMHS \[
H^{j}(\gs)''\subseteq H^{j}(\gs)\]
 with $H^{j}(\gs)''\cap W_{2j-1}H^{j}(\gs)=\{0\}$. By semisimplicity,
writing$\linebreak$ $\eta:\, H^{j}(\gs)''\to Gr_{2j}^{W}H^{j}(\gs)$,
we have $Gr_{2j}^{W}H^{j}(\gs)\begin{array}[b]{c}
_{HS}\\
=\end{array}\H'\oplus im(\eta).$ Pulling back the extension $0\to W_{2j-1}H^{j}(\gs)\to H^{j}(\gs)\to Gr_{2j}^{W}\to0$
along $\H'\hookrightarrow Gr_{2j}^{W}$ yields $H^{j}(\gs)'\hookrightarrow H^{j}(\gs)$
with $H^{j}(\gs)\begin{array}[b]{c}
_{MHS}\\
=\end{array}H^{j}(\gs)'\oplus H^{j}(\gs)''.$ Define more generally\[
H^{\ell}(\gs)'':=im\left\{ H^{j}(\gs)''\otimes H^{\ell-j}(\gs)\to H^{\ell}(\gs)\right\} ,\]
$\K_{\H}^{\bb}(\gs)'':=\oplus_{\ell}H^{\ell}(\gs)''[-\ell].$ A generalization
of the above argument (involving both projections and pullbacks of
extensions) produces ($\forall\ell$) $\oplus$-decompositions of
$H^{\ell}(\gs)$, hence $\K_{\H}^{\bb}(\gs)=\K_{\H}^{\bb}(\gs)'\oplus\K_{\H}^{\bb}(\gs)''$.
This induces a $\oplus$-decomposition on relevant $Ext$ groups (hence
also projections and inclusions of $Ext$'s). Note that $H^{\ell}(\gs)'\supseteq W_{\ell+j-1}H^{\ell}(\gs)\supseteq W_{\ell}H^{\ell}(\gs)=\underline{H^{\ell}}(\gs).$

Since $\K_{\H}^{\bb}(\gs)'\supseteq\underline{\K_{\H}^{\bb}}(\gs)$,
and the projection $\K_{\H}^{\bb}(\gs)\twoheadrightarrow\K_{\H}^{\bb}(\gs)'$
is compatible with $\cup$-products, (6.3) reduces to

\small \begin{equation} \\
\begin{matrix} { im \left( \begin{matrix} Ext^i_{^{D^bMHS}} ( \QQ(-p), H^{2p-m-i}(X) \otimes \underline{\K_{\H}^{\bb}} (\gs) ) \rTo^{(\epsilon^j_{\H})'}  \\ Ext^i_{^{D^bMHS}} ( \QQ(-p), H^{2p-m-i}(X) \otimes \K_{\H}^{\bb}(\gs)' ) \end{matrix} \right) \; \; \; \bigcap } \\ \\ { im \left( \begin{matrix} Ext^{i-j}_{^{D^bMHS}}(\QQ(-p+j),H^{2p-m-i}(X)\otimes \K_{\H}^{\bb}(\gs)') \otimes Ext^j_{^{D^bMHS}}(\QQ(-j),\K_{\H}^{\bb}(\gs)') \\ \rTo^{(\theta_{\H}^j)'} Ext^i_{^{D^bMHS}}(\QQ(-p),H^{2p-m-i}(X)\otimes \K_{\H}^{\bb}(\gs)') \end{matrix} \right) } \end{matrix} .
\end{equation} \normalsize\\
\\
By construction, $Hom_{MHS}(\QQ(-j),H^{j}(\gs)')=0$. Hence $(\theta_{\H}^{j})'$
reduces to

\small \begin{equation} \\
\begin{matrix} Hom_{MHS} \left( \QQ(-p+j),H^{2p-m-i}(X)\otimes H^{i-j}(\gs)' \right) \otimes Ext^1_{MHS} \left( \QQ(-j), H^{j-1}(\gs) \right) \\ \\ \longrightarrow Ext^1_{MHS}\left( \QQ(-p), H^{2p-m-i}(X) \otimes H^{i-1}(\gs)' \right) . \end{matrix} \\
\end{equation} \normalsize\\
One easily shows that (in (6.5)) the $Hom=Hom_{MHS}(\QQ(-p+j),$$\linebreak$$F_{h}^{p-i}H^{2p-m-i}(X)\otimes H^{i-j}(\gs)').$
(This is straightforward and has nothing to do with $H^{i-j}(\gs)$
vs. $H^{i-j}(\gs)'$.) Hence (6.5) factors through\[
\begin{array}{c}
Ext_{MHS}^{1}\left(\QQ(-p),F_{h}^{p-i}H^{2p-m-i}(X)\otimes H^{i-1}(\gs)'\right)\\
\\\hookrightarrow\,\, Ext_{MHS}^{1}\left(\QQ(-p),H^{2p-m-i}(X)\otimes H^{i-1}(\gs)'\right)\end{array}\]
and so (6.4) is contained in\[
\begin{array}{c}
im\left(\begin{array}{c}
Ext_{MHS}^{1}\left(\QQ(-p),H^{2p-m-i}(X)\otimes\underline{H^{i-1}}(\gs)\right)\\
\to Ext_{MHS}^{1}\left(\QQ(-p),H^{2p-m-i}(X)\otimes H^{i-1}(\gs)'\right)\end{array}\right)\\
\\\bigcap\,\,\,\, Ext_{MHS}^{1}\left(\QQ(-p),F_{h}^{p-i}H^{2p-m-i}(X)\otimes H^{i-1}(\gs)'\right)\end{array}.\]
Using semisimplicity to write $H^{2p-m-i}(X)\cong F_{h}^{p-i}\oplus\frac{H^{2p-m-i}(X)}{F_{h}^{p-i}}$
as HS, we see that the above lies in\[
im\left\{ \begin{array}{c}
Ext_{MHS}^{1}\left(\QQ(-p),F_{h}^{p-i}H^{2p-m-i}(X)\otimes\underline{H^{i-1}}(\gs)\right)\\
\to Ext_{MHS}^{1}\left(\QQ(-p),F_{h}^{p-i}H^{2p-m-i}(X)\otimes H^{i-1}(\gs)'\right)\end{array}\right\} ,\]
which obviously $\subseteq$ $im\,(6.2)$.\\
\\
(ii) If $j>i$, the $j$-decomposables lie in $\L^{i+1}$; if $j>m$,
there aren't any; and if $m>i$, the $Hom$ in (6.5) vanishes.
\end{proof}
Setting $i=j\,:=1$ and $p:=q$ in Theorem $6.1$, we have for $m\neq1$
(using (ii)) \[
AJ_{\text{{ind}}}:\,\L^{1}\underline{CH}_{1\text{{-ind}}}^{q}(X_{k},m)\to Ext_{MHS}^{1}\left(\QQ(-q),H^{2p-m-1}(X)\right)=J^{p,m}(X)\]
(since $F_{h}^{p}H^{2p-m-1}(X)=\{0\}$) and for $m=1$ (using (i))\[
AJ_{\text{{ind}}}:\,\L^{1}\underline{CH}_{1\text{{-ind}}}^{q}(X_{k},1)\to Ext_{MHS}^{1}\left(\QQ(-q),\frac{H^{2p-2}(X)}{Hg^{p-1}(X)}\right)\]
(since $F_{h}^{p-1}H^{2p-2}(X)\cong Hg^{p-1}(X)$). It seems natural
to call a cycle $\W\in CH_{hom}^{q}(X_{k},m)$ with $AJ_{\text{{ind}}}(\W)\neq0$
{}``strongly $AJ$-indecomposable'', but evidently (from the form
of $AJ_{\text{{ind}}}$ in the $2$ cases) this is the same thing
as {}``regulator-indecomposable'' (as in \cite{RS}, \cite{CF},
\cite{Le2}, etc.), so we will use that terminology.

\subsection{First application: exterior products}

The real point of Theorem 6.1, however, is to detect indecomposable
higher cycles in the kernel of $AJ$. Here is a simple application
to classes in the image of the exterior product\[
\L^{i}CH^{d_{1}}((Y_{1})_{K})\otimes\L^{1}CH^{p-d_{1}}((Y_{2})_{k},m)\to\L^{i+1}\underline{CH^{p}}((Y_{1}\times Y_{2})_{K},m),\]
where $K/k$ is finitely generated of transcendence degree $i$, and
$Y_{1},Y_{2}$ are smooth projective $/k$ of resp. dimensions $d_{1},d_{2}$.

\begin{thm}
Let $\W\in CH_{hom}^{p-d_{1}}((Y_{2})_{k},m)$ be regulator-indecomposable,
and $\V\in\L^{i}CH^{d_{1}}((Y_{1})_{K})$ a $0$-cycle whose complete
$k$-spread $\bar{\sV}\in Z^{d_{1}}((\ms\times Y_{1})_{k})$ induces
a nontrivial map $\Omega^{i}(Y_{1})\to\Omega^{i}(\ms)$ of holomorphic
$i$-forms. (See \cite{Ke3} for how to construct such $0$-cycles.)

Then $\Z:=\V\times\W\in CH^{p}((Y_{1}\times Y_{2})_{K},m)$ is strongly
indecomposable.
\end{thm}
\begin{proof}
This is so similar to the proof of Thm. 1 in \cite{Ke3} that we merely
indicate briefly how $\Z$ is detected and the technical lemma on
HS's needed.

Applying the map of Thm. 6.1 in the appropriate case (substituting
the present $i+1$ for $i$ and $1$ for $j$) sends $\Z$ to a quotient
(in general) of $[AJ(\sZ)]_{i}^{tr}$. Since $\Z$ has complete spread
$\bar{\sZ}:=\W\times\bar{\sV}$, this can be computed (as in $\S1$)
using $(2\pi i)^{p-m}R_{\bar{\sZ}}'$ (as a functional on forms).
Now $T_{\bar{\sZ}}=T_{\W}\times\bar{\sV}$ $\implies$ $\partial^{-1}T_{\bar{\sZ}}:=(\partial^{-1}T_{\W})\times\bar{\sV},$
and $R_{\bar{\sZ}}=R_{\W}\cdot(2\pi i)^{d_{1}}\delta_{\bar{\sV}}$;
so $\frac{1}{(2\pi i)^{d_{1}}}R_{\bar{\sZ}}'=R_{\W}'\cdot\delta_{\bar{\sV}}$,
and (the quotient of) $[AJ(\sZ)]_{i}^{tr}$ may be computed by the
image of $[\bar{\sV}]_{i}\otimes AJ_{\text{{(ind)}}}(\W)$. (We are
assuming $AJ_{\text{{ind}}}(\W)\neq0$.)

Define sub-Hodge structures\[
\H_{1}:=H^{2d_{1}-i}(Y_{1})\otimes H^{2p-2d_{1}-m-1}(Y_{2})\otimes H^{i}(\ms)\]
\[
\G_{0}:=H^{2p-i-m-1}(Y_{1}\times Y_{2})\otimes N^{1}H^{i}(\ms)\,+\, F_{h}^{p-i-1}H^{2p-i-m-1}(Y_{1}\times Y_{2})\otimes H^{i}(\ms)\]
of $\H_{0}:=H^{2p-i-m-1}(Y_{1}\times Y_{2})\otimes H^{i}(\ms)$, and\[
\H_{\bar{\sV}}:=\QQ\left\langle [\bar{\sV}]_{i}\right\rangle \otimes H^{2p-2d_{1}-m-1}(Y_{2})\]
\[
\G_{1}:=H^{2d_{1}-i}(Y_{1})\otimes F_{h}^{p-d_{1}-1}H^{2p-2d_{1}-m-1}(Y_{2})\otimes H^{i}(\ms)\]
sub-HS of $\H_{1}$. From the last paragraph,

\begin{equation} \\
\left[\bar{\sV}\right]_i \otimes AJ_{\text{ind}}(\W) \; \in \; \QQ \left[ \bar{\sV}\right]_i \otimes J^{p,m}_{\text{ind}}(Y_2) \; \cong \; Ext^1_{MHS} \left( \QQ(-p), \frac{\H_{\bar{\sV}}}{\G_1 \cap \H_{\bar{\sV}}} \right) \\ 
\end{equation}\\
and the Thm. 6.1 quotient of $[AJ(\sZ)]_{i}^{tr}$ in $Ext_{MHS}^{1}(\QQ(-p),\H_{0}/\G_{0})$,
map to the same element of

\begin{equation} \\
\frac{Ext^1_{MHS} \left( \QQ(-p), \H_1/\G_1 \right)}{Ext^1_{MHS}\left( \QQ(-p),(\H_1\cap\G_0)/(\G_1\cap \G_0)\right)} . \\
\end{equation}\\
One shows that \[
\H_{\bar{\sV}}\cap(\H_{1}\cap\G_{0})\subseteq\G_{1},\]
so that the map from (6.6) to (6.7) is an injection and we are done.
\end{proof}

\subsection{Relation to a criterion of Lewis}

Given $\ms'$ smooth projective $/k$ and $\alpha\in H^{\dim(\ms')-i,\,\dim(\ms')}(\ms',\CC),$
let $cl_{\ms',\alpha}^{\{ p-m,i-m,0\}}$ denote the composition\\
\small \xymatrix{CH^{p-i}(X\times \ms') \ar [r]^{cl} & H^{2p-2i}(X\times \ms',\CC) \ar @{->>} [r]^{\text{Hodge-}\mspace{40mu}}_{\text{K\"unneth} \mspace{40mu}} & H^{p-i,p-m}(X,\CC) \otimes H^{i,0}(\ms',\CC) \ar [d]^{\cup \alpha} \\ & & H^{p-i,p-m}(X,\CC). } \normalsize\\

The next application (Theorem 6.3 below) strengthens and generalizes
a result of Collino and Fakhruddin \cite{CF}. To put this in context,
we first explain how Theorem 6.1 essentially subsumes the invariant
they use for indecomposables in $\ker(AJ)$, due to \cite{Le2}. This
may be expressed as a map \begin{equation} \L^i \underline{CH}^p_{m\text{-ind}}(X_K,m) \to \frac{H^{p-i,p-m}(X,\CC)}{H^{\{p-m,i-m,0\}}(X)} \otimes H^{i-1,0}(\ms,\CC) \end{equation} 
where as usual $k(\ms)\cong K$, $m\geq1$, and\[
H^{\{ p-m,i-m,0\}}(X):=\sum_{\ms',\alpha}\text{{im}}\left(cl_{\ms',\alpha}^{\{ p-m,i-m,0\}}\right).\]
To describe (6.8), let $H$ and $G$ be quotient Hodge structures
(resp.) of $H^{2p-m-i}(X,\QQ)$ and $H^{i-1}(\ms,\QQ)$, such that
$G_{\CC}^{i-1,0}\cong H^{i-1,0}(\ms,\CC)$; let $\H:=H\otimes G$
and write $pr_{\H}:H^{2p-m-1}(X\times\ms,\QQ)\rOnto\H$. Consider
the composition $\R_{\H}^{\{ p,i,m\}}:$\\
\small \xymatrix{\ext (\QQ(-p),\H) \ar @{->>} [d] & CH^p(X\times \ms ,m) \ar [l]_{pr_{\H}\circ AJ} \\  \text{Ext}^1_{_{\RR\text{-MHS}}} (\RR(-p),\H_{\RR}) \ar @{=} [r] & \H_{\RR(p-1)}\cap \{\H^{p-m,p-1}_{\CC} \oplus \cdots \oplus \H^{p-1,p-m}_{\CC}\} \ar [d] \\ H^{p-i,p-m}_{\CC} \otimes H^{i-1,0}(\ms ,\CC) &  \H^{p-1,p-m}_{\CC} \ar @{->>} [l] \\ } \normalsize \\
\\
Let $\H_{1}:=H^{2p-m-i}(X)\otimes H^{i-1}(\ms)$, $\H_{2}:=\frac{H^{2p-m-i}(X)}{F_{h}^{p-i}}\otimes\underline{H}^{i-1}(\gs)$
be choices of $\H$ as above, and consider the commutative diagram\\
\xymatrix{ CH^p(X \times \ms ,m) \ar [r]_{\R^{\{p,i,m\}}_{\H_1}\mspace{50mu}} \ar [rd]_{\overline{\R^{\{p,i,m\}}_{\H_1}}\mspace{0mu}} \ar [rdd]_{\R^{\{p,i,m\}}_{\H_2}\mspace{0mu}} \ar [dd]_{\overline{pr_{\H_2}}\circ AJ} & H^{p-i,p-m}(X,\CC) \otimes H^{i-1,0}(\ms,\CC) \ar @{->>} [d] \\ & {\frac{H^{p-i,p-m}(X,\CC)}{H^{\{p-m,i-m,0\}}(X)}} \otimes H^{i-1,0}(\ms,\CC) \ar @{->>} [d]^{(*)} \\ {\ext \left( \QQ(-p), \H_2 \right) \ar [r]_{(**)\mspace{100mu}}} & {\left( \frac{H^{2p-m-i}(X)}{F^{p-i}_h} \right) ^{(p-i,p-m)}_{\CC} \otimes H^{i-1,0}(\ms,\CC).} \\ }  \\
\\
The proof in \cite{Le2} implies that $\overline{\R_{\H_{1}}^{\{ p,i,m\}}}$
is well-defined on the subquotient $\L^{i}\underline{CH}_{m\text{{-ind}}}^{p}(X_{K},m)$
of $CH^{p}(X\times\ms,m)$, and this gives (6.8). Moreover, $\overline{{pr_{\H_{2}}}}\circ AJ$
is just the map of our Theorem $6.1$ (so this is also well-defined
on $\L^{i}\underline{CH}_{m\text{{-ind}}}^{p}$). Now the GHC says
that $F_{h}^{p-i}H^{2p-m-i}(X)=N^{p-i}H^{2p-m-i}(X)$; and Lewis shows
that the Hard Lefschetz Conjecture (implied by GHC) gives $\frac{H^{p-i,p-m}(X,\CC)}{H^{\{ p-m,i-m,0\}}(X)}\cong\left(\frac{H^{2p-m-i}(X)}{N^{p-i}}\right)_{\CC}^{(p-i,p-m)}.$
Hence GHC $\implies$ $(*)$ is an isomorphism, which means that $\overline{pr_{\H_{2}}}\circ AJ$
factors $\overline{\R_{\H_{1}}^{\{ p,i,m\}}}$. The point is that
we conjecturally aren't losing any information via $(*)$, while we
definitely are losing data along $(**)$.

\subsection{Second application: Pontryagin cycles}

We now describe the improvement of \cite[Prop. 1.1]{CF}. Let $C$
be a smooth projective curve $/k$ of genus $g$. Fix $i\in\mathbb{{N}}$;
$q_{1},\ldots,q_{i}\in\C(k)$; $o\in J(\C)(k)$ {[}origin of group
law{]}; and $P:=(p_{1},\ldots,p_{i})\in\C^{\times i}(\CC)$ very general,
i.e. $K:=k(P)\subset\CC$ has $trdeg(K/k)=i$. We define \[
\tilde{B}_{i}:\, Z_{\RR}^{p}(J(C),m)\to Z_{\RR}^{p}(J(C)_{K},m)\]
by $\Z\mapsto\Z*\left(AJ(p_{1}-q_{1})-o\right)*\cdots*\left(AJ(p_{i}-q_{i})-o\right)=:\Z(i).$
This induces\[
B_{i}:\, CH_{hom}^{p}(J(C),m)\to\L^{i+1}\underline{CH}^{p}(J(C)_{K},m)\]
 since the $Pontryagin$ $cycle$ $\Z(i)$ is in $\L^{i+1}$ by the
parallelogram law \cite{Bl3}, \cite[3.2]{Ke2}.

\begin{thm}
Assume \begin{equation} i\leq 2p-g-m-1, \end{equation}  and let $\left\langle \Z\right\rangle \in CH_{hom}^{p}(J(C),m)$
be given.\\
(a) If $AJ(\Z)\neq0,$ then $AJ_{i}(\Z(i))\neq0,$ hence\[
0\neq\left\langle \Z(i)\right\rangle \in Gr_{\L}^{i+1}\underline{CH}^{p}(J(C)_{K},m).\]
\\
(b) Let $m\geq1$. If $\Z$ is regulator-indecomposable, then\[
0\neq\overline{\left\langle \Z(i)\right\rangle }\in Gr_{\L}^{i+1}\underline{CH}_{1\text{{-ind}}}^{p}(J(C)_{K},m);\]
in particular $\Z(i)$ is (strongly) indecomposable.
\end{thm}
\begin{rem}
The case corresponding to \cite{CF} is (b) $m=1$. No assumption
on the sub-Hodge structures of the $H^{i}(J(C),\QQ)$ is required.
\end{rem}
We establish this with a parade of definitions and lemmas. Define
cycles $/k$\small \[
\mathfrak{A}_{i}:=\sum_{j=0}^{i}(-1)^{j}\sum_{|J|=j}\left\{ (u;c_{1},\ldots,c_{i})\,\left|\, u=\sum_{\ell\in J}AJ_{C}(c_{\ell}-q_{\ell})\right.\right\} \in Z^{g}(J(C)\times C^{\times i}),\]
\normalsize \[
\mathfrak{G}:=\left\{ (u_{1},u_{2},u_{3})\,|\, u_{1}+u_{2}=u_{3}\right\} \in Z^{2g}(J(C)^{\times3}).\]
Let $\pi_{I},\, I\subseteq\{1,2,3,4\}$, denote various projections
from $J(C)^{\times3}\times C^{\times i}$ (where factors $1,2,3$
are copies of $J(C)$, and $C^{\times i}$ is the $4^{th}$ factor).
In the otherwise self-explanatory diagram\\
\xymatrix{ & H^{2p-m-1}(J(C)\times C^{\times i}) \ar @{->>} [d]^{\text{K\"unneth}} \\ H^{2p-m-1}(J(C)) \ar [ru]^{\beta_i} \ar [r]^{\bar{\beta}_i \mspace{50mu}} \ar [rd]^{\bar{\beta}_i'} \ar @{->>} [d]^{\mathsf{p}} & H^{2p-m-i-1}(J(C))\otimes H^i(C^{\times i}) \ar @{=} [r] \ar @{->>} [d]^{\mathsf{p}'} & :\H \\ {\frac{H^{2p-m-1}(J(C))}{F^{p-1}_h}} \ar [rd]^{\bar{\beta}_i''} & {\frac{H^{2p-m-i-1}(J(C))}{F^{p-i}_h}}\otimes \underline{H}^i(\eta_{_{C^i}}) \ar @{->>} [d]^{\mathsf{p}''} \\  & {\frac{H^{2p-m-i-1}(J(C))}{F^{p-i-1}_h}} \otimes \underline{H}^i(\eta_{_{C^i}}) \\ }
\\
define\[
\beta_{i}(\xi):=\pi_{34}{}_{_{*}}\{\pi_{123}^{^{*}}[\mathfrak{G}]\cup\pi_{24}^{^{*}}[\mathfrak{A}_{i}]\cup\pi_{1}^{^{*}}\xi\}.\]
Since $\beta_{i}$ is induced by action of a correspondence, all arrows
are morphisms of (M)HS. (Note that well-definedness of $\bar{\beta}_{i}''$
is clear once one knows $\bar{\beta}_{i}'$ is a morphism of HS.)
We also write \\
\xymatrix{H^{2p-m-1}(J(C),\CC) \ar [rd]^{\bar{\beta}_i} \ar [r]^{\bar{\beta}_i^{(i,0)} \mspace{50mu}} & H^{2p-m-i-1}(J(C),\CC)\otimes H^{i,0}(C^{\times i},\CC) \\ & H^{2p-m-i-1}(J(C),\CC)\otimes H^i(C^{\times i},\CC) . \ar @{->>} [u]_{pr^{(i,0)}} }\\

\begin{lem}
If (6.9) holds, $\bar{\beta}_{i}^{(i,0)}$ and $\bar{\beta}_{i}$
are injective.
\end{lem}
\begin{proof}
It suffices to examine $\bar{\beta}_{i}^{(i,0)}$. Note that pullback
along $\{ p\mapsto AJ(p-q)\}:\, C\to J(C)$ is independent of $q\in C$,
and call this $AJ^{*}$. Let $\{\Omega_{\ell}\}_{\ell=1}^{g}\subset\Omega^{1}(J(C))$
be a basis such that $\{\omega_{\ell}:=AJ^{*}\Omega_{\ell}\}\subset\Omega_{C}^{1}$
satisfies $\int_{C}\omega_{k}\wedge\overline{\omega_{\ell}}=\delta_{k\ell}$.
Write $\Omega_{k}:=\bigwedge_{k\in K}\Omega_{k}$, and let $\tau_{IJ}:=\Omega_{I}\wedge\overline{\Omega_{J}}$
for any pair of multi-indices $I,J\subseteq\{1,\ldots,g\}$ with $|I|+|J|=2g-2p+m+1$.

Given nonzero $\xi\in H^{2p-m-1}(J(C))$, $\exists$ $I,J$ such that
$\left\langle \xi,\tau_{IJ}\right\rangle \neq0$. Since (6.9) holds,
we may take $K=\{ k_{1},\ldots,k_{i}\}\subseteq\{1,\ldots,g\}\m I\cup J$.
Write $\pi_{0}:\, J(C)\times C^{\times i}\rOnto J(C)$, $\pi_{j}:\, J(C)\times C^{\times i}\rOnto C\,\,\,(j=1,\ldots i)$
for the obvious projections. By an easy calculation (essentially the
same as that in \cite[Lemma 1.3]{CF}),\[
\left\langle \beta_{i}(\xi),\,\pi_{0}^{^{*}}(\Omega_{K}\wedge\tau_{IJ})\wedge\pi_{1}^{^{*}}(\overline{\omega_{k_{1}}})\wedge\cdots\wedge\pi_{i}^{^{*}}(\overline{\omega_{k_{i}}})\right\rangle =\pm\left\langle \xi,\tau_{IJ}\right\rangle ;\]
hence $\bar{\beta}_{i}^{(i,0)}(\xi)\neq0$.
\end{proof}
For a map of (pure $\QQ$-)Hodge structures $\H_{1}\rTo^{\theta}\H_{2}$
of weight $<2p$, write $Ext^{1}(\theta)$ for the induced map \small \[
\ext(\QQ(-p),\H_{1})\cong\frac{\H_{1,\CC}}{F^{p}\H_{1,\CC}+\H_{1}}\rTo^{\bar{\theta}}\frac{\H_{2,\CC}}{F^{p}\H_{2,\CC}+\H_{2}}\cong\ext(\QQ(-p),\H_{2}).\]
\normalsize ${}$

\begin{lem}
The following diagram \SMALL \\
\xymatrix{ CH^p_{hom}(J(C),m) \ar [r]_{B_i} \ar [d]^{AJ} & {\L^{i+1}} \underline{CH}^p(J(C)_K,m) \ar [d]^{AJ^{tr}_i} \\ \ext (\QQ(-p),H^{2p-m-1}(J(C))) \ar @{->>} [d]^{Ext^1(\mathsf{p})} \ar [r]_{Ext^1(\bar{\beta}_i') \mspace{50mu}} & \ext \left( \QQ(-p),\frac{H^{2p-m-i-1}(J(C))}{F^{p-i}_h}\otimes \underline{H}^i (\eta_{_{C^i}})\right) \ar @{->>} [d]^{Ext^1(\mathsf{p}'')} \\ \ext \left( \QQ(-p) , \frac{H^{2p-m-1}(J(C))}{F^{p-1}_h} \right) \ar [r]_{Ext^1(\bar{\beta}_i'')\mspace{50mu}} & \ext \left( \QQ(-p), \frac{H^{2p-m-i-1}(J(C))}{F^{p-i-1}_h}\otimes \underline{H}^i(\eta_{_{C^i}}) \right) \\ } \\
\normalsize \\
is commutative.
\end{lem}
\begin{proof}
Let $\iota_{P}:\, J(C)[\times\{ P\}]\rInto J(C)\times C^{\times i}$
be the inclusion (defined $/K$). Define a map $\tilde{\mathfrak{B}}_{i}$
factoring $\tilde{B}_{i}$\\
\xymatrix{& & Z^p_{\RR}(J(C),m) \ar [rd]_{\tilde{\mathfrak{B}}_i} \ar [r]^{\tilde{B}_i} & Z^p_{\RR}(J(C)_K,m) \\ & & & Z^p_{\RR}(J(C)\times C^{\times i},m) \ar [u]_{\iota_P^*} }\\
\\
by sending $\Z\mapsto\pi_{34}^{^{*}}\{\pi_{123}^{^{*}}\mathfrak{G}\cdot\pi_{24}^{^{*}}\mathfrak{A}_{i}\cdot\pi_{1}^{^{*}}\Z\}=:\overline{\sZ(i)}$.
Clearly $\tilde{\mathfrak{B}}_{i}$ descends modulo $\rateq$. Since
$\overline{\sZ(i)}$ is a complete spread for $\Z(i)$, $AJ_{i}^{tr}(\Z(i))$
is a projection of $AJ_{J(C)\times C^{i}}(\overline{\sZ(i)})$ and
it suffices to show \small \\
\xymatrix{CH^p_{hom}(J(C),m) \ar [d]^{AJ} \ar [r]_{\mathfrak{B}_i} & CH^p_{hom}(J(C)\times C^{\times i},m) \ar [d]^{AJ} \\ \ext \left( \QQ(-p), H^{2p-m-1}(J(C)) \right) \ar [r]_{Ext^1(\beta_i)\mspace{50mu}} & \ext \left( \QQ(-p), H^{2p-m-1}(J(C)\times C^{\times i}) \right) }\\
\normalsize \\
commutes. Now $\overline{\sZ(i)}$ is a (signed) sum of families of
translates (over $C^{\times i}$) of $\Z$. Translating $\Z$ on $J(C)$
gives a new (higher) cycle whose KLM currents $\delta_{T},\,\Omega,\, R$
are simply push-forwards of $\delta_{T_{\Z}},\,\Omega_{\Z},\, R_{\Z}$
along the translation. The same thing happens with the closed current
$R':=R-(2\pi i)^{*}\delta_{\d^{-1}T}+d^{-1}\Omega$, whose class in
$H^{2p-m-1}(J(C),\CC)$ descends to $AJ_{J(C)}(\Z)$. More precisely,
we can choose\[
R_{\overline{\sZ(i)}}':=\pi_{34}{}_{_{*}}\{\pi_{123}^{^{*}}\delta_{\mathfrak{G}}\cdot\pi_{24}^{^{*}}\delta_{\mathfrak{A}_{i}}\cdot\pi_{1}^{^{*}}R_{\Z}'\},\]
 and taking classes gives the result.
\end{proof}
Consider the sub-Hodge structures\[
\N:=H^{2p-m-i-1}(J(C))\otimes N^{1}H^{i}(C^{\times i}),\]
\[
\F_{0}:=F_{h}^{p-i}H^{2p-m-i-1}(J(C))\otimes H^{i}(C^{\times i}),\]
\[
\F_{1}:=F_{h}^{p-i-1}H^{2p-m-i-1}(J(C))\otimes H^{i}(C^{\times i})\]
of $\H$; clearly $\mathsf{{p}}'$ and $\mathsf{{p}}''\circ\mathsf{{p}}'$
are (resp.) quotients of $\H$ by $\N+\F_{0}$, $\N+\F_{1}$.

\begin{lem}
Provided (6.9) holds, we have:\[
(\N+\F_{1})\cap\text{{im}}(\bar{\beta}_{i})\subseteq\bar{\beta}_{i}(F_{h}^{p-i}),\]
\[
(\N+\F_{0})\cap\text{{im}}(\bar{\beta}_{i})=\{0\}.\]

\end{lem}
\begin{proof}
Consider $\xi\in H^{2p-m-1}(J(C),\CC)$ with $\bar{\beta}_{i}(\xi)\in\N_{\CC}+\F_{\ell,\CC}$,
$\ell=0\text{{\, or\,}}1$. We have $pr^{(i,0)}(\N_{\CC})=\{0\}$,
$F_{h}^{p-i-\ell}H^{2p-m-i-1}(J(C),\CC)\subset F^{p-i-\ell}$, hence
\begin{equation} \bar{\beta}_i^{(i,0)}(\xi) \in F^{p-i-\ell}H^{2p-m-i-1}(J(C),\CC)\otimes H^{i,0}(C^{\times i},\CC) . \end{equation}
Since $\bar{\beta}_{i}$ is a morphism of HS, it decomposes ($\otimes\CC$)
by type; hence $\bar{\beta}_{i}^{(i,0)}$ must decompose into a direct
sum (over $a\in\ZZ$) of maps\[
H^{a,2p-m-a-1}(J(C),\CC)\to H^{a-i,2p-m-a-1}(J(C),\CC)\otimes H^{i,0}(C^{\times i},\CC).\]
These are injective (because $\bar{\beta}_{i}^{(i,0)}$ is), so (6.10)
$\implies$ $\xi\in F^{p-\ell}H^{2p-m-1}(J(C),\CC).$

Now $\bar{\beta}_{i}^{-1}(\N+\F_{\ell})$ is necessarily a subHS of
$H^{2p-m-1}(J(C))$. The last paragraph shows $\bar{\beta}_{i}^{-1}(\N_{\CC}+\F_{\ell,\CC})\subset F^{p-\ell}$
($\ell=0,1$); hence\[
\bar{\beta}_{i}^{-1}(\N+\F_{\ell})\subseteq F_{h}^{p-\ell}H^{2p-m-1}(J(C)).\]
Applying $\bar{\beta}_{i}$ to this gives the result, noting (for
$\ell=0$) that $F_{h}^{p}H^{2p-m-1}(J(C))=\{0\}$.
\end{proof}
${}$

\begin{proof}
(of Theorem 6.3). Since $\bar{\beta}_{i}$ is injective (Lemma 6.5),
Lemma 6.7 implies $\bar{\beta}_{i}',\,\bar{\beta}_{i}'',$ hence $Ext^{1}(\bar{\beta}_{i}'),\, Ext^{1}(\bar{\beta}_{i}''),$
are injective too. By assumption, (a) $AJ(\Z)$ resp. (b) $Ext^{1}(\mathsf{{p}})(AJ(\Z))$
is nonzero. So Lemma 6.6 $\implies$ nonvanishing of (a) $AJ_{i}^{tr}(\Z(i))$
resp. (b) $Ext^{1}(\mathsf{{p}}'')(AJ_{i}^{tr}(\Z(i)))$, as desired.
\end{proof}

\section{\textbf{Reduced higher Abel-Jacobi maps}}

We now give a description of part of each higher Abel-Jacobi invariant
in the spirit of Griffiths-Green \cite{GG}. As in $\S5.1$ we avoid
double complexes, but we don't restrict to the {}``nice case'' examined
there. 

Roughly speaking: reduced higher $AJ$ maps are, like the $AJ_{i}^{tr}$
we have been using, Hodge-theoretically defined quotients of the previously
defined $AJ_{i}$. Only we quotient somewhat farther, in order to
get something which is explicitly computable in terms of generalized
membrane integrals on $X$ --- even when one does not have a homologically
trivial complete spread. Our expectation is that they will be useful
in special cases of arithmetic interest (as have their analogues for
$0$-cycles, see \cite[sec. 16]{Ke1}, \cite[sec. 6.3]{Ke2}), rather
than for proving results of a general flavor (such as Theorems 6.2
and 6.3 above). Assume $m>0$ for the discussion.

Set $V=F^{d-p+i}H^{2d-2p+m+i}(X,\CC)$. The {}``period lattice''\[
\Lambda:=im\left\{ H^{2p-m-i}(X,\QQ(p))\to\frac{H^{2p-m-i}(X,\CC)}{F^{p-i+1}H^{2p-m-i}(X,\CC)}\cong V^{\vee}\right\} \]
has in general $\dim_{\QQ}\Lambda\neq2\dim_{\CC}V$. There are (well-defined)
maps

\begin{equation} \\
Gr^i_{\L} \underline{CH^p}(X_K,m) \rTo^{\overline{AJ_{i-1}}} Hom_{\QQ} \left( \overline{H_{i-1}}(\gs,\QQ), V^{\vee}/\Lambda \right) \\
\end{equation}\\
where $\overline{H_{i-1}}(\gs):=coim\{ H_{i-1}(\gs)\to H_{i-1}(\ms)\},$
and (for cycles without complete spread)

\begin{equation} \\
Gr^i_{\L} CH^p(X_K,m) \rTo^{\overline{cl_i}} Hom_{\CC} \left( V, F^iW_{m+i}H^i(\gs,\CC) \right), \\
\end{equation}\\
\begin{equation} \\
ker(\overline{cl_i}) \rTo^{\overline{AJ_{i-1}}} Hom_{\QQ} \left( H_{i-1}(\gs,\QQ), V^{\vee}/\Lambda \right) . \\
\end{equation}\\
We note that as $ker(\overline{cl_{i}})\supseteq ker(cl_{i})$ is
not in general an equality, $\overline{AJ_{i-1}}$ may be defined
{}``before'' $AJ_{i-1}$ (of $\S5.1$) is.

Whether or not $\Z\in\L^{i}$ has complete $k$-spread $\bar{\sZ}$
($\partial_{\B}$-closed in $Z^{p}(X\times\ms,m)$) does not matter
so much for the formulas for $\overline{AJ_{i-1}}$. In either case
one has (for $U\subset\ms$ {}``sufficiently small'' Zariski open)
$\sZ_{U}\in Z_{\RR}^{p}(X\times U,m)$, $\partial_{\B}$-closed, with
irreducible components equidimensional $/U$, and restricting to the
class of $\sZ=\mathfrak{{s}}(\Z)\in CH^{p}(X\times\gs,m)$. Given
$[\gamma]\in H_{i-1}(\gs,\QQ)$ resp. $\overline{H_{i-1}}(\gs,\QQ)$,
its image in $H_{i-1}(U)$ resp. $\overline{H_{i-1}}(U)$ may be represented
by $\gamma\in Z_{i-1}^{top}(U,\QQ)$ such that $\gamma\times X\times\square^{m}$
intersects properly all $\sZ_{U}\cap\{\ms\times X\times\partial_{f,\RR}^{I,\ell}\square^{m}\}$.
(See $\S8.2$ for a definition of $\partial_{f,\RR}^{I,\ell}\square^{m}$.)
Finally, we take classes $\in V$ to be represented by test forms
$\omega\in Z_{d}(X,F^{d-p+i}\Omega_{X^{\infty}}^{2d-2p+m+i}).$

Writing $\pi_{X}:\, X\times\ms\to X$, $\pi_{U}:X\times U\to U$,
define\[
T_{\sZ}(\gamma):=\pi_{X}((\gamma\times X)\cap T_{\sZ})\in Z_{2d-2p+m+i}^{top}(X,\QQ),\]
\[
R_{\sZ}(\gamma):=\pi_{X_{*}}(\delta_{\gamma\times X}\cdot R_{\sZ})\in\D^{2p-m-i}(X)\]
(i.e. $(R_{\sZ}(\gamma))(\omega):=\int_{\gamma\times X}\pi_{X}^{^{*}}\omega\wedge R_{\sZ}$
defined by appropriate $\lim_{\epsilon\to0}\int$'s), and \[
\Omega_{\sZ}(\omega):=\pi_{U_{*}}(\pi_{X}^{^{*}}\omega\wedge\Omega_{\sZ})\in\Gamma_{d-\text{{cl.}}}(U,\D_{U}^{i,0})=\Omega^{i}(U).\]
Note $T_{\sZ}(\gamma)\homeq0$ since $\sZ\in\L^{i}$ $\implies$ $[\sZ]_{i-1}=0$.
The following result says how to compute (7.1),(7.2),(7.3).

\begin{thm}
(i) Given $\left\langle \Z\right\rangle \in\L^{i}CH^{p}(X_{K},m),$
$[\omega]\in V$, $\overline{cl_{i}}\left\langle \Z\right\rangle [\omega]$
is computed by $(2\pi i)^{p-m}\Omega_{\sZ}(\omega)\in\Omega^{i}(U)$.\\
(ii) Given $\left\langle \Z\right\rangle \in\{ ker(\overline{cl_{i}})\subseteq\L^{i}CH^{p}(X_{K},m)\}$
resp. $\L^{i}\underline{CH^{p}}(X_{K},m),$ $[\gamma]\in H_{i-1}(\gs,\QQ)$
resp. $\overline{H_{i-1}}(\gs,\QQ),$ $\overline{AJ_{i-1}}\left\langle \Z\right\rangle [\gamma]$
is computed by \[
(2\pi i)^{p-m}R_{\sZ}(\gamma)+(2\pi i)^{p}\delta_{\partial^{-1}T_{\sZ}(\gamma)}=:(2\pi i)^{p-m}R_{\sZ}(\gamma)'\in\D^{2p-m-i}(X),\]
viewed as a functional on $V$.
\end{thm}
\begin{proof}
In the complete spread case (7.1), the Theorem essentially follows
from writing $\overline{AJ_{i-1}}$ as the composition of $AJ_{i-1}$
with a projection and using the end of $\S5.1$. For the projection,
start with\[
Ext_{_{MHS}}^{1}\left(\QQ(-p),\underline{H^{i-1}}(\gs)\otimes\frac{H^{2p-m-i}(X)}{F_{h}^{p-i+1}H^{2p-m-i}(X)}\right)\cong\]
\[
\frac{\underline{H^{i-1}}(\gs,\CC)\otimes\frac{H^{2p-m-i}(X,\CC)}{F_{h}^{p-i+1}H^{2p-m-i}(X,\CC)}}{F^{p}\{\text{{num}}\}+H^{i-1}(\gs,\QQ)\otimes\frac{H^{2p-m-i}(X,\QQ(p))}{F_{h}^{p-i+1}H^{2p-m-i}(X,\QQ(p))}}.\]
Since $F_{h}^{p-i+1}H^{2p-m-i}(X,\CC)\subseteq F^{p-i+1}H^{2p-m-i}(X,\CC),$
and$\linebreak$ $F^{p}\left(H^{i-1}(\gs,\CC)\otimes\frac{H^{2p-m-i}(X,\CC)}{F^{p-i+1}H^{2p-m-i}(X,\CC)}\right)=0$,
this projects to\[
\frac{\underline{H^{i-1}}(\gs,\CC)\otimes\frac{H^{2p-m-i}(X,\CC)}{F^{p-i+1}H^{2p-m-i}(X,\CC)}}{im\left\{ \underline{H^{i-1}}(\gs,\QQ)\otimes\frac{H^{2p-m-i}(X,\QQ(p))}{F_{h}^{p-i+1}}\right\} }\cong\underline{H^{i-1}}(\gs,\QQ)\otimes V^{\vee}/\Lambda\]
 and thus to the r.h.s. of (7.1).

In the incomplete spread case (7.3), the result is considerably deeper
and involves a study of meromorphic $V^{\vee}/\Lambda$-valued differential
characters (sufficiently parallel to $\S\S10,11$ of \cite{Ke1} that
we won't give details).
\end{proof}
There are three special cases or extensions. First, there is a version
of (7.1) for indecomposables, namely (for $1\leq j\leq m\leq i\leq p$)\small \[
\L^{i}\underline{CH}_{j\text{{-ind}}}^{p}(X_{K},m)\rTo^{\overline{AJ}_{i-1,\text{{ind}}}}Hom_{\QQ}\left(\overline{H_{i-1}}(\gs,\QQ),\frac{\{ F^{d-p+i+1}H^{2d-2p+m+i}(X,\CC\}^{\vee}}{im\{ H_{2d-2p+m+i}(X,\QQ(p))\}}\right),\]
 \normalsize again computed by $(2\pi i)^{p-m}R_{\sZ}(\gamma)'$
(only viewed as a functional {[}mod periods{]} on a smaller space).
For example, for $\L^{2}\underline{CH}_{\text{{ind}}}^{3}(X_{K},1)$
with $X$ a $3$-fold, the target is $Hom\left(\overline{H_{1}}(\gs),\frac{\{\Omega^{3}(X)\}^{\vee}}{im\{ H_{3}(X,\QQ(3))\}}\right)$.

Second, there is the heretofore neglected case $m=0$ of algebraic
cycles. Although these have complete spread $\bar{\sZ}$, it is not
possible in general to arrange $\bar{\sZ}\homeq0$; so this is more
like the situation of (7.2),(7.3) than (7.1). One has a map\small \[
\{ ker(cl_{i})\subseteq\L^{i}CH^{p}(X_{K})\}\rTo^{\overline{AJ_{i-1}}}Hom_{\QQ}\left(\overline{H_{i-1}}(\gs,\QQ),\frac{\{ F^{d-p+i}H^{2d-2p+i}(X,\CC)\}^{\vee}}{\Lambda}\right),\]
 \normalsize with $\overline{AJ_{i-1}}\left\langle \Z\right\rangle [\gamma]$
computed by $(2\pi i)^{p}\int_{\partial^{-1}\sZ(\gamma)}(\,\cdot\,)$
since $R_{\sZ}=0$ and $T_{\sZ}=\sZ$. (Here $\sZ(\gamma):=\pi_{X}\{\sZ\cap(\gamma\times X)\}$.)

Finally, there is the case of $X=\text{{Spec}}K$, a point. In order
for $V$ to be $\CC$ and not $\{0\}$, we must have $2d-2p+m+i=0$
(with $d=0$), hence $i=2p-m$. Moreover, $d-p+i\leq0$ (again with
$d=0$) $\implies$ $p\leq m$. (7.2) maps $Gr_{\L}^{2p-m}CH^{p}(\text{{Spec}}K,m)$
to (dropping the weight) $F^{2p-m}H^{2p-m}(\gs,\CC)$, and is given
by $\Z(\mapsto\sZ)\mapsto\Omega_{\sZ}$. This is zero if $trdeg(K/k)=2p-m-1$,
in which case (7.3) reduces to\[
Gr_{\L}^{2p-m}CH^{p}(\text{{Spec}}K,m)\to Hom_{\QQ}\left(H_{2p-m-1}(\gs,\QQ),\CC/\QQ(p)\right)\cong J^{p,m}(\gs),\]
represented by $(2\pi i)^{p-m}R_{\sZ}'$ (a $d$-closed $(2p-m-1)$-current
on some $U\subset\ms$). This is just the \cite{KLM} formula for
$AJ$ on $\gs.$

\section{\textbf{Abel-Jacobi for motivic cohomology of \\ relative and singular
varieties}}

In this section we essentially reverse the Gysin spectral sequences
of $\S\S3,5$. This allows us to study, for example, $AJ$ maps on
$H_{\M}^{2p-n}(\Y,\QQ(p))$ rather than $H_{\M,\Y}^{2p-n+2}(\X,\QQ(p+1))\cong CH^{p}(\Y,n)$
when $\Y\subset\X$ is a complete NCD. This is technically more involved,
since pulling back (cycles, currents) is considerably harder than
pushing forward. The reward is in the arithmetically fascinating examples
of $\S8.4$, as well as \cite{DK} and other forthcoming papers.

\subsection{Currents with pullback}

We follow \cite[sec. 4]{Ki} for the first part of this subsection.

Let $X$ be a smooth projective analytic variety over $\CC$. Recall
that a current $T\in\D_{X}^{r}(U)$ on an open subset of $X$ is a
continuous linear functional on $C^{\infty}$ forms compactly supported
on $U$.

\begin{defn}
\emph{Normal currents} on open $U\subset X$ are defined by \small 

\[
\N_{X}^{r}(U):=\left\{ T\in\D_{X}^{r}(U)\left|\begin{array}{c}
T,\, dT\text{{ admit extensions to}}\\
\text{{continuous linear functionals on}}\\
continuous\text{{ forms compactly supported on }}U\end{array}\right.\right\} ;\]
\normalsize $\N_{X}^{r}$ is their sheafification.
\end{defn}
Now let $V\subset X$ be a smooth closed subvariety of codimension
$q$, and set\[
\Psi:=\text{{ characteristic function of }}X\m V.\]
Call $T\in\N_{X}^{r}(U)$ of $V$-$residue$ $type$ if $T=\Psi T$,
and $V$-$transversal$ if also $dT=\Psi dT$. We first treat the
case of $q=1$:

\begin{defn}
Assuming $V\cap U\subset U$ defined by $z=0$,\[
\N_{X}^{r}\{ V\}(U):=\left\{ T\in\N_{X}^{r}(U)\left|\begin{array}{c}
T\text{{ is }}V\text{{-transversal, and }}\exists\\
S,\mathfrak{{S}}\in\N_{X}^{r+1}(U)\text{{ of }}V\text{{-residue type}}\\
\text{{such that }}T\wedge dz=zS,\,\, dT\wedge dz=z\mathfrak{{S}}\end{array}\right.\right\} \]
are the \emph{normal currents of $V$-intersection type} on $U$,
whose sheafification yields $\N_{X}^{r}\{ V\}$.
\end{defn}
Intuitively, we require the {}``singularities'' of our currents
to properly intersect $V$. 

\begin{rem}
Note that $S,\,\mathfrak{{S}}$ are unique if they exist; we think
of $S$ as {}``$T\wedge\frac{dz}{z}$''. By definition $S$ operates
on continuous forms; hence $T$ operates on continuous $\log(V)$
forms, $a$ $fortiori$ on $C^{\infty}$ $\log(V)$ forms.
\end{rem}
The map \[
\iota_{V}^{^{*}}:\,\N_{X}^{r}\{ V\}\to\N_{V}^{r}\]
 is given by sending $T\mapsto(1-\Psi)dS=:\iota_{V}^{^{*}}(T).$ One
may think of $\iota_{V}^{^{*}}(T)$ as {}``$Res_{V}(T\wedge\frac{dz}{z})$'';
showing this is indeed (the push-forward of) a current on $V$ uses
normality of $T$.

Next, for arbitrary codimension $q$, let $\sigma:\,\tilde{X}\to X$
be the blowup along $V$, with $\sigma^{-1}(V)=:\tilde{V}$ and $\sigma^{-1}(U)=:\tilde{U}$.
(Note that the pullback of currents is defined along the projective-bundle
projection $\sigma|_{\tilde{V}}$.)

\begin{defn}
More generally, $\N_{X}^{r}\{ V\}(U):=$ \[
\sigma_{_{*}}\left\{ \tilde{T}\in\N_{\tilde{X}}^{r}\{\tilde{V}\}(\tilde{U})\left|\iota_{\tilde{V}}^{^{*}}(T)=(\sigma|_{\tilde{V}})^{^{*}}\tau\right.\text{{\, for some }}\tau\in\N_{V}^{r}(U\cap V)\right\} .\]
 
\end{defn}
Finally, we consider the case where $V=\bigcup_{i=1}^{N}V_{i}$ is
a normal crossing divisor (NCD) on $X$. (Here we depart from \cite{Ki}.)
Writing $I,J,K\subseteq\{1,\ldots,N\}$ for multi-indices, each $V_{I}:=\bigcup_{i\in I}V_{i}$
must be smooth of codimension $|I|$. Note that if $T$ is $V_{i}$-transversal
($\forall i$), then $T$ is $V_{I}$-transversal ($\forall I$).

Take $U\subset X$ sufficiently small that $\vspace{2mm}$\\
(a) $\exists$ $K=\{ k_{1},\ldots,k_{m}\}$, $m\leq d:=dim(X)$, such
that only $V_{k_{1}},\ldots,V_{k_{m}}$ intersect $U$; and$\vspace{1mm}$\\
(b) $\exists$ analytic chart $\{ z_{k_{1}},\ldots,z_{k_{m}};\underline{w}\}$
on $U$ such that $U\ni\{\underline{0}\}$, $U\cap V_{k}=\{ z_{k}=0\}$
and $U\cap V=\{ z_{k_{1}}\cdots z_{k_{m}}=0\}$.$\vspace{2mm}$\\
We write locally $dz_{I}:=\bigwedge_{i\in I\cap K}dz_{i}$, $z_{I}:=\prod_{i\in I\cap K}z_{i}$.

\begin{defn}
For $V=$NCD, $T\in\N_{X}^{r}\{ V\}(U)$ $\Longleftrightarrow$ the
following hold:\\
(a) $T\in\N_{X}^{r}(U)$ is $V_{j}$-transversal ($\forall j$)\\
(b) $\forall I$, $\exists$ $S_{I},\mathfrak{{S}}_{I}\in\N_{X}^{r+|I|}(U)$

(i) satisfying $T\wedge dz_{I}=z_{I}S_{I}$, $dT\wedge dz_{I}=z_{I}\mathfrak{{S}}_{I}$;
and

${}$ (ii) of $V_{j}$-residue type $\forall j\in I$,

${}$ (iii) $V_{j}$-transversal $\forall j\notin I$,

${}$ (iv) $V_{J}$-transversal $\forall J$ with $|J|>1$.
\end{defn}
\begin{rem}
Again, if such $S_{I},\,\mathfrak{{S}}_{I}$ exist then they are unique;
the sheaves $\N_{X}^{r}\{ V\}$ are called (once more) normal currents
of $V$-intersection type.
\end{rem}
To describe the pullbacks, write $X=V_{\emptyset}$, and put\[
V^{I}:=V_{I}\bigcap\cup_{j\notin I}V_{j}\subseteq V_{I}\]
(e.g., $V^{\emptyset}=V$ itself). Then $\N_{V_{I}}^{r}\{ V^{I}\}$
has pullback maps to any $\N_{V_{J}}^{r}\{ V^{J}\}$ for $J\supseteq I$,
and pullbacks commute. More precisely, the assignment $T\mapsto(1-\Psi_{j})dS_{j}=:\iota_{V_{j}}^{^{*}}(T)$
generalizes for $J=\{ j_{1},\ldots,j_{\ell}\}$ to $\iota_{V_{J}}^{^{*}}:\N_{X}^{r}\{ V\}\to\N_{V_{J}}^{r}\{ V^{J}\}$
by taking $\iota_{V_{J}}^{^{*}}:=\iota_{V_{J}}^{^{*}}\circ\iota_{V_{J\m\{ j_{\ell}\}}}^{^{*}}\circ\cdots\circ\iota_{V_{\{ j_{1},j_{2}\}}}^{^{*}}\circ\iota_{V_{j_{1}}}^{^{*}}.$
That this is independent of the ordering of the elements of $J$ follows
from the existence of the $S_{J'}$ for all $J'\subseteq J$.

Here is the Poincar\'e lemma for intersection currents.

\begin{lem}
The inclusions\[
\Omega_{X}^{\bb}\subset\Omega_{X^{\infty}}^{\bb}\subset\N_{X}^{\bb}\{ V\}\subset\N_{X}^{\bb}\subset\D_{X}^{\bb}\]
are $F^{\bb}$-filtered quasi-isomorphisms in the cases considered
above (Def. 8.2, 8.4, 8.5).
\end{lem}
The main idea here is to show that King's Cauchy slice-kernels $K_{i}$
(which furnish the coboundary between $\cap$-currents and holomorphic
forms) preserve $F^{p}\N_{X}^{\bb}\{ V\}$ at the stalk level. Roughly
speaking, given $S$, $\mathfrak{{S}}$ for $T$, one must produce
normal currents $S'$, $\mathfrak{{S}}'$ for $K_{i}(\text{{cutoff function}}\cdot T)$.
For the non-NCD cases, this is done in \cite{Ki}; for $V=$NCD the
formulas for $S'$, $\mathfrak{{S}}'$ are more complicated (we don't
give them here), but the other details of the proof are the same.

Now assume $V\subset X$ is a NCD (with both $X$ and $V$ complete),
and write $C_{top}^{r}\{ V\}(X)$ for real-codimension-$r$ $C^{\infty}$-chains
properly intersecting all the $V_{I}$. Set\[
C_{\D}^{\bb}(X,\QQ(p))_{V}:=Cone\left\{ \begin{array}{c}
C_{top}^{\bb}\{ V\}(X,\QQ(p))\\
\oplus F^{p}\N_{X}^{\bb}\{ V\}(X)\end{array}\longrightarrow\N_{X}^{\bb}\{ V\}(X)\right\} [-1].\]
It is an obvious corollary of Lemma 8.7 that this computes Deligne
cohomology of $X$; that is,\[
C_{\D}^{\bb}(X,\QQ(p))_{V}\rInto^{\simeq}C_{\D}^{\bb}(X,\QQ(p)).\]
Now define a double complex $\mathfrak{{K}}_{V}^{\bb,\bb}(p)$ by
\[
\mathfrak{{K}}_{V}^{\ell,m}(p):=\oplus_{|I|=\ell}C_{\D}^{2p+m}(V_{I},\QQ(p))_{V^{I}}\]
 with differentials $D:\,\mathfrak{{K}}_{V}^{\ell,m}(p)\to\mathfrak{{K}}_{V}^{\ell,m+1}(p)$
(Deligne homology differential) and

\begin{equation} \sum_{|I|=\ell} \sum_{i\notin I} (-1)^{<i>_I} (\iota_{V_{I\cup i}\subset V_I})^{^*} = \mathfrak{I} : \, \mathfrak{K}^{\ell,m}_V (p) \to \mathfrak{K}^{\ell+1,m}_V(p), \end{equation} 
where $<i>_{I}$ means the position in which $i$ appears in $\{1,\ldots,N\}\setminus I$.
Write $s^{\bb}$ for associated simple (total) complex, and $\DD$
for total differential.

\begin{prop}
Cohomology of\[
C_{\D}^{2p+\bb}(V,\QQ(p)):=s^{\bb}(\mathfrak{{K}}_{V}^{\bb\geq1,\bb}(p))[1]\]
\[
\text{{resp.\,\,\,}}C_{\D}^{2p+\bb}((X,V),\QQ(p)):=s^{\bb}(\mathfrak{{K}}_{V}^{\bb\geq0,\bb}(p))\]
computes Deligne cohomology of $V$ resp. $(X,V)$.
\end{prop}
\begin{rem}
We emphasize the distinction (quite significant for singular $V$)
between Deligne cohomology $H_{\D}^{*}(V,\QQ(p))$ of $V$, and $H_{\D,V}^{*}(X,\QQ(p))$
--- which is Deligne homology of $V$.
\end{rem}
Noting that $\mathfrak{{K}}_{V}^{0,\bb}(p)=C_{\D}^{2p+\bb}(X,\QQ(p))_{V}$,
we see that\small \[
C_{\D}^{2p+\bb}((X,V),\QQ(p))\rTo^{\cong}Cone\left\{ C_{\D}^{2p+\bb}(X,\QQ(p))_{V}\to C_{\D}^{2p+\bb}(V,\QQ(p))\right\} [-1];\]
\normalsize this gives rise to a long-exact sequence\small 

\begin{equation} \to H^{q-1}_{\D}(V,\QQ(p)) \to H^q_{\D}((X,V),\QQ(p)) \to H^q_{\D}(X,\QQ(p)) \to H^q_{\D}(V,\QQ(p)) \to . \end{equation} \normalsize For
$Y:=X,\, V,\,$or $(X,V)$ one has as usual $H_{\D}^{q}(Y,\QQ(p))\cong H_{\H}^{q}(Y,\QQ(p))$
whenever $W_{2p}H^{q}(Y,\QQ)=H^{q}(Y,\QQ)$, and this is the case
here when $q\leq2p$.

\subsection{Higher Chow chains with pullback}

Let $X$ be a nonsingular projective algebraic variety over a field
$K$, with $\QQ\subseteq K\subseteq\CC$. For technical reasons, we
assume $K$ has a finitely generated field extension $L\subseteq\CC$
invariant under complex conjugation.

Given $I\subseteq\{1,\ldots,n\}$, $f:I\to\{0,\infty\}$, $j\in\hat{I}\cup\{0\}$
(where $\hat{I}:=\{1,\ldots,n\}\m I$), define subsets \[
\d_{f}^{I}\bxn:=\cap_{\ell\in I}\{ z_{\ell}=f(\ell)\},\]
\[
\d_{f,\RR}^{I,j}\bxn:=\left\{ \begin{array}{cc}
[\cap_{\ell\in\{1,\ldots,j\}\cap\hat{I}}T_{z_{\ell}}]\cap\d_{f}^{I}\bxn & \text{{if }}j\geq1\\
\\\d_{f}^{I}\bxn & \text{{if }}j=0\end{array}\right.\]
of $\bxn$. Write $\iota_{f}^{I}:\bx^{n-|I|}\rInto\bxn$, $\pi_{I}:\bxn\rOnto\bx^{n-|I|}$
for the obvious inclusions and projections.

\begin{defn}
$Z_{\RR}^{p}(X,n)$ is the free abelian group $(\otimes\QQ)$ generated
by irreducible $\sZ\in Z^{p}((X\times\bxn)_{K})$ with $\sZ_{\CC}$
meeting all $X_{\CC}\times\d_{f,\RR}^{I,j}\bxn$ properly as real
analytic varieties, modulo degenerate cycles. Write $\d_{f}^{I}\sZ:=(\iota_{f}^{I})^{^{*}}\sZ\in Z_{\RR}^{p}(X,n-|I|)$,
$\db\sZ:=\sum_{j=1}^{n}(-1)^{j}(\d_{o}^{j}-\d_{\infty}^{j})\sZ\in Z_{\RR}^{p}(X,n-1)$.
\end{defn}
\begin{rem}
Recall $\sZ\in Z^{p}(X,n)$ if the intersections $\sZ\cap(X\times\d_{f}^{I}\bxn)$
are proper. The following description is useful for checking whether
a higher Chow chain is in good real position: given irreducible $\sZ\in Z^{p}(X,n)$,
$\sZ$ belongs to $Z_{\RR}^{p}(X,n)$ iff (i) $\sZ_{\CC}$ intersects
all $X_{\CC}\times(T_{z_{1}}\cap\cdots\cap T_{z_{j}})$ properly $/\RR$
and (ii) each $\d_{f}^{I}\sZ_{\CC}$ has the same property (in $X\times\bx^{n-|I|}$).
\end{rem}
We describe $V$-intersection chains first in the case where $V\subset X$
is a nonsingular $K$-subvariety of codimension $q$.

\begin{defn}
$Z_{\RR}^{p}(X,n)_{V}$ has generators as in Def. 8.10, with the additional
requirement that $\sZ_{\CC}$ properly intersect all $V\times\d_{f,\RR}^{I,j}\bxn$
properly $/\RR$. (For $Z^{p}(X,n)_{V}$ we just demand proper intersections
with $V\times\d_{f}^{I}\bxn$ in addition to $X\times\d_{f}^{I}\bxn$.)
\end{defn}
\begin{rem}
One can also define the {}``normalized'' subcomplex $C^{p}(X,\bb)_{V}\subset Z^{p}(X,\bb)_{V}$
consisting of cycles with all facet intersections $zero$ except for
that with $X\times\d_{0}^{n}\bxn$; this subcomplex is well-known
to be quasi-isomorphic. It is convenient to write $\d^{-}\bxn:=\left(\cup_{i=1}^{n-1}\d_{0}^{i}\bxn\right)\bigcup\left(\cup_{i=1}^{n}\d_{\infty}^{i}\bxn\right).$
Put $C_{\RR}^{p}(X,n)_{V}:=C^{p}(X,n)\cap Z_{\RR}^{p}(X,n)_{V}$. 
\end{rem}
\begin{lem}
\emph{(}\emph{\noun{Moving Lemma I}}\emph{)} The inclusions \\
\xymatrix{& & & Z^p_{\RR}(X,\bb)_V \ar @{^(->} [r] \ar @{^(->} [d] & Z^p(X,\bb)_V \ar @{^(->} [d] \\ & & & Z^p_{\RR}(X,\bb) \ar @{^(->} [r] & Z^p(X,\bb) \\ }
\\
are all quasi-isomorphisms.
\end{lem}
For the right vertical arrow, this is due to \cite[Cor 3.2]{Lv1}.%
\footnote{\cite{Lv1} proves this for the normalized complexes; use Rem. 8.13
to transfer the result.%
} It will therefore suffice to examine the top horizontal inclusion,
which we do in an appendix at the end of the subsection. Easy homological
algebra shows one can tack a subscript {}``$\RR$'' onto the moving
lemmas of Bloch \cite{Bl1} and Levine:

\begin{cor}
We have $\frac{Z^{p}(X,\bb)}{Z_{\RR}^{p}(X,\bb)}\simeq0$, $Z_{\RR}^{p}(X,\bb)_{V}\rInto^{\simeq}Z_{\RR}^{p}(X,\bb)$
and $\frac{Z_{\RR}^{p}(X,\bb)}{Z_{\RR}^{p-q}(D,\bb)}\rTo^{\simeq}Z_{\RR}^{p}(X\m D,\bb)$
(for $D\subset X$ a closed subvariety of pure codimension $q$).
\end{cor}
The old moving lemmas can also be (with some more algebra) combined
in the following way, which we will require for an application:

\begin{lem}
\emph{\noun{(Moving Lemma II)}} Given $D\subset X$ a NCD,\[
\frac{Z^{p}(X,\bb)_{V}}{Z^{p-1}(D,\bb)_{V\cap D}}\rTo^{\simeq}Z^{p}(X\m D,\bb)_{V\m V\cap D}\]
provided $V\cap D$ is proper and avoids the codimension-$2$ skeleton
of $D$. (One may also replace $Z$ by $Z_{\RR}$ here.)
\end{lem}
Next we want to augment $V$ to a collection $\ms:=\{ S_{1},\ldots,S_{n}\}$
of closed (possibly singular) subvarieties of $X$. By imposing the
requirements of Def. 8.12 for each $S_{i}$ one obtains $Z_{(\RR)}^{p}(X,\bb)_{\ms}$.
All the moving lemmas still hold with $\ms$ replacing $V$.

Given $V\subset X$ a NCD, write $\V^{I}$ for the collection of all
$V_{J}$ with $J\supsetneq I$ (see $\S8.1$). We take $Z_{V}^{\ell,m}(p):=\oplus_{|I|=\ell}Z^{p}(V_{I},-m)_{\V^{I}}$
with differentials $\db:Z_{V}^{\ell,m}(p)\to Z_{V}^{\ell,m+1}(p)$
and $\mathfrak{{I}}:Z_{V}^{\ell,m}(p)\to Z_{V}^{\ell+1,m}(p)$ (essentially
identical to (8.1) above), and total differential $\DD$. Since there
is already a notion of $Z^{p}(V,-\bb)$ (computing $K'$-theory),
we set\[
\widetilde{Z}^{p}((X,V),-\bb):=s^{\bb}(Z_{V}^{\bb\geq0,\bb}(p))\,,\,\,\,\widetilde{Z}^{p}(V,-\bb):=s^{\bb}(Z_{V}^{\bb\geq1,\bb}(p))[1].\]
These have respective $(-n)^{\text{{th}}}$ cohomologies\[
H_{\M}^{2p-n}((X,V),\QQ(p))=:\widetilde{CH}^{p}((X,V),n)\]
and \[
H_{\M}^{2p-n}(V,\QQ(p))=:\widetilde{CH}^{p}(V,n)\]
($\neq CH^{p}(V,n)$ in general). One also defines $Z_{V,\RR}^{\ell,m}(p)$
and $\widetilde{Z}_{\RR}^{p}$'s in the obvious way, with the same
cohomologies (by Moving Lemma I). Comments about cones as in $\S8.1$
yield a version of (8.2) with $\M$ replacing $\D$. Note that negative
$n$ produces nontrivial $\widetilde{CH}$ groups in both cases: e.g.
for $N>m$, $\widetilde{CH}^{p}\left((X\times\bx^{N},X\times\d\bx^{N}),\,-m\right)\cong CH^{p}\left(X,N-m\right)\cong\widetilde{CH}^{p}\left(X\times\d\bx^{N+1},-m\right).$

Finally, the following notion of relative cycles is convenient:\[
RZ^{p}(X,V):=\ker(\mathfrak{{I}})\subseteq\left\{ Z_{V}^{0,0}(p)=Z^{p}(X)_{\V^{\emptyset}}\right\} .\]
We also define \[
RZ^{p}(V):=\ker(\mathfrak{{I}})\subseteq\left\{ Z_{V}^{1,0}(p)=\oplus_{i=1}^{N}Z^{p}(V_{i})_{\V^{i}}\right\} .\]
Obviously $X$, $V$ need not be complete, and $V$ need not have
normal crossings (or live in some $X$) for these to be defined. \\
\\
\textbf{Appendix: proof of Moving Lemma I.} In the following proof,
we will abuse notation by writing $Z^{p}(X)_{V}$ for $Z^{p}(X)_{\V^{\emptyset}}$,
i.e. using a NCD to denote the collection of its own {}``faces'';
the same goes for {}``real analytic NCD's''. So as a subscript,
$\bxn_{\RR}$ will mean the set of all $\d_{f,\RR}^{I,j}\bx^{n}$
(including $|I|,j=0$), $\d\bxn$ the collection of $\d_{f}^{I}\bxn$
($|I|\geq1$), and so on. We write $(\cdot)_{X}$ for $X\times(\cdot)$.

In the square of inclusions\\
\xymatrix{& & & C^p_{\RR}(X,\bb)_V \ar @{^(->} [r]_{\I} \ar @{^(->} [d]_{\J} & Z^p_{\RR}(X,\bb)_V \ar @{^(->} [d]_{\J'} \\ & & & C^p(X,\bb) \ar @{^(->} [r]_{\I'} & Z^p(X,\bb) , \\ }\\
it is well-known that $\I'$ is a quasi-isomorphism. We will show
here that \textbf{(a}) $\J$ is a quasi-isomorphism and \textbf{(b)}
$\I$ is a quasi-surjection%
\footnote{For $\alpha\in Z_{\RR}^{p}(X,n)_{V,\,\db\text{{-cl.}}}$, $\exists$
$\beta\in Z_{\RR}^{p}(X,n+1)_{V}$ and $\gamma\in C_{\RR}^{p}(X,n)_{V,\,\db\text{{-cl.}}}$
s.t. $\alpha=\db\beta+\I(\gamma)$.%
} (in fact, it is a $\simeq$). It follows that $\J'$ is a quasi-isomorphism;
applying this to \\
\xymatrix{& & Z^p_{\RR}(X,\bb)_V \ar @{^(->} [rr] \ar @{^(->} [rd] & &  Z^p_{\RR}(X,\bb) \ar @{^(->} [ld] \\ & & & Z^p(X,\bb) \\ }\\
for $V=V,\emptyset$ shows $Z_{\RR}^{p}(X,\bb)_{V}\rTo^{\simeq}Z_{\RR}^{p}(X,\bb)$
as required. $\vspace{2mm}$\\
\textbf{(a)} We must show \SMALL \[
\frac{RZ^{p}\left(\bxnx,\d\bxnx\right)_{\bxn_{V}\cup\bxn_{\RR,X}}}{\d_{0}^{n+1}\left\{ RZ^{p}\left(\bx_{X}^{n+1},\d^{-}\bx_{X}^{n+1}\right)_{\bxn_{V}\cup\bxn_{\RR,X}}\right\} }\rTo^{\J}\frac{RZ^{p}\left(\bxnx,\d\bxnx\right)}{\d_{0}^{n+1}\left\{ RZ^{p}\left(\bx_{X}^{n+1},\d^{-}\bx_{X}^{n+1}\right)_{\d_{0}^{n+1}\bx_{X}^{n+1}}\right\} }\]
\normalsize is an isomorphism. This makes heavy use of \cite{Lv1}.
Set\[
\II:=\left\{ (I,f)\,|\, I\subseteq\{1,\ldots,n\},\, f:I\to\{0,\infty\}\right\} ,\]
and note $|\II|=3^{n}$; \cite{Lv1} constructs a scheme via iterated
doubles, which as a topological space is\[
\T^{n}:=\II\times\bxnx\left/\left\langle \begin{array}{c}
(I,f,x)\sim(I',f',x')\text{{ if}}\\
I\subset I',\, f=f'|_{I},\, x=x'\in\d_{\overline{f'}}^{I'\m I}\bxnx\end{array}\right\rangle \right..\]
The preimage of a subvariety (or real analytic subspace) $W\subset\bxnx$
under the obvious projection $\T^{n}\rOnto\bxnx$, is denoted $\T^{n}(W)$.
We may combine the inclusions $\jmath_{f}^{I}:\bxnx\rInto\T^{n}$
into\[
\jmath^{*}:=\sum_{I,f}(-1)^{|I|}(\jmath_{f}^{I})^{^{*}}:\, Z^{p}(\T^{n})\to Z^{p}(\bxnx).\]

Given $\xi\in RZ^{p}(\bxnx,\d\bxnx)$, Levine's argument in \cite[Thm. 2.3-4]{Lv1}
(viewed through the lens of cycles, and his 2.7, 3.1) implies $\exists$
of $\alpha\in RZ^{p}(\bx_{X}^{n+1},\d^{-}\bx_{X}^{n+1})_{\d_{0}^{n+1}\bx_{X}^{n+1}}$,
$f:\T^{n}\rInto H=\,$homogeneous space for $GL_{n}(K)$, and $\eta\in Z^{p}(H)_{f(\T^{n})}$
such that the following holds. Writing $\fs$ for the composition\[
Z^{p}(H)_{f(\T^{n})}\rTo^{f^{*}}RZ^{p}(\T^{n})\rTo^{\jmath^{*}}RZ^{p}(\bxnx,\d\bxnx),\]
$\fs\eta=\xi+\d_{0}^{n+1}\alpha$. Moreover, if $g\in GL_{n}(K)$
has $g^{*}\eta\in Z^{p}(H)_{f(\T^{n})}$, then $\fs g^{*}\eta=\xi+\d_{0}^{n+1}\alpha'$.

Our point is that by extending scalars to some $L\supset K$, we can
choose $g$($\in GL_{m}(L)$) so that $g^{*}\eta\in Z^{p}(H_{L})_{f(\T^{n})_{L}\cup f(\T^{n}(\bxn_{V}))_{L}\cup f(\T^{n}(\bxn_{\RR,X}))},$
which implies $\fs g^{*}\eta\in RZ^{p}(\bxn_{X_{L}},\d\bxn_{X_{L}})_{\bxn_{V_{L}}\cup\bxn_{\RR,X}}$.
Namely, choose $L$ such that $[L:L\cap\RR]=2$, and {}``realify''
the entire situation, replacing all $L$-varieties by underlying $L\cap\RR$-varieties
of twice the dimension and $GL_{m}(L)$ by its isomorphic copy in
$GL_{2m}(L\cap\RR)$. (The realification of $H_{L}$ is a homogeneous
space for this isomorphic copy, not for $GL_{2m}(L\cap\RR)$.) Then
the $f(\T^{n}(\bxn_{\RR,X}))$ are now algebraic (our $L\cap\RR$)
and we may apply Kleiman transversality (as quoted in \cite{Lv1})
to obtain $g$ as described.

In fact, if $[L:K]<\infty$ then one can choose $g$ such that \[
\sigma^{^{*}}(g^{^{*}}\eta)=(\sigma^{^{*}}g)^{^{*}}\eta\,\in\, Z^{p}(H_{L})_{f(\T^{n})_{L}\cup f(\T^{n}(\bxn_{V}))_{L}\cup f(\T^{n}(\bxn_{\RR,X}))}\]
$\forall\sigma\in Gal(L/K)$, so that \[
Norm(\fs g^{^{*}}\eta):=\frac{1}{[L:K]}\sum_{\sigma}\fs(\sigma^{^{*}}(g^{^{*}}\eta))\,\in\, RZ^{p}\left(\bxnx,\d\bxnx\right)_{\bxn_{V}\cup\bxn_{\RR,X}},\]
completing the argument for surjectivity of $\J$ in that case. If
$L/K$ is merely finitely generated, then we can combine the norm
argument with specialization of transcendentals.

For injectivity, given $\Xi\in RZ^{p}(\bxnx,\d\bxnx)_{\bxn_{V}\cup\bxn_{\RR,X}}$
such that $\exists$ $\alpha\in RZ^{p}(\bx_{X}^{n+1},\d^{-}\bx_{X}^{n+1})_{\d_{0}^{n+1}\bx_{X}^{n+1}}$
with $\Xi=\d_{0}^{n+1}\alpha$, we construct $\linebreak$$\beta\in RZ^{p}(\bx_{X}^{n+1},\d^{-}\bx_{X}^{n+1})_{\bxn_{V}\cup\bxn_{\RR,X}}$
with $\d_{0}^{n+1}\beta=\Xi$. Write \[
\bx_{X}^{n+1}\m|\Xi|\times\bx\,=:\,{}^{\circ}\bx_{X}^{n+1}\,,\,\,\,\,\,\,\,\,\,^{\circ}\d\bx_{X}^{n+1}\,:=\,\d\bx_{X}^{n+1}\cap{}^{\circ}\bx_{X}^{n+1}\,,\]
\[
{}^{\circ}\T^{n+1}\,:=\,\T^{n+1}\m\T^{n+1}(|\Xi|\times\bx),\]
 etc.

The argument in the proof of \cite[Thm. 2.7]{Lv1} implies $\exists$
of $f:{}^{\circ}\T^{n+1}\rInto H${[}$=$homogeneous space for $GL_{m}(K)${]}
and $\eta\in Z^{p}(H)_{f({}^{\circ}\T^{n+1})}$, such that if $\fs$
denotes the composition\\
\xymatrix{& Z^p(H)_{f({}^{\circ} \T^{n+1})} \ar [r]^{f^*} & RZ^p({}^{\circ} \T^{n+1}) \ar [r]^{\jmath^* \mspace{50mu}} & RZ^p({}^{\circ}\bx^{n+1}_X, {}^{\circ} \d \bx^{n+1}_X ) \ar [d]_{\tiny \begin{matrix} \text{Zariski closure} \\ \text{+ add degenerate} \\ \text{cycle} \end{matrix}} \\ & & & RZ^p(\bx^{n+1}_X,\d^-\bx^{n+1}_X ), \\   }\\
then $\d_{0}^{n+1}(\fs\eta)=\Xi$. Moreover, $\d_{0}^{n+1}\circ\fs$
of any $g^{^{*}}$-translate of $\eta$ also maps to $\Xi$.

Make an extension to $L$ as before, choose a translate of $\eta$
{[}with all $Gal(L/K)$-conjugates{]} properly intersecting $f({}^{\circ}\T^{n+1}(\bx_{V}^{n+1}))$
$/L$ and (in the realification) $f(\T^{n+1}(\bx_{\RR,X}^{n+1}))$,
apply $\fs$ and take norm; this yields $\beta$ as desired. $\vspace{2mm}$\\
\textbf{(b)} We first recall some definitions:\[
d^{p}(X,n):=\sum_{i=1}^{n}im\left\{ \pi_{i}^{^{*}}:\, Z^{p}(\bx_{X}^{n-1})_{\d\bx_{X}^{n-1}}\to Z^{p}(\bxnx)_{\d\bxnx}\right\} \]
\[
d_{\RR}^{p}(X,n):=d^{p}(X,n)\cap Z^{p}(\bxnx)_{\bxn_{\RR,X}}=\sum_{i=1}^{n}\pi_{i}^{^{*}}\left(Z^{p}(\bx_{X}^{n-1})_{\bx_{\RR,X}^{n-1}}\right)\]
\[
Z^{p}(X,n)=\frac{Z^{p}(\bxnx)_{\d\bxnx}}{d^{p}(X,n)}\,\,,\,\,\,\,\,\,\, Z_{\RR}^{p}(X,n)=\frac{Z^{p}(\bxnx)_{\bxn_{\RR,X}}}{d_{\RR}^{p}(X,n)}.\]
Now given%
\footnote{We omit the subscript $V$ throughout, but all steps are still valid
with this reinstated.%
} $\Z\in Z_{\RR}^{p}(X,n)$ with $\db\Z=0$, let $\tilde{\Z}_{1}\in Z^{p}(\bxnx)_{\bxn_{\RR,X}}$
be a representative, and define inductively\[
\tilde{\Z}_{i+1}:=\tilde{\Z}_{i}-\pi_{i}^{^{*}}\d_{\infty}^{i}\tilde{\Z}_{i}\mspace{50mu}\text{{\, for\,\,}}i=1,\ldots n.\]
This produces a repesentative $\sZ:=\tilde{\Z}_{n+1}\in RZ^{p}(\bxnx,\d_{\infty}\bxnx)_{\bxn_{\RR,X}}\subset Z^{p}(\bxnx)_{\bxn_{\RR,X}},$
where $\d_{\infty}\bxnx:=\cup_{j=1}^{n}\d_{\infty}^{j}\bxnx.$

Consider the projections\\
\xymatrix{\bxnx & (x;z_1,\ldots ,z_{i-1},z_i,z_{i+1},z_{i+2},\ldots,z_n) \\ \bx^{n+1}_X \ar [u]^{\pi_{i+2}} \ar [d]_{\pi_{\{i,i+1\}}} & (x;z_1,\ldots,z_{i-1},z_i,z_{i+1},z_{i+1}',z_{i+2},\ldots ,z_n) \ar @{|->} [u] \ar @{|->} [d] \\ \bx^{n-1}_X & (x;z_1,\ldots ,z_{i-1},z_{i+1}',z_{i+2},\ldots ,z_n )\\}\\
and set $Y_{i}:=\{(1-z_{i})(1-z_{i+1})=(1-z_{i+1}')\}\in Z^{1}(\bx_{X}^{n+1}).$
Sending $\W\mapsto(\pi_{i+2})_{_{*}}\left(Y_{i}\cdot\pi_{\{ i,i+1\}}^{^{*}}(\d_{0}^{i}\W)\right)$
yields a map \[
M_{i}\in End\left(RZ^{p}\left(\bxnx\,,\,\cup_{j=1}^{i-1}\d_{0}^{j}\bxnx\bigcup\d_{\infty}\bxnx\right)_{\bxn_{\RR,X}}\right).\]
(It is easy to check that real good position is preserved. For example,
to show $\dim_{\RR}(T_{z_{1}}\cap\cdots\cap T_{z_{j}}\cap M_{i}(\W))\leq2(\dim(X)+n-p)-j\,,$
one uses the facts that $z_{i},z_{i+1}\in\RR^{-}$ $\implies$ $z_{i+1}'=1-(1-z_{i})(1-z_{i+1})\in\RR^{-}$,
and that $\d_{0}^{i}\W$ is in good real position if $\W$ is.) Moreover,\[
\d_{0}^{i}M_{i}(\W)=\d_{0}^{i}\W=\d_{0}^{i+1}M_{i}(\W)\]
and so\[
\W-M_{i}(\W)\in RZ^{p}\left(\bxnx\,,\,\cup_{j=1}^{i}\d_{0}^{j}\bxnx\bigcup\d_{\infty}\bxnx\right)_{\bxn_{\RR,X}}.\]

We claim there is a cycle\[
\H_{i}(\W)\in RZ^{p}\left(\bx_{X}^{n+1}\,,\,\cup_{j=1}^{i-1}\d_{0}^{j}\bx_{X}^{n+1}\bigcup\d_{\infty}\bx_{X}^{n+1}\right)_{\bx_{\RR,X}^{n+1}}\]
with \begin{equation} \db \H_i(\W) = (-1)^{i+1}M_i(\W)-\H_i(\db \W). \end{equation} 
Consider more projections:\\
\xymatrix{\bx^{n+1}_X & (x;z_1,\ldots ,z_{i-1},z_i,z_{i+1},z_{i+1}',z_{i+2},\ldots,z_n) \\ \bx^{n+2}_X \ar [u]^{\pi_{i+3}} \ar [d]_{\pi_{\{i,i+1,i+2\}}} & (x;z_1,\ldots,z_{i-1},z_i,z_{i+1},z_{i+1}',z_{i+1}'',z_{i+2},\ldots ,z_n) \ar @{|->} [u] \ar @{|->} [d] \\ \bx^{n-1}_X & (x;z_1,\ldots ,z_{i-1},z_{i+1}'',z_{i+2},\ldots ,z_n )\\}\\
and set $\Y_{i}:=\{(1-z_{i})(1-z_{i+1})(1-z_{i+1}')=(1-z_{i+1}'')\}\in Z^{1}(\bx_{X}^{n+2});$
then $\H_{i}(\W):=(\pi_{i+3})_{_{*}}\left(\Y_{i}\cdot\pi_{\{ i,i+1,i+2\}}^{^{*}}(\d_{0}^{i}\W)\right)$
does the job.

Now using (8.3) with $\db\W=0$, inductively define $\tilde{\sZ}_{1}:=\sZ$
and\[
\tilde{\sZ}_{i+1}:=\tilde{\sZ}_{i}-M_{i}(\tilde{\sZ}_{i})=\tilde{\sZ}_{i}\pm\db\H_{i}(\tilde{\sZ}_{i})\mspace{50mu}\text{{\, for\,\,}}i=1,\ldots,n-1.\]
Clearly $\tilde{\sZ}_{n}$ differs from $\sZ$ by $\db$ of a cycle
in $Z^{p}(\bx_{X}^{n+1})_{\bx_{\RR,X}^{n+1}},$ and belongs to $RZ^{p}(\bxnx,\d^{-}\bxnx)_{\bxn_{\RR,X}}$;
since $\tilde{\sZ}_{n}$ is $\db$-closed, it actually lives in $RZ^{p}(\bxnx,\d\bxnx)_{\bxn_{\RR,X}}=\ker(\db)\subseteq C_{\RR}^{p}(X,n)$.

\subsection{Abel-Jacobi maps}

An irreducible $\sZ\in Z_{\RR}^{p}(X,n)$ produces {}``KLM currents''
$R_{\sZ},\,\o_{\sZ},\,\delta_{T_{\sZ}}$ \cite[sec. 5.4]{KLM}. For
these to operate on continuous forms, the integrals described there
must be bounded, and this in fact follows from the boundedness against
$C^{\infty}$ forms. Since $d$ of these currents is given in terms
of $R_{\db\sZ},\,\o_{\db\sZ},\,\delta_{T_{\db\sZ}}$, it follows that
they are normal.

Let\[
\mathfrak{{R}}_{\sZ}^{i,j,I,f}:=\overline{\left\{ \begin{array}{c}
\text{{subset of }}X_{\CC}^{an}\text{{ over which }}\sZ_{\CC}\cap(X_{\CC}\times\d_{f,\RR}^{I,j}\bxn)\\
\text{{is equidimensional of real dimension }}i\end{array}\right\} },\]
and consider $V\subset X$ smooth of codimension $1$. To say $\sZ\in Z_{\RR}^{p}(X,n)_{V}$,
is equivalent to stipulating that all $\mathfrak{{R}}_{\sZ}^{i,j,I,f}$
properly intersect $V$. Roughly speaking, this ensures that $R_{\sZ},\,\o_{\sZ},\,\delta_{T_{\sZ}}$
have their most {}``serious'' kinds of singularities or discontinuities
(deltas, log-branch jumps) transversely to $V$. More precisely, the
restriction of the KLM formula $\linebreak$$\sZ\mapsto(-2\pi i)^{p+\bb}\left((2\pi i)^{-\bb}T_{\sZ},\,\o_{\sZ},\, R_{\sZ}\right)$
\cite[sec. 5.5]{KLM} to $V$-intersection higher Chow chains, produces
a map of complexes\[
Z_{\RR}^{p}(X,-\bb)_{V}\rTo^{\widetilde{AJ}}C_{\D}^{2p+\bb}(X,\QQ(p))_{V}\]
 commuting with the pullback. This extends automatically to $codim_{X}(V)>1$
since the pullback of $\sZ\in Z_{\RR}^{p}(X,n)_{V}$ under $V$-blowup
$\sigma:\,\tilde{X}\to X$ belongs to $Z_{\RR}^{p}(\tilde{X},n)_{\tilde{V}}$.

The upshot of all this is the following

\begin{prop}
Let $V\subset X$ be a NCD with smooth components; assume $X,\, V$
complete. The KLM formula induces a map of double complexes\[
Z_{V,\RR}^{\ell,m}(p)\rTo^{\widetilde{AJ}}\mathfrak{{K}}_{V}^{\ell,m}(p),\]
hence maps of total complexes\[
\tilde{Z}_{\RR}^{p}((X,V),-\bb)\to C_{\D}^{2p+\bb}((X,V),\QQ(p))\,,\,\,\,\,\,\tilde{Z}_{\RR}^{p}(V,-\bb)\to C_{\D}^{2p+\bb}(V,\QQ(p))\]
computing $AJ$ on motivic cohomology.
\end{prop}
This yields a map of long-exact sequences\small \\
\xymatrix{{} \ar [r] & \widetilde{CH}^p((X,V),n) \ar [r] \ar [d]^{AJ} & CH^p(X,n) \ar[r] \ar [d]^{AJ} & {\widetilde{CH}^p(V,n)} \ar [r]^{\;\;\;\;\;\;\;\;\; n\to n-1} \ar [d]^{AJ} & {} \\ {} \ar [r] & H^{2p-n}_{\D}((X,V),\QQ(p)) \ar [r] & H_{\D}^{2p-n}(X,\QQ(p)) \ar [r] & H^{2p-n}_{\D}(V,\QQ(p)) \ar [r] & {} \\ }${}$\normalsize 
\\
\\

Much of the foregoing construction extends to more general singular
or relative varieties (over our field $K\subseteq\CC$). Levine \cite[sec. IV.3]{Lv2}
describes a procedure (partly due to Hanamura) for associating to
any reduced finite-type $K$-scheme $V$ a cohomological motive $\ZZ_{V}$
in $D_{mot}^{b}(Sm_{K})$.%
\footnote{triangulated motivic category, defined in \cite[I.2.1.4]{Lv2}.%
} This is a diagram of motives of smooth quasi-projective varieties
(arrows opposite to the actual morphisms), obtained via cubical hyperresolutions;
it can be used to produce $\widetilde{Z}_{(\RR)}^{p}(V,\bb)$ computing
$H_{\M}^{*}(V,\QQ(p))$ ($V$=NCD as done above is essentially a special
case). What we have done for $H_{\D}$ (and $\widetilde{AJ}$) will
also extend in a straightforward way under the following conditions.
For each $V_{\alpha}$ occurring in the hyperresolution, (a) $V_{\alpha}$
is complete and (b) the irreducible components of images of morphisms
(in the hyperresolution) into $V_{\alpha}$ are $locally$ components
of {[}coskeletons of{]} a NCD. Here are the two simplest examples
that go beyond $V=$NCD; another is given in $\S8.4$.

\begin{example}
$V$ a degenerate elliptic curve of Kodaira type $I_{1}$, i.e. a
rational curve with ordinary double point $p$. In $D_{mot}^{b}(Sm_{K})$,

\[
\ZZ_{V}\cong Cone\{\ZZ_{\PP^{1}}\begin{array}[b]{c}
_{\iota_{0}^{^{*}}-\iota_{\infty}^{^{*}}}\\
\longrightarrow\end{array}\ZZ_{\{ p\}}\}[-1].\]
For us, \[
\widetilde{Z}_{\RR}^{p}(V,-\bb):=Cone\{ Z_{\RR}^{p}(\PP^{1},-\bb)_{\{0,\infty\}}\begin{array}[b]{c}
_{\iota_{0}^{^{*}}-\iota_{\infty}^{^{*}}}\\
\longrightarrow\end{array}Z_{\RR}^{p}(\{ p\},-\bb)\}[-1]\]
computes $H_{\M}^{2p+*}(V,\QQ(p))$. If $V$ is embedded in some smooth
projective $X$, then motivic cohomology of $(X,V)$ is computed by
\[
\widetilde{Z}_{\RR}^{p}((X,V),-\bb):=Cone\{ Z_{\RR}^{p}(X,-\bb)_{\{ V,p\}}\to\widetilde{Z}_{\RR}^{p}(V,-\bb)\}[-1].\]

\end{example}
${}$

\begin{example}
$V$ a $K3$ surface with rational singularity at $p$. There exists
a minimal (crepant) resolution $b:\,\tilde{V}\to V$, with $b^{-1}(p)=:D$
an $A$-$D$-$E$ curve \cite[sec. III.2]{BHPV}, hence a NCD. The
weak $2$-resolution \\
\xymatrix{& & & & & D \ar @{^(->} [r]^{\iota_D} \ar [d]^{b|_D} & \tilde{V} \ar [d]^{b} \\ & & & & & {p} \ar @{^(->} [r]^{\iota_p} & V \\ }\\
produces an isomorphism $\ZZ_{V}\cong Cone\{\ZZ_{\tilde{V}}\oplus\ZZ_{\{ p\}}\rightarrow\ZZ_{D}\}[-1]$
in $\linebreak$$D_{mot}^{b}(Sm_{K})$, and\small \[
\widetilde{Z}_{\RR}^{p}(V,-\bb):=Cone\left\{ Z_{\RR}^{p}(\tilde{V},-\bb)_{\D^{\emptyset}}\oplus Z_{\RR}^{p}(\{ p\},-\bb)\rTo^{\iota_{D}^{*}-(b|_{D})^{*}}\widetilde{Z}_{\RR}^{p}(D,-\bb)\right\} [-1]\]
\normalsize computes $H_{\M}$ of $V$.
\end{example}
In both examples, cochains for $H_{\D}$ are constructed similarly,
and $\widetilde{AJ}$ goes through.

\subsection{Regulator period computations on singular varieties}

We start by setting up a pairing between homology and cohomology of
$(X,V)$ (a similar procedure works for $V$). Consider the double
complexes $\D_{V}^{\ell,m}(p):=\oplus_{|I|=\ell}\N_{V_{I}}^{2p+m}\{ V^{I}\}(V_{I})$
with differentials $d,\,\mathfrak{{I}}$ and $\D_{\ell,m}^{V}(-p):=\oplus_{|I|=\ell}\D_{V_{I}}^{2d-2p-2\ell-m}(V_{I})$
with differentials $d$, $Gy:=$$\linebreak$$2\pi i\sum_{|I|=\ell}\sum_{i\in I}(-1)^{<i>_{I\m i}}(\iota_{V_{I}\subset V_{I\m i}})_{_{*}}.$
The associated simple complexes (with resp. total differentials $d\pm\mathfrak{{I}}$,
$d\pm Gy$ both written $\DD$)\[
C_{DR}^{2p+\bb}((X,V),\CC):=s^{\bb}(\D_{V}^{\bb,\bb}(p))\,,\,\,\, C_{2p+\bb}^{DR}((X,V),\CC):=s^{\bb}(D_{\bb,\bb}^{V}(-p))\]
will be our (co)chains; write $Z_{DR},\, Z^{DR}$ for the $\DD$-(co)cycles.
In particular, put $Z_{2p+M}^{DR,\infty}((X,V),\CC):=$\[
\ker(\DD)\cap\left\{ \oplus_{\ell}\left(\oplus_{|I|=\ell}\Omega_{V_{I}^{\infty}}^{2d-2p-\ell-M}\left\langle \log V^{I}\right\rangle (V_{I})\right)\right\} \subset s^{M}(\D_{\bb,\bb}^{V}(-p))\]
for $\DD$-cycles given by $C^{\infty}$ log forms (modulo $\DD$-boundary,
they all are). The entrywise pairings\[
\N_{V_{I}}^{2p+\ell+M}\{ V^{I}\}(V_{I})\otimes\o_{V_{I}^{\infty}}^{2(d-\ell)-(2p+\ell+M)}\left\langle \log V^{I}\right\rangle (V_{I})\rTo\CC\]
(which exist by Rem. 8.3), $summed$ $over$ $the$ $diagonal$ (in
the double complex), yield\[
Z_{DR}^{2p+M}((X,V),\CC)\otimes Z_{2p+M}^{DR,\infty}((X,V),\CC)\rTo\CC,\]
 which descends to \begin{equation}  H_{DR}^{2p+M}((X,V),\CC) \otimes H^{DR}_{2p+M}((X,V),\CC ) \rTo \CC . \end{equation} 

Assume henceforth that $n\geq1$, and also that $n\geq p$ or $p>d\,(=\dim(X))$.
We describe in this case a (bilinear) {}``regulator period'' map
\begin{equation} \widetilde{CH}^p((X,V),n)\otimes H^{DR}_{2p-n-1}((X,V),\QQ) \rTo \CC/\QQ(p) \end{equation} 
in terms of $\xi=\{\xi_{I}\in Z_{\RR}^{p}(V_{I},n+\ell)\}_{\Tiny\begin{array}{c}
\ell,I\\
|I|=\ell\end{array}}$ and $\alpha=\{\alpha_{I}\}_{_{\ell,I}}\in$$\linebreak$$Z_{2p-n-1}^{DR,\infty}((X,V),\CC)$
representing classes in the respective $\otimes$-factors. Since $n\geq1$,
$F^{p}H^{2p-n}((X,V),\CC)\cap H^{2p-n}((X,V),\QQ)=0$; since also
$n\geq p$ or $p>d$, $H_{\D}^{2p-n}((X,V),\QQ(p))\cong\frac{H^{2p-n-1}((X,V),\CC)}{H^{2p-n-1}((X,V),\QQ(p))}.$
So \[
\widetilde{AJ}(\xi)=\left\{ (-2\pi i)^{p-n-\ell}\left((2\pi i)^{n+\ell}T_{\xi_{I}},\o_{\xi_{I}},R_{\xi_{I}}\right)\right\} _{_{\ell,I}}\]
 can be moved (with only rational ambiguity) by $\DD$-coboundary
to \[
\mathcal{{R}}:=\left\{ (-2\pi i)^{p-n-\ell}(0,0,R_{\xi_{I}}')\right\} _{_{\ell,I}}\in im\left\{ Z_{DR}^{2p-n-1}((X,V),\CC)\right\} .\]
 Now (8.5) of $\xi\otimes\alpha$ is given by (8.4) of $\R\otimes\alpha$,
which is \[
(-2\pi i)^{p-n}\sum_{\ell,I}\int_{V_{I}}R_{\xi_{I}}'\wedge\alpha_{I}.\]
 By taking limits we can replace the $\log(V^{I})$ $C^{\infty}$-forms
$\alpha_{I}$ by relative topological cycles $\gamma_{I}\in Z_{2p-n-\ell-1}^{top}((V_{I},V^{I}),\QQ)$.
(Here $\gamma=\{\gamma_{I}\}$ gives a $\DD$-cycle in $C_{2p-n-1}^{DR}((X,V),\CC)$
by taking the associated $\delta$-currents; it is these $\delta$-currents
which are limits of $\{\alpha_{I}^{\epsilon}\}$'s, and by definition
$\int_{\gamma}R:=\lim_{\e\to0}\int R\wedge\alpha^{\e}$.) In fact,
unless $p=n$ and $p\leq d$, the $\o_{\xi_{I}}$ are all $0$. So
the $\DD$-coboundary involves only topological chains, and \begin{equation} (-2\pi i)^{p-n}\int_{\gamma}R_{\xi} := (-2\pi i)^{p-n} \sum_{\ell} \sum_{|I|=\ell} \int_{\gamma_I} R_{\xi_I} \;\; \in \;\; \CC/\QQ(p) \end{equation} 
computes (8.5). We find these periods in two examples, in which $V$
replaces $(X,V)$ and $p>\dim(V)$ suffices (in lieu of $p>\dim(X)$).

\begin{example}
\textbf{($\widetilde{CH}^{2}(\text{{---}},2)${[}$\cong H_{\M}^{2}(\text{{---}},\QQ(2))${]}
of a singular curve.)}\\
Let $\zeta_{1}$, $\zeta_{2}$ be roots of unity (other than $-1$)
with $\zeta_{1}+\zeta_{2}\neq1$. In $\AA^{2}$, let $V_{1}^{\circ}=\{ x=\zeta_{1}\}$,
$V_{2}^{\circ}=\{ y=\zeta_{2}\}$, $V_{3}^{\circ}=\{ x+y=1\}$; their
closures in $\PP^{2}$ are denoted $V_{i}$. Set $V=\cup V_{i}$ and
$V^{\circ}=V_{1}\cup V_{2}\cup V_{3}^{\circ}$; $V$ is defined (and
we work) over $K=\QQ(\zeta_{1},\zeta_{2})$. We first want to construct
a higher Chow cocycle\[
\xi\in\ker(\DD)\subset\widetilde{Z}_{\RR}^{2}(V,2)=\oplus_{i}Z_{\RR}^{2}(V_{i},2)_{\V^{i}}\bigoplus\oplus_{i<j}Z_{\RR}^{2}(V_{ij},3),\]
and then use (8.6) to pair $\widetilde{AJ}(\xi)$ with the $\gamma=\{\gamma_{I}\}$
shown: \\
\includegraphics[%
  bb=-6cm 0bp 195bp 129bp,
  clip,
  scale=0.6]{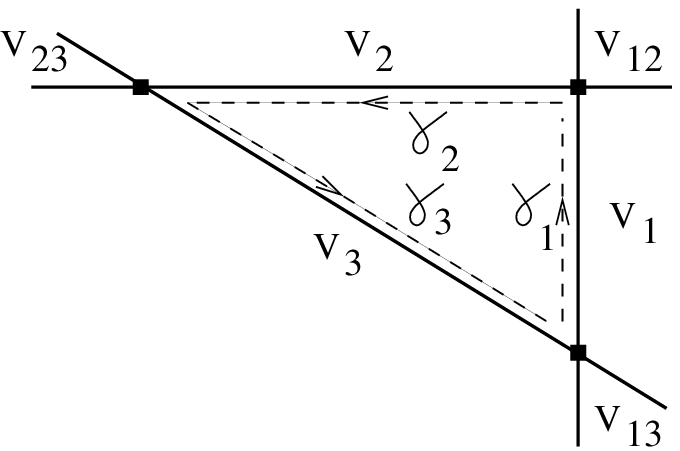}\\
Denote the graph of $(x,y)$ in $Z_{\RR}^{2}(\AA_{K}^{2},2)_{(\V^{\circ})^{\emptyset}}$
by $\{ x,y\}$; the $Tame$ symbols of its restrictions to the $\eta_{V_{i}}$
are torsion ($=0$ for us), and its pullback $\xi_{3}^{\circ}$ to
$V_{3}^{\circ}$ is $\db$-closed. The Zariski closures of $\{ x,y\}|_{V_{3}^{\circ}}=:\xi_{i}^{\circ}$
($i=1,2$) may be {}``completed'' to $\db$-closed $\xi_{i}\in Z_{\RR}^{2}(V_{i},2)_{\V^{i}}$
by adding appropriate $\a_{1},\,\a_{2}$ supported resp. on $\{ y=0,\infty\}\subset V_{1}$,
$\{ x=0,\infty\}\subset V_{2}$. (If $\zeta_{i}=1$, there is nothing
to do.) Hence we have\[
\xi^{\circ}:=(\xi_{1},\,\xi_{2},\,\xi_{3}^{\circ};\,0,\,0,\,0)\in\widetilde{Z}_{\RR}^{2}(V^{\circ},2).\]
Applying Moving Lemma II to $\xi_{3}^{\circ}$ yields a cochain in
$Z_{\RR}^{2}(V_{3},2)_{\V^{3}}$ with restriction $\xi_{3}^{\circ}+\db\beta_{3}$
to $V_{3}^{\circ}$; we may complete it into a $\db$-cocycle $\xi_{3}$
by adding $\alpha_{3}$ supported on $\{(\infty,\infty)\}\in V_{3}$.
Keeping track of $\DD\beta_{3}$ (not just $\db\beta_{3}$) in this
{}``move'', yields the desired\[
\xi:=(\xi_{1},\,\xi_{2},\,\xi_{3},\,\pm\beta_{3}|_{V_{13}},\,\pm\beta_{3}|_{V_{23}},\,0),\]
with regulator period $\int_{\gamma}R_{\xi}=\{\sum_{i=1}^{3}\int_{\gamma_{i}}R_{\xi_{i}}\}\pm\frac{1}{2\pi i}\{ R_{(\beta|_{V_{13}})}+R_{(\beta|_{V_{23}})}\}.$
But general principles (Stokes's theorem for currents) ensure that
$R_{\xi^{\circ}}$ must have (mod $\QQ(2)$) the same period, which
is easy to compute:\[
\int_{\gamma}R_{\xi^{\circ}}\,\,=\,\,\sum_{i=1}^{3}\int_{\gamma_{i}}R\{ x,y\}\]
\[
=\,\,\int_{y=1-\zeta_{1}}^{y=\zeta_{2}}\log\zeta_{1}\dlog y\,\,+\,\,\int_{x=\zeta_{1}}^{x=1-\zeta_{2}}(-2\pi i)\log\zeta_{2}\delta_{T_{x}}\]
\[
\mspace{100mu}+\,\,\int_{y=\zeta_{2}}^{y=1-\zeta_{1}}\{\log(1-y)\dlog y-2\pi i\log y\delta_{T_{1-y}}\}\]
${}$\[
=\,\,-\log\zeta_{1}\log(1-\zeta_{1})\,\,+\,\,0\,\,+\,\,\{\text{{Li}}_{2}(\zeta_{2})-\text{{Li}}_{2}(1-\zeta_{1})\}\]
${}$\[
\equiv\,\,\text{{Li}}_{2}(\zeta_{1})\,+\,\text{{Li}}_{2}(\zeta_{2})\,\,\,\in\,\,\,\CC/\QQ(2).\]
This really boils down to a (lifted) Borel regulator map computation
for an element of $K_{3}^{ind}(K)$, in the sense that one may explicitly
construct $\rho$ so that\\
\xymatrix{& & H^2_{\M}(V,\QQ(2)) \ar [rr]^{\rho}_{\cong} \ar [d]^{AJ} & & CH^2(K,3) \ar [d]^{AJ} & \\ & & H^2_{\D}(V,\QQ(2)) & H^1(V,\CC/\QQ(2)) \ar [l]^{\cong} \ar [r]_{\cong}^{\int_{\gamma}} & \CC/\QQ(2) \\ }\\
commutes.
\end{example}
${}$

\begin{example}
\textbf{($\widetilde{CH}^{3}(\text{{---}},3)$ of a singular surface.)}\\
In $(\PP^{1})^{3}$, let $V_{1}=\{ x=1\}$, $V_{2}=\{ y=1\}$, $V_{3}=\{ z=1\}$,
$V_{4}=\{ z=\frac{(x-1)(y-1)}{xy}\}$, and $V:=\cup V_{i}$. This
is not quite a NCD as $V_{14}$ and $V_{24}$ are reducible. For this
reason we write $V^{4}\,(=V_{14}\cup V_{24}\cup V_{34})=:W=\cup W_{j}$,
with\[
W_{1}=\{ x=1,\, z=0\},\,\, W_{2}=\{ y=1,\, z=0\},\,\, W_{3}=\{ x=1,\, y=0\},\]
\[
W_{4}=\{ x=0,\, y=1\},\text{{ and }}W_{5}=\{ z=1\}\cap V_{4}.\]
Let $D=D_{1}\cup D_{2}$ where $D_{1}:=\{ x=0,\, z=\infty\}\subset V_{4}$,
$D_{2}:=\{ y=0,\, z=\infty\}\subset V_{4}$, and set $V_{4}^{\circ}:=V_{4}\m D$,
$W_{j}^{\circ}:=W_{j}\cap V_{4}^{\circ}$, $W^{\circ}=\cup W_{j}^{\circ}$.
Finally let $C_{1}=\{ x=\infty\}\cap V_{4}$, $C_{2}=\{ y=\infty\}\cap V_{4}$.
We work over $K=\QQ$.

We will compute the regulator period of $\xi\in\ker(\DD)\subset\widetilde{Z}_{\RR}^{3}(V,3)$
over the topological $2$-cycle $\gamma=\{\gamma_{I}\}$ with $\gamma_{1}=\{1\}\times[0,1]\times[0,1]$,
$\gamma_{2}=[0,1]\times\{1\}\times[0,1],$ $\gamma_{3}=(\{ x+y\geq1\}\subset[0,1]^{2})\times\{1\},$
$\gamma_{4}=(\{ x+y\geq1\}\subset[0,1]^{2})\cap V_{4}.$\\
 \\
\includegraphics[%
  clip,
  scale=0.6]{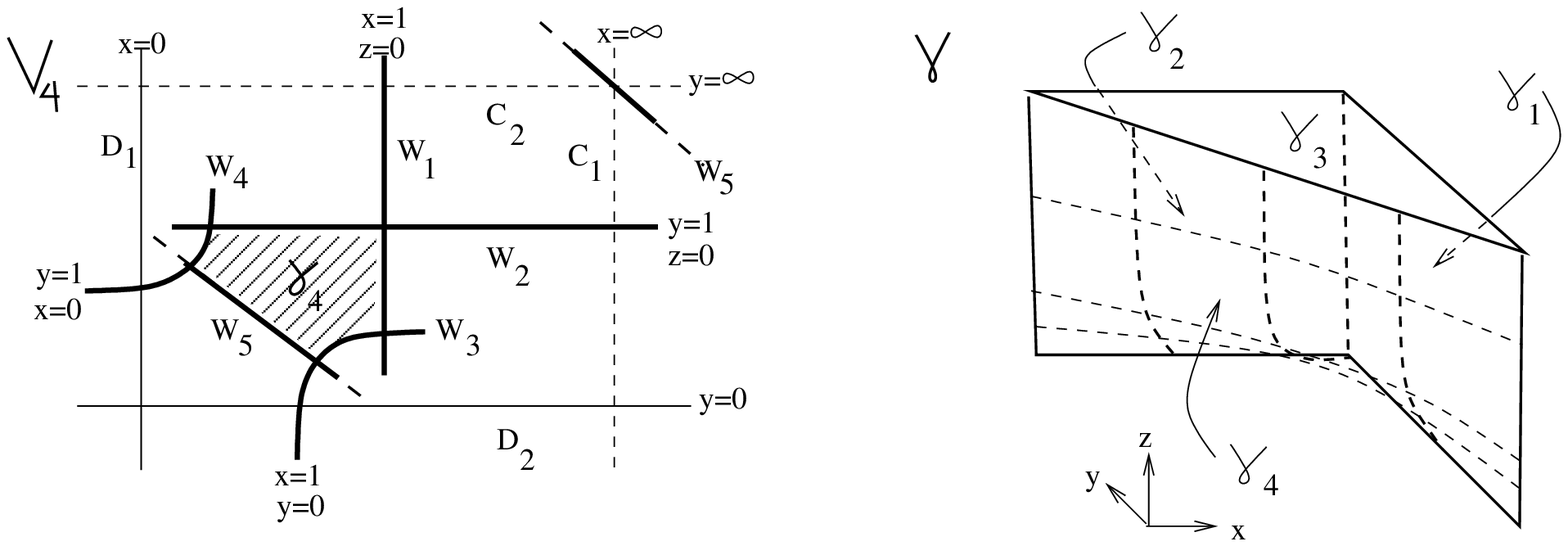}\\
To construct $\xi$, the idea is to move and complete the pullback
of $\{ x,y,z\}\in Z_{\RR}^{3}((\GG_{m})^{3},3)$ to $\widetilde{Z}_{\RR}^{3}(V\cap(\GG_{m})^{3},3)$
(which has trivial $Tame$ symbols on components). Since $\{ x,y,z\}$
is identically zero on $V_{1},\, V_{2},$ and $V_{3}$, we do this
construction in the subcomplex $\widetilde{Z}_{\RR}^{3}((V_{4},W),-\bb)$
of $\widetilde{Z}_{\RR}^{3}(V,-\bb).$

The images of the maps\[
\PP^{1}\m\{0\}\,\times\,\PP^{1}\m\{0,1,\infty\}\,\,\rInto\,\, V_{4}^{\circ}\times\bx^{3}\]
\[
(t,\, u)\,\,\mapsto\,\,(\infty,\, t,\,1-t^{-1};\, u^{-1},\,1-\frac{1}{tu},\,1-u)\]
\[
\mspace{50mu}\mapsto\,\,(\infty,\, t,\,1-t^{-1};\, u^{-1},\,1-\frac{1}{u},\,1-u)\]
\[
\mspace{50mu}\mapsto\,\,(t,\,\infty,\,1-t^{-1};\, u^{-1},\,1-\frac{1}{tu},\,1-u)\]
\[
\mspace{50mu}\mapsto\,\,(t,\,\infty,\,1-t^{-1};\, u^{-1},\,1-\frac{1}{u},\,1-u)\]
yield higher Chow cochains $\A_{1}$, $\A_{2}$ resp. $\A_{3}$, $\A_{4}$
supported on $C_{1}$ resp. $C_{2}$, such that $\xi_{4}^{\circ}:=\{ x,y,z\}|_{V_{4}^{\circ}}-\A_{1}+\A_{2}+\A_{3}-\A_{4}$
is $\db$-closed on $V_{4}^{\circ}$, and\[
\xi_{rel}^{\circ}:=(\xi_{4}^{\circ};\,0,\ldots,0)\in\ker(\DD)\subset\widetilde{Z}_{\RR}^{3}((V_{4}^{\circ},W^{\circ}),\,3).\]
That one (if so inclined) can then $repeatedly$ apply the Moving
Lemma II to produce $\xi_{rel}\in\ker(\DD)\subset\widetilde{Z}_{\RR}^{3}((V_{4},W),\,3)$,
is a consequence of exactness of \begin{equation} \begin{matrix} H^3_{\M}((V_4,W),\QQ(3)) \rTo H^3_{\M}((V_4^{\circ},W^{\circ}),\QQ(3)) \rTo^{Res} \\ \\ \{ H^2_{\M}((D_2\m D_{12},D_2\cap W_3),\QQ(3))  \oplus H^2_{\M}((D_1\m D_{12},D_1 \cap W_4),\QQ(3)) \} \end{matrix} \end{equation} 
at the middle term. We could compute $Res(\xi_{rel}^{\circ}),$ but
this is unnecessary: middle-exactness of\[
H_{\M}^{1}(\{1\},\QQ(2))\to H_{\M}^{2}((\AA_{\QQ}^{1},\{1\}),\,\QQ(2))\to H_{\M}^{2}(\AA_{\QQ}^{1},\QQ(2))\]
and vanishing of $CH^{2}(\AA_{\QQ}^{1},2)\cong K_{2}^{M}(\QQ)_{\QQ}$
and $CH^{2}(Spec(\QQ),3)$, imply the second line of (8.7) is zero.

To compute $\int_{\gamma}R_{\xi}$, note that by construction it is
just $\int_{\gamma_{4}}R_{\xi_{rel}}$, which must be $\int_{\gamma_{4}}R_{\xi_{rel}^{\circ}}=\int_{\gamma_{4}}R\{ x,y,z\}$
by (repeated) application of Stokes. That is, $\int_{\gamma}R_{\xi}=$\[
\int_{\gamma_{4}}\left\{ \log(x)\dlog(y)\wedge\dlog(z)+2\pi i\log(y)\dlog(z)\delta_{T_{x}}+(2\pi i)^{2}\log(z)\delta_{T_{x}\cap T_{y}}\right\} .\]
 Writing $\Delta:=\{ x+y\geq1\}\subset[0,1]^{2}$, the latter is simply\[
\int_{\Delta}\log(x)\dlog(y)\wedge\dlog\{\frac{(x-1)(y-1)}{xy}\}.\]
since $\gamma_{4}\cap T_{x}\subset W_{4}$ (and $\log(y)\dlog(z)$,
$\log(z)\delta_{T_{y}}$ are zero on $W_{4}$). So we get\[
\int_{\Delta}\log(x)\dlog(y)\wedge\dlog\{\frac{(x-1)}{x}\}\]
\[
=\,\,-\int_{x=0}^{1}\log(x)\left(\int_{y=1-x}^{1}\dlog y\right)[\dlog(x-1)-\dlog(x)]\]
\[
=\,\,\int_{x=0}^{1}\log(x)\log(1-x)[\dlog(1-x)-\dlog(x)]\]
\[
=\,-2\int_{u=0}^{1}\log(u)\log(1-u)\dlog(u)\]
\[
=\,\,2\sum_{n>0}\frac{1}{n}\int_{0}^{1}u^{n-1}\log(u)du\,\,=\,\,2\sum_{n>0}\frac{1}{n}\left(-\frac{1}{n^{2}}\right)\,\,=\,\,-2\zeta(3)\]
${}$\\
in $\CC/\QQ(3)$.

This can be interpreted as a Borel regulator computation (for $K_{5}(\QQ)_{\QQ}^{(3)}$)
as follows. The element $\Xi:=\left\{ \frac{z}{z-1},\,\frac{(x-1)(y-1)}{(x-1)(y-1)-xyz},\,\frac{x}{x-1}\right\} \in$$\linebreak$$CH^{3}\left((\PP^{1})^{3}\m V,\,3\right)$
has $Res_{V}(\Xi)\in CH^{2}(V,2)$ with cycle-class in $\linebreak$$\hm(\QQ(-3),\, H_{V}^{4}((\PP^{1})^{3},\QQ))$
given by $\d T_{\Xi}$ --- which is precisely $\gamma$. Hence integration
against $\gamma$ splits $H^{2}(V)$, that is, gives a morphism of
MHS $H^{2}(V,\QQ)\to\QQ(0)$. One may construct a compatible {}``splitting''
of motivic cohomology, in the sense that \\
\xymatrix{&  H^3_{\M}(V,\QQ(3)) \ar @{->>} [rr]^{\rho} \ar [d]^{AJ} & & CH^3(Spec(\QQ),5) \ar [d]^{AJ} \\ &  H^3_{\D}(V,\QQ(3)) & H^2(V,\CC/\QQ(3)) \ar [l]_{\cong} \ar @{->>} [r]^{\int_{\gamma}} & \CC/\QQ(3) }\\
commutes. This phenomenon will be treated in a subsequent work.
\end{example}

${}$\\
Matt Kerr$\mspace{20mu}$$\mathsf{{matkerr@math.uchicago.edu}}$\\
Department of Mathematics, University of Chicago, Chicago, IL 60637,
USA\\
\\
James D. Lewis$\mspace{20mu}$$\mathsf{{lewisjd@gpu.srv.ualberta.ca}}$\\
Department of Mathematics, University of Alberta, Edmonton, Alberta,
T6G 2G1, Canada
\end{document}